\documentclass[a4paper,10pt]{article}
\usepackage[a4paper,left=2cm,right=2cm,top=2.5cm,bottom=2.5cm]{geometry}
\usepackage{amsmath,amssymb}
\usepackage{hyperref}
\usepackage{setspace}
\usepackage[utf8]{inputenc}
\usepackage[english]{babel}
\usepackage{booktabs,multirow}
\usepackage{blkarray}
\usepackage{lscape}
\usepackage[T1]{fontenc}
\usepackage{authblk,color}
\usepackage[english]{babel}
\usepackage{amsfonts, amstext, amsthm, booktabs, dcolumn}
\usepackage[titletoc,toc,title]{appendix}
\usepackage{times}

\definecolor{gr}{rgb}{0.7, 0.0, 0.49}
\definecolor{gr1}{rgb}{0.2, 0.1, 0.4}
\newcommand{\matindex}[1]{\mbox{\scriptsize#1}}

\newtheorem{theorem}{\bf Theorem}[section]

\newtheorem{corollary}{\bf Corollary}[section]
\newtheorem{remark}{\bf Remark}[section]
\newtheorem{remarks}{\bf Remarks}[section]
\newtheorem{lemma}{\bf Lemma}[section]

\newcommand\floor[1]{\left\lfloor#1\right\rfloor}

\begin{document}
\title{On Discrete Gibbs Measure Approximation to Runs}
\author[*]{A. N. Kumar}
\author[**]{N. S. Upadhye}
\affil[ ]{\small Department of Mathematics,}
\affil[ ]{\small Indian Institute of Technology Madras,}
\affil[ ]{\small Chennai-600036, India.}
\affil[*]{\small Email: amit.kumar2703@gmail.com}
\affil[**]{\small Email: neelesh@iitm.ac.in}
\date{}
\maketitle

\begin{abstract}

\noindent
A Stein operator for the runs is derived as a perturbation of an operator for discrete Gibbs measure. Due to this fact, using perturbation technique, the approximation results for runs arising from identical and non-identical Bernoulli trials are derived via Stein's method. The bounds obtained are new and their importance is demonstrated through an interesting application.

\end{abstract}

\noindent
\begin{keywords}
Runs; discrete Gibbs measure; perturbations; probability generating function; Stein operator; Stein's method.
\end{keywords}\\
{\bf MSC 2010 Subject Classifications :} Primary : 62E17, 62E20 ; Secondary :  60C05, 60E05.
\maketitle

\section{Introduction}
Runs and patterns is an important topic in the areas related to probability and statistics, such as reliability theory, meteorology and agriculture, statistical testing and quality control among many others (see Balakrishnan and Koutras \cite{BK}, Kumar and Upadhye \cite{KU2} and Dafnis {\em et al.} \cite{DAP} for more details). The research on this topic started with runs related to success/failure (see Philippou {\em et al.} \cite{PGP} and Philippou and Makri \cite{PM}) and  series of articles later followed in this area, see Aki \cite{AKI}, Aki {\em et al.} \cite{AKH}, Antzoulakos {\em et al.} \cite{ABK}, Antzoulakos and Chadjiconstantinidis \cite{AC}, Balakrishnan and Koutras \cite{BK} and Makri {\em et al.} \cite{MPP} and references therein. Huang and Tsai \cite{HT} and Dafnis {\em et al.} \cite{DAP} extended the study by considering failures and successes together which is known as ($k_1,k_2$)-runs or modified distribution of order $k$. In this paper, we consider the approximation problem related to $(k_1,k_2)$-runs.\\
Let $\eta_1,\eta_2,\dotsc,\eta_n$ be a finite sequence of independent Bernoulli trials with ${\mathbb P}(\eta_i=1)=p_i=1-q_i={\mathbb P}(\eta_i=0)$. Define, $\eta_0=1$,
\begin{equation}
{\mathbb I}_l:=\eta_{l-1}\left[\left(\prod_{i=1}^{k_1}(1-\eta_{l+i-1})\right)\left(\prod_{j=1}^{k_2}\eta_{l+k_1+j-1}\right)\right](1-\eta_{l+k_1+k_2}),~1 \le l \le n-k_1-k_2~{\rm and}~M_{k_1,k_2}^n:=\hspace{-0.05cm}\sum_{l=1}^{n-k_1-k_2}{\mathbb I}_l. \label{three:formulation}
\end{equation} 
Then, $M_{k_1,k_2}^n$ denotes the number of occurrences of exactly $k_1$ consecutive failures followed by exactly $k_2$ consecutive successes in $n$ trials with $k_1, k_2 \ge 1$. We denote $M_{k_1,k_2}^n=S_{k_1,k_2}^n$ if $p_i=p$ (i.e., identical Bernoulli trials). The distributional properties, such as probability generating function (PGF), probability mass function (PMF) and moments, are studied by Dafnis {\em et al.} \cite{DAP} for $S_{k_1,k_2}^n$ and are intractable for $M_{k_1,k_2}^n$, in general.\\
Probabilistic approximations play a crucial role in understanding and comparing the behavior of the distributions. Therefore, it is suggested to approximate the distribution of $M_{k_1,k_2}^n$ to the set of some well-known distributions so that the distribution of $M_{k_1,k_2}^n$ can be characterized using these distributions at the cost of error in approximation. Recently, approximation problems, such as Poisson approximation (Vellaisamy \cite{V}), binomial convoluted Poisson approximation (Upadhye {\em et al.} \cite{UCV}), pseudo-binomial approximation and negative binomial approximation (Kumar and Upadhye \cite{KU1,KU2}), related to runs are studied in the literature.\\
In this paper, we study three approximation problems related to the distribution of $M_{k_1,k_2}^n$, namely, discrete Gibbs measure (DGM) approximation to $S_{k_1,k_2}^n$, $S_{k_1,k_2}^n$ approximation to $M_{k_1,k_2}^n$ and DGM approximation to $M_{k_1,k_2}^n$, using Stein's method. It is important to mention that Stein's method is one of the most powerful tool to study probabilistic approximations, due to its applicability for dependent setups. However, the application of Stein's method is heavily dependent on finding a suitable Stein operator. There have been several approaches for identifying a suitable operator (see Stein \cite{stein}, Barbour and G\"{o}tze \cite{Ba,G} and Diaconis and Zabell \cite{DZ}). Recently, an interesting approach is presented by Upadhye {\em et al.} \cite{UCV}, that uses PGF to obtain a Stein operator. This approach is most suitable for the distributions that arise out of dependent setups. We use PGF approach to derive Stein operator for DGM and the distribution of $S_{k_1,k_2}^n$. Note that DGM is a general class of distributions which include Poisson, binomial, geometric, negative binomial and logarithmic series distribution among many other distributions. To the best of our knowledge, these approximation results are new and not studied in the literature. So it is worth exploring this problem, as it may benefit the readers to understand the properties of the distribution of $M_{k_1,k_2}^n$.\\
This paper is organized as follows. In Section \ref{three:KR}, we discuss some known results related to Stein's method along with some notations. In Section \ref{three:AR}, we present new error bounds for total variation distance between $S_{k_1,k_2}^n$, $M_{k_1,k_2}^n$ and DGM, and also obtain Poisson, pseudo-binomial and negative binomial approximation results as corollary to main results. In Section \ref{three:APP}, we show the importance of approximation results to two-stage start-up demonstration testing. In section \ref{ARR}, we derive some useful results related to the distributional properties of $S_{k_1,k_2}^n$. Finally, in Section \ref{three:pr}, we derive the proofs for the results presented in Section \ref{three:AR}.

\section{Preliminaries and Notations}\label{three:KR}
In this section, we discuss some known results related to Stein's method and also define notations to simplify presentation of the paper.\\
Let $X$ and $Y$ be any two random variables, concentrated on ${\mathbb Z}$ (the set of integers), defined on some probability space. Suppose that $X$ has well-known distribution and the distribution of $Y$ is intractable, then Stein's method can be used to study $X$-approximation to $Y$. The following are the important steps involved in this approximation. 
\begin{enumerate}
\item We obtain a Stein operator (denoted by ${\cal A}_X$ for a random variable $X$) which acts on a large class of functions ${\cal G}_X$ such that
$${\mathbb E} \left[{\cal A}_X g(X)\right]=0, \quad {\rm for}~g \in {\cal G}_X,$$
where ${\cal G}_X = \left\{g:~g \in {\cal G}~{\rm such~that}~g(0)=0~{\rm and} ~g(x)=0, ~{\rm for}~x \notin {Supp(X)}  \right\}$, $\cal G$ be the set of all bounded functions and $Supp(X)$ denotes the support of a random variable $X$.
\item We obtain the solution to Stein equation
\begin{equation}
{\cal A}_X g(m) = f (m) - {\mathbb E}f (X),\quad m \in {\mathbb Z} ~{\rm and}~ f \in {\cal G}.\label{three:seq}
\end{equation}
\item We replace $m$ with a random variable $Y$ in Stein equation and take expectations and supremum to find
$$d_{TV}\left(X,Y\right) := \sup_{f \in {\cal J}}\left|{\mathbb E}f(X) - {\mathbb E}f(Y)\right| = \sup_{f \in {\cal J}}\left|{\mathbb E}\left[{\cal A }_Xg(Y)\right]\right|,$$
where ${\cal J} = \{{\bf 1}(A)|~ A \subseteq{\mathbb Z}\}$ and ${\bf 1}(A)$ is the indicator function of the set $A$.
\end{enumerate} 
For more details, see Barbour \cite{Ba}, Barbour and Chen \cite{BAC}, Barbour {\em et al.} \cite{BCX,BCL,BHJ}, \v{C}ekanavi\v{c}ius \cite{CV}, Chen {\em et al.} \cite{CGS}, Eichelsbacher and Reinert \cite{ER}, Nourdin and Peccati \cite{NP}, Ley {\em et al.} \cite{LRS}, Reinert \cite{R1} and references therein. \\
Throughout this paper, let $X_1$, $X_2$ and $X_3$ have Poisson (with parameter $\lambda$), pseudo-binomial (with parameter $\check{\alpha}$ and $\check{p}$) and negative binomial (with parameter $\hat{\alpha}$ and $\hat{p}$) distribution, respectively, with PMFs
\begin{align}
{\mathbb P}(X_1=m)&= \frac{e^{-\lambda}\lambda^m}{m!}, \quad m=0,1,2\dotsc,\label{three:p}\\
{\mathbb P}(X_2=m)&= \frac{1}{R} {\check{\alpha} \choose m} \check{p}^m \check{q}^{\check{\alpha}-m}, \quad m=0,1,2,\dotsc,\floor{\check{\alpha} },\label{three:pb}\\
{\mathbb P}(X_3= m) &= {\hat{\alpha}+m-1 \choose m} \hat{p}^{\hat{\alpha}} \hat{q}^m, \quad \quad m = 0, 1, \dotsc,\label{three:nb}
\end{align}
where $\lambda, \check{\alpha},\hat{\alpha} > 0$ and $0 < \check{p},~\hat{p}<1$ with $\check{q}=1-\check{p}$, $\hat{q}=1-\hat{p}$, $\floor{\check{\alpha}}$ is the greatest integer part of $\check{\alpha}$ and $R = \sum_{m=0}^{\floor{\check{\alpha}}}{\check{\alpha} \choose m} \check{p}^m \check{q}^{\check{\alpha}-m}$. From (34) of Upadhye {\em et al.} \cite{UCV}, the bounds for the solution to the Stein equation for Poisson, pseudo-binomial and negative binomial distributions, respectively, are given by
\begin{equation}
\|\Delta g\| \le \frac{2 \|f\|}{\max(1,\lambda)},\quad \|\Delta g\| \le \frac{2 \|f\|}{\floor{\check{\alpha}} \check{p}\check{q}}\quad {\rm and} \quad \|\Delta g\| \le \frac{2 \|f\|}{\hat{\alpha}\hat{q}}.\label{three:bound}
\end{equation}
where $\|\Delta g\|:=\sup_{j\in{\mathbb Z}_+}|\Delta g(j)|$ and $\Delta g(j)=g(j+1)-g(j)$.\\
Next, let $Z_1$, $Z_2$ and $Z_3$ be three random variables defined on common probability space and we are interested in studying $Z_2$ approximation to $Z_3$. The distribution and Stein equation for $Z_1$ is well-known and the bounds are available in the literature. Suppose a Stein operator for $Z_2$ can be derived and observed as a perturbation of operator for $Z_1$. Then, the following result suggests the method to obtain error bounds between $Z_2$ and $Z_3$.
\begin{lemma}\label{three:perturbation} [Lemma $3.1$, Upadhye {et al.} \cite{UCV}]
Let $Z_1$ be a random variable with support ${\cal S}$, Stein operator ${\cal A}_{Z_1}$ and $g_0$ be the solution to Stein equation \eqref{three:seq} satisfying 
$$\|\Delta g_0\|\le w_1 \|f\| \min(1,\gamma^{-1}),$$
where $w_1,\gamma>0$. Also, let $Z_2$ be a random variable with Stein operator ${\cal A}_{Z_2}={\cal A}_{Z_1}+U_1$ and $Z_3$ be another random variable such that, for $g\in {\cal G}_{Z_1}\cap {\cal G}_{Z_2}$, 
$$\|U_1g\|\le w_2\|\Delta g\|,\quad |{\mathbb E}{\cal A}_{Z_2}g(Z_3)|\le \varepsilon \|\Delta g\|,$$
where $w_1w_2<\gamma$. Then
$$d_{TV}(Z_2,Z_3)\le \frac{\gamma}{2(\gamma-w_1 w_2)}\left(\varepsilon w_1 \min(1,\gamma^{-1})+2 {\mathbb P}(Z_2\in {\cal S}^c)+2 {\mathbb P}(Z_3\in {\cal S}^c)\right),$$
where ${\cal S}^c$ denote the complement of set ${\cal S}$.
\end{lemma}
\noindent
For more details about these results, we refer the reader to Barbour {\em et al.} \cite{BCX}, \v{C}ekanavi\v{c}ius and Roos \cite{CR}, Upadhye {\em et al.} \cite{UCV} and Vellaisamy {\em et al.} \cite{VUC} and references therein.\\
Next, we say that the distribution of the random variable $Z$ belongs to DGM, a family of discrete distributions, if it has the PMF which can be represented as
\begin{equation}
{\mathbb P}(Z=m):=\Lambda(m)=\frac{1}{\beta} \frac{e^{U(m)}w^m}{m!},\quad m\in{\cal S},\label{three:DGM1}
\end{equation}
where ${\cal S}\equiv Supp(Z)$, $w>0$ is fixed, $U:{\cal S} \to {\mathbb R}$ be a function and $\beta=\sum_{m\in {\cal S}}\frac{e^{U(m)}w^m}{m!}$. \\
Next, we give some examples of well-known distributions that belong to DGM family.
\begin{itemize}
\item[{\bf (O1)}] If $\beta=e^{U(b)}/b! $, $w=1$ and ${\cal S}=\{b\}$, where $b$ is a constant then $Z$ follows degenerate distribution.
\item[{\bf (O2)}] If $\beta=1$, $U(m) = -\lambda$, $w=\lambda$ and ${\cal S}=\{0,1,2,\dotsc\}$ then $Z$ follows Poisson distribution with parameter $\lambda$.
\item[{\bf (O3)}] If $\beta=e^\lambda-1$, $U(m) = 0$, $w=\lambda$ and ${\cal S}=\{1,2,3,\dotsc\}$ then $Z$ follows zero-truncated Poisson distribution with parameter $\lambda$.
\item[{\bf (O4)}] If $\beta=1$, $U(m)=\ln (p m!)$, $w=(1-p)$ and ${\cal S}=\{0,1,2,\dotsc\}$ then $Z$ follows geometric distribution with parameter $p$.
\item[{\bf (O5)}] If $\beta=1$, $U(m) = \ln[n(n-1)\dotsb(n-m+1)]+n \ln q$, $w=p/q$ and ${\cal S}=\{0,1,2,\dotsc,n\}$ then $Z$ follows binomial distribution with parameters $n \in {\mathbb N}=\{1,2,\dotsc\}$ and $p$.
\item[{\bf (O6)}] If $\beta=R$, $U(m) = \ln[n(n-1)\dotsb(n-m+1)]+n \ln q$, $w=p/q$ and ${\cal S}=\left\{0,1,2,\dotsc,\floor{n}\right\}$ then $Z$ follows pseudo-binomial distribution with parameters $n \in (0,\infty)$ and $p$ with $R$ as a normalizing constant.
\item[{\bf (O7)}] If $\beta=1$, $U(m)=\ln[n(n+1)\dotsb(n+m-1)]+n \ln p$, $w=1-p$ and ${\cal S}=\{0,1,2,\dotsc\}$ then $Z$ follows negative binomial distribution with parameters $\alpha$ and $p$. 
\item[{\bf (O8)}] If $\beta=-\ln(1-p)$, $U(m)=\ln(m-1)!$, $w=p$ and ${\cal S}=\{1,2,3,\dotsc\}$ then $Z$ follows logarithmic series distribution with parameter $p$.
\end{itemize}
\noindent
Note also that the representation is not unique, for example, other representation for Poisson distribution is $w=\lambda,~ U(m)=0,~\beta = e^{\lambda}$.\\ 
Eichelsbacher and Reinert \cite{ER} derived the following Stein operator for DGM using generator approach (Barbour and G\"{o}tze \cite{Ba,G}). 
\begin{equation}
{\cal A}_{Z}g(m) = w e^{U(m+1)-U(m)}g(m+1)-mg(m),\quad m\in{\cal S},\label{bbb}
\end{equation}
where ${\cal S}=\{0,1,2,\dotsc,N\}$ and $N$ can take the value infinity. This operator can also be derived using PGF approach as follows:\\
The PGF of $Z$, whenever exists, is given by 
\begin{equation}
G(t)=\sum_{m=0}^{N}\Lambda(m)t^m.\label{three:pgfDGM}
\end{equation} 
Therefore, 
\begin{align}
G^\prime(t) =
\sum_{m=0}^{N}m \Lambda (m)t^{m-1}=
\sum_{m=0}^{N}(m+1)\Lambda(m+1)t^m
=
\sum_{m=0}^{N}w e^{U(m+1)-U(m)}\Lambda(m)t^m.\label{three:prime}
\end{align}
Comparing the coefficients of $t^m$, we get
\begin{equation}
(m+1)\Lambda(m+1)=w e^{U(m+1)-U(m)}\Lambda(m).\label{5:jjj}
\end{equation}
Let $g \in {\cal G}_Z$, then $\displaystyle{\sum_{m=0}^{N}g(m+1)(m+1)\Lambda(m+1)=\sum_{m=0}^{N}g(m+1)w e^{U(m+1)-U(m)}\Lambda(m).}$

\noindent
This implies $\displaystyle{{\mathbb E}\left[{\cal A}_{Z}g(Z)\right]=\sum_{m=0}^{N}[w e^{U(m+1)-U(m)}g(m+1)-mg(m)]\Lambda(m)=0.}$

\noindent
Hence, \eqref{bbb} follows.

\begin{remark}
Note that the relation \eqref{5:jjj} can be computed from \eqref{three:DGM1} and therefore, PGF approach is not more relevant as the Stein operator can be computed directly.
\end{remark}

\noindent
Next, we introduce some notations to improve the readability of the paper. Define $k:=k_1+k_2$,
$$a_1:=1,\quad a_2:=-1,\quad a_3:=qp, \quad d_1=d_3:=n-k-2, \quad d_2:=n-k-1,$$
$$b_1(n):=n+1,\quad b_3(n):=n\quad {\rm and}\quad b_2(n):=\left\{
          \begin{array}{ll}
          -q(k+2) & n = k+1,\\
          n+1-q & n \ge k+2.\\
          \end{array}\right. $$
Also, define
\[c_{n,k}^{(i)} := \left\{
          \begin{array}{ll}
           (n-2k-2)(k+1)+\frac{(k+1)^{n-k+1}}{(k+2)^{n-k-1}}, & i=1;\\
           n-3k+k\left(\frac{k+1}{k+2}\right)^{n-2k}, & i=2;\\
           n-5k+k(nk+6k+4-k^2)\frac{(k+1)^{n-3k-1}}{(k+2)^{n-3k+1}}, & i=3;\\
           n(3k+1)-(11k^2+9k+2)+(2nk+k^2+7k+2)\left(\frac{k+1}{k+2}\right)^{n-2k}, & i=4;\\
           n(2k+3)-6k^2-16k-10+(2n+2k+5)\frac{(k+1)^{n-k}}{(k+2)^{n-k-1}},~~~~~~~~ & i=5,
          \end{array}\right.
\]
$$\delta = 2+\delta^*(k+2),\quad\delta_1=\delta^*\left(5k\hspace{-0.07cm}+\hspace{-0.07cm}6\hspace{-0.07cm}+\hspace{-0.07cm}\frac{3k^2\hspace{-0.07cm}+\hspace{-0.07cm}11k\hspace{-0.07cm}+\hspace{-0.07cm}10}{2}\delta^*\right), \delta_2=\delta^*\left(3+\frac{4k+7}{2}\delta^*\right) ~\text{with}~ \delta^*=1+q+qp,~~~~~~~~~~~~~~~~~$$
\begin{equation}
h_1(n,k,p):=(n-k)(4+\delta_1)+\delta_2 c_{n,k}^{(1)}+\delta c_{n,k}^{(2)}+c_{n,k}^{(3)}+\delta^* c_{n,k}^{(4)}+\frac{{\delta^*}^2}{2}c_{n,k}^{(5)} ~{\rm and}~h_2(n,k,p):=(n-k)\delta+\delta^* c_{n,k}^{(1)}+c_{n,k}^{(2)}.\label{three:h1h2}
\end{equation}

\section{Main Results}\label{three:AR}
In this section, we study approximation problem related to $M_{k_1,k_2}^n$, $S_{k_1,k_2}^n$ and DGM using Stein's method. We also show that the approximation results for Poisson, pseudo-binomial and negative binomial distributions follow as special cases. The proofs of the results are given in Section \ref{three:pr}.

\noindent
Note that, it is not possible to obtain a Stein operator for $S_{k_1,k_2}^n$, as a perturbation of an operator for DGM, in general, as Stein operator for DGM contains the term $U(m)$ which is not known in general. However, if the following condition is satisfied
\begin{equation}
 e^{U(m+1)-U(m)}=a+bm,\label{panj}
\end{equation}
\noindent
then a Stein operator can be derived. Also, observe that \eqref{panj} is similar to a well-known Panjer's recursion (see Panjer and Wang \cite{Panj95} for more details) and a large class of distributions, for example, the cases {\bf (O1)-(O8)} among many others, satisfy \eqref{three:DGM1} and \eqref{panj}.

\noindent
Next, we first compute the mean and variance for the distribution of DGM, $S_{k_1,k_2}^n$ and $M_{k_1,k_2}^n$, which can be used for the choice of parameters. From \eqref{three:pgfDGM}, \eqref{three:prime} and \eqref{panj}, we get
$$G^\prime(t)=\frac{wa }{1-wbt}G(t).$$
Therefore, the mean and variance of DGM ($Z$ as defined in \eqref{three:DGM1}) is given by
\begin{align*}
{\mathbb E}(Z) &= G^\prime(1)=\frac{wa}{1-wb}~~{\rm and}~~\mathrm{Var}(Z)=G^\prime(1)+G^{\prime\prime}(1)-(G^\prime(1))^2=\frac{wa}{(1-wb)^2}.
\end{align*}
Also, from \eqref{three:formulation}, it is clear that
$${\mathbb E}(M_{k_1,k_2}^n)=\sum_{i=1}^{n-k}{\mathbb E}({\mathbb I}_i)~~ {\rm and}~~\mathrm{Var}(M_{k_1,k_2}^n)=\sum_{i=1}^{n-k}{\mathbb E}({\mathbb I}_i)+\sum_{i=1}^{n-k}\sum_{\substack{|j-i|\le k+1\\j\neq i}}{\mathbb E}({\mathbb I}_i{\mathbb I}_{j})-\sum_{i=1}^{n-k}\sum_{|j-i|\le k+1}{\mathbb E}({\mathbb I}_i){\mathbb E}({\mathbb I}_j)$$
\vskip -3ex \noindent
and, for $p_i=p,~i=1,2,\dotsc,n$, let $a(p)=q^{k_1}p^{k_2}$, then
\begin{align*}
{\mathbb E}\left(S_{k_1,k_2}^n\right)&=q[1+(n-k-1)p]a(p)~~{\rm and}~~\mathrm{Var}\left(S_{k_1,k_2}^n\right)=q[1+(n-k-1)p]a(p)-s_{n,k},
\end{align*}
where $s_{n,k}=\left[(n(2k+3)-(3k+5)(k+1))q^2p^2-2(k+1)q^3+(2n-2k+1)q^2-2(n-2k)q\right]a(p)^2$.\\

\noindent
Next, we present the approximation results between DGM and $S_{k_1,k_2}^n$. 
\begin{theorem}\label{three:twth2}
Let $n \ge 5k$, ${\mathbb E}(Z)={\mathbb E}\big(S_{k_1,k_2}^n\big)$ and $\varphi=\mathrm{Var}(Z)-\mathrm{Var}\big(S_{k_1,k_2}^n\big)$. Then
\begin{align}
d_{TV}\left(Z,S_{k_1,k_2}^n\right)&\le\|\Delta g\|  \left\{2(2+qp)a(p)^2\left(|1-wb| h_1(n,k,p)a(p)+|wb|h_2(n,k,p)\right)+|\varphi(1-wb)|\right\},\label{three:twoth}
\end{align}
where $Z$ as defined in \eqref{three:DGM1}, and $h_1(n,k,p)$ and $h_2(n,k,p)$ as defined in \eqref{three:h1h2}.
\end{theorem}

\begin{remarks}
\begin{itemize}
\item[{(i)}] Note that bounds obtained in \eqref{three:twoth} are of constant order and new to the best of our knowledge.
Though the bound is of constant order, due to the presence of the term $\left(\frac{k+1}{k+2}\right)^n$, the bounds decrease with increase in the value of $n$ (see Table \ref{tab:table1}).
\item[(ii)] Observe that, if $\mathrm{Var}(Z)=\mathrm{Var}\big(S_{k_1,k_2}^n\big)$ (i.e., $\varphi=0$) then the bound become sharper, as expected.
\item[(iii)] Note that, if $\mathrm{Var}(Z)=\mathrm{Var}\big(S_{k_1,k_2}^n\big)$ (i.e., $\varphi=0$) then the validity of the bound depends on admissibility of parameters. For example, $s_{n,k} > 0$ (mean larger than variance) and $s_{n,k} < 0$ (mean smaller than variance) pseudo-binomial and negative binomial approximations are valid, respectively. However, negative binomial (for $s_{n,k} > 0$) and pseudo-binomial (for $s_{n,k} < 0$) approximations are not valid as  these conditions yield inadmissible parameters.
\end{itemize}
\end{remarks}
\noindent
Next, we present the results for Poisson, pseudo-binomial and negative binomial distributions and the proofs follow directly from ${\bf (O2)}$, ${\bf (O6)}$ and ${\bf (O7)}$, and \eqref{three:bound}. 


\begin{corollary}\label{three:cor2}
Assume the conditions of Theorem \ref{three:twth2} hold. Then, we have the following results
\begin{itemize}
\item[(i)] $\begin{aligned}[t]
d_{TV}\left(X_1,S_{k_1,k_2}^n\right)&\le\frac{1}{\max(1,\lambda)}\left\{2 (2+qp)h_1(n,k,p)a(p)^3+|\varphi|\right\}, ~~{{\rm for}~X_1{\rm~defined~in~\eqref{three:p}}}.
\end{aligned}$
\item[(ii)] $\begin{aligned}[t]
d_{TV}\left(X_2,S_{k_1,k_2}^n\right) & \le\frac{1}{{\floor{\check{\alpha}}}\check{p}\check{q}} \left\{2(2+qp)a(p)^2\left(h_1(n,k,p)a(p)+\check{p}h_2(n,k,p)\right)+|\varphi|\right\},~~{{\rm for}~X_2{\rm~defined~in~\eqref{three:pb}}}.
\end{aligned}$
\item[(iii)] $\begin{aligned}[t]
d_{TV}\left(X_3, S_{k_1,k_2}^n\right) & \le \frac{1}{\hat{\alpha}\hat{q}}\left\{2(2+qp)a(p)^2\left(\hat{p}h_1(n,k,p)a(p)+\hat{q}h_2(n,k,p)\right)+\hat{p}|\varphi|\right\},~~{{\rm for}~X_3{\rm~defined~in~\eqref{three:nb}}}.
\end{aligned}$
\end{itemize} 
\end{corollary}

\begin{remark}
In a similar spirit, the approximation results for various distributions satisfying \eqref{three:DGM1} and  \eqref{panj} can be derived.
\end{remark}

\noindent
Next, we present the approximation results for the distribution of $M_{k_1,k_2}^n$, $S_{k_1,k_2}^n$ and $Z$.


\begin{theorem}\label{three:noniidth}
Let $\bar{g}_0$ be the solution to the Stein equation for DGM ($Z$ as defined in \eqref{three:DGM1}) with
$$\|\Delta \bar{g}_0\|\le \bar{w}_1\|f\|\min(1,\bar{\gamma}^{-1}),$$
where $\bar{w}_1, \bar{\gamma}>0$. 
Assume that $n\ge 5k$, $\bar{w}_1\bar{w}_2<\bar{\gamma}$, ${\mathbb E}\big(S_{k_1,k_2}^n\big)={\mathbb E}\big(M_{k_1,k_2}^n\big)$ and $\tau=\mathrm{Var}\big(M_{k_1,k_2}^n\big)-\mathrm{Var}\big(S_{k_1,k_2}^n\big)$. Then
\begin{align*}
d_{TV}\left(Z,M_{k_1,k_2}^n\right)&\le\frac{\bar{\gamma} }{2(\bar{\gamma} - \bar{w}_1 \bar{w}_2)}\left( \bar{w}_1 \bar{\varepsilon}^*\min\left(1,\bar{\gamma}^{-1}\right)+2{\mathbb P}(S_{k_1,k_2}^n>N)+2{\mathbb P}(M_{k_1,k_2}^n>N)\right)+d_{TV}(Z, S_{k_1,k_2}^n),
\end{align*}
where $d_{TV}(Z, S_{k_1,k_2}^n)$ is given in \eqref{three:twoth}, and 
\begin{align*}
\bar{w}_2&=(2+qp)a(p)\left\{(n-k)(|1-wb|a(p)\delta+|wb|) +|1-wb|a(p)\left(\delta^*c_{n,k}^{(1)}+c_{n,k}^{(2)}\right) \right\}~{\rm and}\\
\bar{\varepsilon}^*&\hspace{-0.07cm}=\hspace{-0.07cm} 2\Bigg\{|1\hspace{-0.07cm}-\hspace{-0.07cm}wb|\Bigg[\sum_{i=1}^{n-k}\sum_{|j-i|\le k+1}\hspace{-0.07cm}\Bigg(\hspace{-0.15cm} \left(\sum_{u=i-2k-2}^{j-1}\hspace{-0.1cm}+\hspace{-0.1cm}\sum_{u=i+k+2}^{i+2k+2}\hspace{-0.07cm}\right)[{\mathbb E}({\mathbb I}_i){\mathbb E}({\mathbb I}_j {\mathbb I}_u)\hspace{-0.07cm}+\hspace{-0.07cm}{\mathbb E}({\mathbb I}_i{\mathbb I}_j{\mathbb I}_u)]\hspace{-0.07cm}+\hspace{-0.07cm}[{\mathbb E}({\mathbb I}_i){\mathbb E}({\mathbb I}_j)\hspace{-0.07cm}+\hspace{-0.07cm}{\mathbb E}({\mathbb I}_i{\mathbb I}_j)]\hspace{-0.4cm}\sum_{|u-i|\le 2k+2}\hspace{-0.4cm}{\mathbb E}({\mathbb I}_u)\hspace{-0.07cm}\Bigg)\\
&~~~+\frac{|\tau|}{2}\hspace{-0.07cm}+\hspace{-0.07cm}(2\hspace{-0.07cm}+\hspace{-0.07cm}qp)h_1(n,k,p)a(p)^3\Bigg]\hspace{-0.07cm}+\hspace{-0.07cm}|wb|\Bigg[\sum_{i=1}^{n-k}\Bigg(\hspace{-0.07cm}\sum_{|j-i|\le k+1}\hspace{-0.3cm}{\mathbb E}({\mathbb I}_i{\mathbb I}_j)\hspace{-0.07cm}+\hspace{-0.07cm}{\mathbb E}({\mathbb I}_i)\hspace{-0.4cm}\sum_{|j-i|\le k+1}\hspace{-0.3cm}{\mathbb E}({\mathbb I}_j)\Bigg)\hspace{-0.07cm}+\hspace{-0.07cm}(2+qp)h_2(n,k,p)a(p)^2\Bigg] \hspace{-0.07cm}\Bigg\}\hspace{-0.02cm}.
\end{align*}
\end{theorem}

\begin{remarks}
\begin{itemize}
\item[(i)] The results for Poisson, pseudo-binomial and negative binomial approximations follow from Theorem \ref{three:noniidth}, Corollaries \ref{three:cor2}, ${\bf (O2)}$, ${\bf (O6)}$ and ${\bf (O7)}$, and \eqref{three:bound}. Also, note that similar bounds can be obtained for other distributions satisfying \eqref{three:DGM1} and  \eqref{panj}.

\item[(ii)] 
Observe that ${\mathbb P}(S_{k_1,k_2}^n>N)$ and ${\mathbb P}(M_{k_1,k_2}^n>N)$ are zero for Poisson and negative binomial distributions. Also, if we take $\check{\alpha}>\floor{n/k}$ then these probabilities are zero for pseudo-binomial distribution.

\item[(iii)] Observe that, if $\mathrm{Var}(S_{k_1,k_2}^n)=\mathrm{Var}\big(M_{k_1,k_2}^n\big)$ (i.e., $\tau=0$) then the bound become sharper, as expected.

\item[(iv)] Using the condition ${\mathbb E}\big(S_{k_1,k_2}^n\big)={\mathbb E}\big(M_{k_1,k_2}^n\big)$ $\left(i.e.,~q(1+(n-k-1)p)a(p)=\sum_{i=1}^{n-k}{\mathbb E}({\mathbb I}_i)\right)$, for any values of $n$, $k$ and $p_i$, we can estimate the value of $p$ or $q$ (see Table \ref{tab:table3}). Hence, with the estimated value of $p$ or $q$, the bounds can be easily computed.
\end{itemize}
\end{remarks}

\section{An Application to Two-stage Start-up Demonstration Testing}\label{three:APP}
Let us consider a scenario in which a customer is interested in buying a certain equipment. In such a case, Balakrishnan and Chan \cite{BC} recommend two-stage start-up demonstration test (see also Balakrishnan and Koutras \cite{BK} and references therein). During the second stage of this test, we need to count specific pattern of consecutive failures (say $k_1$) and successes (say $k_2$). The customer accepts the equipment under testing, whenever a recommended count of this pattern is observed. Using this mechanism, the customer can check the performance of the equipment, and can decide to reject bad equipment at early (first) stage or put strict norms for acceptance of the equipment in second stage.
In particular, the decision criteria in two-stage start-up test, for the customer, proposed by Balakrishnan and Chan \cite{BC} (see also Balakrishnan and Koutras \cite{BK}, p.-281) is as follows:
\begin{itemize}
\item[(i)] Accept the equipment (in the first stage) if a run of $k_1^*$ successes occurs before $\ell_1^*$ failures.
\item[(ii)] Accept the equipment if a run of $k_1^*$ successes does not occur before $\ell_1^*$ failures, but a run of $k_1$ successes followed by $k_2$ failures is observed;
and
\item[(iii)] Reject the item if a run of $k_1^*$ successes does not occur before $\ell_1^*$ failures, and also a run of $k_1$ successes followed by $k_2$ failures is not observed.
\end{itemize}
Assume that there is no run of $k_1^*$ successes followed by $\ell_1^*$ failures in the first stage.
In second stage, the problem studied in this paper becomes relevant, whenever $k_1$ consecutive failures followed by $k_2$ consecutive successes are considered. Then, by interchanging the role of success and failure, the setup can be adapted to our setting and the approximation results can be applied rather than computing complicated distributions. Now, we compare approximation results for particular values of $n,k_1,k_2$ and $p$, and demonstrate the closeness of the target distribution with approximate distribution. Also, we choose $\check{p}=0.001$ and $\hat{p}=0.999$ for pseudo-binomial and negative binomial distributions, respectively, in the following table.
\begin{table}[ht!]
  \centering
  \caption{Bounds for Poisson, pseudo-binomial and negative binomial distributions under iid setup.}
   \label{tab:table1}
\vspace{0.5cm}
  \begin{tabular}{cccccccccc}
    \toprule
Approximation & ($k_1,k_2$) & $n$  & $p=0.10$ & $p=0.30$ & $p=0.50$ & $p=0.70$ & $p=0.90$\\
    \midrule
Poisson           & \multirow{3}{*}{(3,4)}  & \multirow{3}{*}{50} & $8.0 \times 10^{-8}$  & 0.0042107 & 0.0804057 & 0.0344901 & 0.0000222\\
Pseudo-binomial   &                         &                     & 0.0013278 & 0.1601320 & 0.9396260 & 0.5864460 & 0.0103829 \\
Negative binomial &                         &                     & 0.0013028 & 0.1619750 & 0.9406920 & 0.5878640 & 0.0099867 \\
\midrule
Poisson           & \multirow{3}{*}{(3,5)}  & \multirow{3}{*}{150}& $4.7 \times 10^{-10}$ & 0.0005335 & 0.0453505 & 0.0528375 & 0.0000712\\
Pseudo-binomial   &                         &                     & 0.0010183 & 0.0207954 & 0.3265200 & 0.3915960 & 0.0111253\\
Negative binomial &                         &                     & 0.0010156 & 0.0227762 & 0.3281970 & 0.3932090 & 0.0108391 \\
\midrule

Poisson           & \multirow{3}{*}{(4,5)}  & \multirow{3}{*}{250}& $5.5\times 10^{-10}$  & 0.0003651 & 0.0124497 & 0.0036038 & $2.2\times 10^{-7}$\\
Pseudo-binomial   &                         &                     & 0.0010176 & 0.0126423 & 0.1027660 & 0.0449730 & 0.0011935 \\
Negative binomial &                         &                     & 0.0010171 & 0.0146311 & 0.1046660 & 0.0469300 & 0.0011844\\
    \bottomrule
  \end{tabular}
\end{table}
\newpage
\noindent
Observe that the bounds for pseudo-binomial and negative binomial distributions depend on free parameter ($\check{\alpha}$ or $\check{p}$ and $\hat{\alpha}$ or $\hat{p}$, respectively) and can be obtain much sharper form according to the value of free parameter. \\ 
Next, let us consider the case when $M_{k_1,k_2}^n$ arising from non-identical Bernoulli trials with $p_i$ as follows:
\begin{table}[ht!]
  \centering
  \begin{tabular}{|c|c|c|c|c|c|c|c|c|c|}
\hline
$i$ & $p_i$ & $i$ & $p_i$ & $i$ & $p_i$ & $i$ & $p_i$ & $i$ & $p_i$\\
\hline
1-30 & 0.15 & 61-90 & 0.17 & 121-150 & 0.19 & 181-210 & 0.21 &241-250 &0.23\\
\hline
31-60 & 0.16 & 91-120 & 0.18 & 151-180 & 0.20 & 211-240 & 0.22 & &\\
\hline
    \end{tabular}
\end{table}\\
Now, we demonstrate the results for Poisson approximation under non-identical setup and similarly, we can obtain bounds for other random variables which satisfy \eqref{panj} with $\bar{w}_1\bar{w}_2<\bar{\gamma}$. To compute $d_{TV}\big(S_{k_1,k_2}^n,M_{k_1,k_2}^n\big)$ and $d_{TV}\big(X_1,S_{k_1,k_2}^n\big)$, we estimate the value of $q$ from the relation $q(1+(n-k-1)p)a(p)=\sum_{i=1}^{n-k}{\mathbb E}({\mathbb I}_i)$ and then obtain the bounds by adding these two bounds.
\begin{table}[ht!]
  \centering
  \caption{Poisson bounds for non-identical trials.}
   \label{tab:table3}
\vspace{0.5cm}
  \begin{tabular}{cccccccccc}
    \toprule
($k_1,k_2$) & $n$ & $q$ & $d_{TV}\big(X_1,S_{k_1,k_2}^n\big)$ & $d_{TV}\big(S_{k_1,k_2}^n,M_{k_1,k_2}^n\big)$ & $d_{TV}\big(X_1,M_{k_1,k_2}^n\big)$\\
    \midrule
\vspace{0.1cm}
(3,4)  &50 & 0.0652339752 & $1.05\times 10^{-9}$ &$1.44\times 10^{-6}$ & $1.44105\times 10^{-6}$\\
\vspace{0.1cm}
(3,5)  & 150& 0.0289666557 &  $2.07\times 10^{-15}$ &$6.79\times 10^{-9}$ & $6.79000\times 10^{-9}$\\
(4,5)  &250& 0.0526121138 &  $1.20\times 10^{-12}$ &$3.10\times 10^{-9}$ & $3.10120\times 10^{-9}$\\
\bottomrule
    \end{tabular}
\end{table}\\
Note that the bounds for non-identical trials are larger than the bounds for identical trials, as expected.

\section{Auxiliary Results}\label{ARR}
\noindent
In this section, we derive some results related to the distribution of $S_{k_1,k_2}^n$. We first derive exact distribution of $S_{k_1,k_2}^n$ (using Markov chain approach, see Fu and Koutras \cite{FK}, Koutras \cite{koutras}, Balakrishnan and Koutras \cite{BK} and Dafnis {\em et al.} \cite{DAP}) in the following lemma.
\begin{lemma}\label{three:le1}
Let $\phi_n(\cdot)$, $\Phi(\cdot,\cdot)$ and $p_{\cdot,n}={\mathbb P}(S_{k_1,k_2}^n=\cdot)$ denote single PGF, double PGF and PMF of $S_{k_1,k_2}^n$. Then,
\begin{itemize}
\item[(i)] $\begin{aligned}[t]
\Phi(t,z)=\sum_{n=0}^{\infty}\phi_n(t)z^n= \frac{1+a(p)z^k(1-qz)(1-t)}{1-z+a(p)z^k(1-t)(1-qz)(1-pz)}.
\end{aligned}$
\item[(ii)] $\begin{aligned}[t]
p_{m,n}=\left\{
          \begin{array}{ll}
          0 & n \le k,~m>0\\
          1 & n \le k,~m=0\\
          qa(p) & n=k+1,~m=1,\\
         1-qa(p) & n=k+1,~m=0,\\
        p_{m,n-1}-a(p)[(p_{m,n-k}-p_{m-1,n-k})\\
       -(p_{m,n-k-1}-p_{m-1,n-k-1})+qp(p_{m,n-k-2}-p_{m-1,n-k-2})] & n \ge k+2,~m \ge 0.
         \end{array}\right.
\end{aligned}$
\end{itemize}
\end{lemma}
\noindent
{\bf Proof}. Let $\ell_n =\sup\{x:{\mathbb P}(S_{k_1,k_2}^n=x)>0\}= \floor{n/k}$ and $C_x=\{c_{x,0},c_{x,1},\dotsc,c_{x,k_1+k_2+1}\}$, where $c_{x,i}=(x,i),~0 \le i \le k_1+k_2+1$. Also, define a Markov chain $\{Y_t: t \ge 0\}$ on $\Omega=\cup_{x=0}^{\ell_n}C_x$ as $Y_t =(x,j)$ if ($k_1,k_2$)-event has occurred $x$ times in the first $t$ outcomes and
\begin{itemize}
\item $j=0$, if $\eta_t = 1$ and there is no pattern with exactly $k_1$ failures before $\eta_t$.
\item $j = i$, $1 \le  i \le  k_1$, if $\eta_t= \eta_{t-1}=\dotsb = \eta_{t-i+1}=0$ and the $(t-i)$-th outcome is a success (if exists).
\item $j = k_1^+$, if there are more than $k_1$ consecutive failures. i.e., there exists a positive integers $l \ge k_1+1$ such that $\eta_t=\eta_{t-1}= \dotsb = \eta_{t-l+1}=0$.
\item $j = k_1 + i$, $1\le i \le  k_2$, if $\eta_t = \eta_{t-1} = \dotsb = \eta_{t-i+1}=1$, $\eta_{t-i}= \eta_{t-i-1}=\dotsb=\eta_{t-i-(k_1-1)}=0$ and the ($t-i-k_1$)-th outcome is a success (if exists).
\end{itemize}
Also, we say the pattern is complete if exactly $k_1$ consecutive failures followed by exactly $k_2$ consecutive successes and a failure occurs at $(k_1+k_2+1)$-th position. Now, $S_{k_1,k_2}^n$ becomes Markov chain embeddable variable of binomial type (MVB, see Dafnis {\em et al.} \cite{DAP} for more details) with this setup and $\pi_0=(1,0,\dotsc,0)_{1\times (k_1+k_2+2)}$,
$$
A=
\begin{blockarray}{cccccccccccc}
    \matindex{($\cdot,0$)} & \matindex{($\cdot,1$)} &\matindex{($\cdot,2$)}&\matindex{$\cdot$}&\matindex{($\cdot,k_1-1$)}&\matindex{($\cdot,k_1$)}&\matindex{($\cdot,k_1^+$)}&\matindex{($\cdot,k_1+1$)}&\matindex{($\cdot,k_1+2$)}&\matindex{$\cdot$}&\matindex{($\cdot,k_1+k_2-1$)}&\matindex{($\cdot,~k_1+k_2$)} \\
    \begin{block}{(cccccccccccc)}
p & q & 0 & \cdot & 0 & 0 & 0 & 0 & 0 & \cdot & 0 & 0\\
p & 0 & q & \cdot & 0 & 0 & 0 & 0 & 0 & \cdot & 0 & 0\\
p & 0 & 0 & \cdot & 0 & 0 & 0 & 0 & 0 & \cdot & 0 & 0\\
\cdot & \cdot & \cdot & \cdot & \cdot & \cdot & \cdot & \cdot & \cdot & \cdot & \cdot & \cdot\\
p & 0 & 0 & \cdot & 0 & q & 0 & 0 & 0 & \cdot & 0 & 0\\
0 & 0 & 0 & \cdot & 0 & 0 & q & p & 0 & \cdot & 0 & 0\\
p & 0 & 0 & \cdot & 0 & 0 & q & 0 & 0 & \cdot & 0 & 0\\
0 & q & 0 & \cdot & 0 & 0 & 0 & 0 & p & \cdot & 0 & 0\\
0 & q & 0 & \cdot & 0 & 0 & 0 & 0 & 0 & \cdot & 0 & 0\\
\cdot & \cdot & \cdot& \cdot & \cdot & \cdot & \cdot & \cdot & \cdot & \cdot & \cdot & \cdot \\
0 & q & 0 & \cdot & 0 & 0 & 0 & 0 & 0 &\cdot&0 & p\\
p & 0 & 0 & \cdot & 0 & 0 & 0 & 0 & 0 & \cdot & 0 & 0\\
 \end{block}
  \end{blockarray}
$$
where $A$ is $(k_1+k_2+2) \times (k_1+k_2+2)$ matrix and $B$ is the matrix with zero entries except at $(k_1+k_2+2,2)$ where it is equal to $q$.
Therefore, from Theorem $2.1$ of Dafnis {\em et al.} \cite{DAP}, with some calculations, we get the required result. Also, $(ii)$ directly follows from $(i)$.\qed\\
\noindent 
In particular, let us demonstrate the result for $k_1=2$ and $k_2=5$. Here, $\pi_0=(1,0,0,0,0,0,0,0,0)_{1\times 9}$,
$$A=
\begin{pmatrix}
p&q&0&0&0&0&0&0&0\\
p&0&q&0&0&0&0&0&0\\
0&0&0&q&p&0&0&0&0\\
p&0&0&q&0&0&0&0&0\\
0&q&0&0&0&p&0&0&0\\
0&q&0&0&0&0&p&0&0\\
0&q&0&0&0&0&0&p&0\\
0&q&0&0&0&0&0&0&p\\
p&0&0&0&0&0&0&0&0
\end{pmatrix}
\quad \text{and}\quad
B=\begin{pmatrix}
0&0&0&0&0&0&0&0&0\\
0&0&0&0&0&0&0&0&0\\
0&0&0&0&0&0&0&0&0\\
0&0&0&0&0&0&0&0&0\\
0&0&0&0&0&0&0&0&0\\
0&0&0&0&0&0&0&0&0\\
0&0&0&0&0&0&0&0&0\\
0&0&0&0&0&0&0&0&0\\
0&q&0&0&0&0&0&0&0
\end{pmatrix}
$$
Then, from Theorem $2.1$ of Dafnis {\em et al.} \cite{DAP}, we have
$$\Phi(t,z)=\sum_{n=0}^{\infty}\phi_n(t)z^n= \frac{1+(qz)^2(pz)^5(1-qz)(1-t)}{1-z+(qz)^2(pz)^5(1-t)(1-qz)(1-pz)}.$$

\noindent 
Next, using Lemma \ref{three:le1} $(i)$, we can derive the following relation between $\phi_n$ and $\phi_n^\prime$.

\begin{lemma}\label{three:relpgf}
The PGF of $S_{k_1,k_2}^n$, $\phi_n(\cdot)$, satisfies the following recursive relation
$$\displaystyle{\phi_{n}^{\prime}(t)=a(p)\sum_{i=1}^{3}a_i\sum_{s=0}^{d_i}b_i(n-s)C_s(t)\phi_{n-k-s-i+1}(t)},$$
where 
$$\displaystyle{C_s(t)=\sum_{l=0}^{\floor{s/k}}\sum_{m=0}^{\floor{\frac{s-lk}{k+1}}}{s-l(k-1)-mk \choose s-lk-m(k+1),l,m}\frac{(k+1)^{s-lk-m(k+1)}}{(k+2)^{s-l(k-1)-mk+1}}(-1)^m 2^l [a(p)(t-1)]^{l+m}}$$ 
and $\displaystyle{n \choose a, b, c}=\frac{n!}{a! b! c!}$ with $n=a+b+c$.
\end{lemma}
\noindent
{\bf Proof}. From Lemma \ref{three:le1} (i), the double PGF of $S_{k_1,k_2}^n$ can be written as
\begin{equation}
[1-z-a(p)(t-1)z^k(1-qz)(1-pz))]\sum_{n=0}^{\infty}\phi_n(t)z^n=1-a(p)(t-1)z^k(1-qz).\label{three:pgf1}
\end{equation}
Differentiating \eqref{three:pgf1} w.r.t. $t$ and $z$, we have
\begin{align}
&[1-z-a(p)(t-1)z^k(1-z+qpz^2)]\sum_{n=0}^{\infty}\phi_n^\prime(t)z^n-a(p)z^k[1-z+qpz^{2}]\sum_{n=0}^{\infty}\phi_n(t)z^n = -a(p)z^k(1-qz)\label{three:pgf2}\allowdisplaybreaks\\
&[1-z-a(p)(t-1)z^k(1-z+qpz^2)]\sum_{n=0}^{\infty}n\phi_n(t)z^n\hspace{-0.07cm}-\hspace{-0.07cm}[z\hspace{-0.07cm}+\hspace{-0.07cm}a(p)(t\hspace{-0.07cm}-\hspace{-0.07cm}1)z^k(k-(k+1)z+qp(k+2)z^{2})]\sum_{n=0}^{\infty}\phi_n(t)z^n\nonumber\\
&~~~~~~~~~~~~~~~~~~~~~~~~~~~~~~~~~~~~~~~~~~~~~~~~~~~~~~~~~~~~~~~~~~~~~~~~~~~~~~~~~~~~~~~~~~~~=-a(p)(t-1)z^k(k-q(k+1)z).\label{three:pgf3}
\end{align}
Multiplying by $(z+a(p)(t-1)z^k(k-(k+1)z+qp(k+2)z^{2}))$ in \eqref{three:pgf2}, $a(p)z^k(1-z+qpz^{2})$ in \eqref{three:pgf3} and subtracting, we get
\begin{align}
[z+a(p)(t-1)z^k(k&-(k+1)z+qp(k+2)z^{2})]\sum_{n=0}^{\infty}\phi_n^\prime(t)z^n-a(p)z^k(1-z+qpz^2)\sum_{n=0}^{\infty}n\phi_n(t)z^n\nonumber\\
&~~~~~~~~~~~~~~~~~~~~~~~~~~~~~~~~~~~~~~~~~~~~~~=\frac{a(p)z^{k+1}(1-qz)(pa(p)(t-1)z^k(1-qz)-1)}{1-z-a(p)(t-1)z^k(1-qz)(1-pz))}.\label{three:pgf4}
\end{align}
Multiplying by $(k+2)$ in \eqref{three:pgf2} and adding with \eqref{three:pgf4}, we have
\begin{align}
[(k+2)&-(k+1)z-a(p)(t-1)z^k(2-z)]\sum_{n=0}^{\infty}\phi_n^\prime(t)z^n=qa(p)(k+2)z^{k+1}\nonumber\\
&~+a(p)\left[\sum_{n=k+2}^{\infty}(n+1)\phi_{n-k}(t)z^n-\sum_{n=k+2}^{\infty}(n+1-q)\phi_{n-k-1}(t)z^n+qp\sum_{n=k+2}^{\infty}n\phi_{n-k-2}(t)z^n\right].\label{three:pgf6}
\end{align}
Consider
\begin{align*}
&\frac{1}{(k+2)-(k+1)z-a(p)(t-1)z^k(2-z)}\\
&=\frac{1}{(k+2)\left(1-\left(\frac{k+1}{k+2}\right)z-\frac{1}{k+2}a(p)(t-1)z^k(2-z)\right)}\\
&=\frac{1}{(k+2)}\sum_{n=0}^{\infty}\left(\left(\frac{k+1}{k+2}\right)z-\frac{1}{k+2}a(p)(t-1)z^k(2-z)\right)^n\\
&=\sum_{n=0}^{\infty}\frac{1}{(k+2)^{n+1}}\left[\left(k+1\right)-a(p)(t-1)z^{k-1}(2-z)\right]^n z^n\\
&=\sum_{n=0}^{\infty}\left(\sum_{l=0}^{n}{n \choose l}\left(k+1\right)^{n-l}a(p)^l(t-1)^l(2-z)^l\right)\frac{z^{n+l(k-1)}}{(k+2)^{n+1}}\allowdisplaybreaks\\
&=\sum_{n=0}^{\infty}\sum_{l=0}^{n}{n \choose l}\left(k+1\right)^{n-l}a(p)^l(t-1)^l\left(\sum_{m=0}^{l}{l \choose m}2^{l-m}(-1)^m z^m\right)\frac{z^{n+l(k-1)}}{(k+2)^{n+1}}\\
&=\sum_{n=0}^{\infty}\sum_{l=0}^{n}\sum_{m=0}^{l}{n \choose l}{l \choose m}\frac{(k+1)^{n-l}}{(k+2)^{n+1}}(-1)^m 2^{l-m} a(p)^l(t-1)^lz^{n+m+l(k-1)}\\
&=\sum_{l=0}^{\infty}\sum_{n=l}^{\infty}\sum_{m=0}^{l}{n \choose l}{l \choose m}\frac{(k+1)^{n-l}}{(k+2)^{n+1}}(-1)^m 2^{l-m} a(p)^l(t-1)^lz^{n+m+l(k-1)}\\
&=\sum_{l=0}^{\infty}\sum_{n=0}^{\infty}\sum_{m=0}^{l}{n+l \choose l}{l \choose m}\frac{(k+1)^{n}}{(k+2)^{n+l+1}}(-1)^m 2^{l-m} a(p)^l(t-1)^lz^{n+m+lk}\\
&=\sum_{n=0}^{\infty}\sum_{m=0}^{\infty}\sum_{l=m}^{\infty}{n+l \choose l}{l \choose m}\frac{(k+1)^{n}}{(k+2)^{n+l+1}}(-1)^m 2^{l-m} a(p)^l(t-1)^lz^{n+m+lk}\\
&=\sum_{n=0}^{\infty}\sum_{m=0}^{\infty}\sum_{l=0}^{\infty}{n+l+m \choose n,l,m}\frac{(k+1)^{n}}{(k+2)^{n+l+m+1}}(-1)^m 2^{l} a(p)^{l+m}(t-1)^{l+m}z^{n+lk+m(k+1)}\\
&=\sum_{m=0}^{\infty}\sum_{n=m(k+1)}^{\infty}\sum_{l=0}^{\infty}{n+l-mk \choose n-m(k+1),l,m}\frac{(k+1)^{n-m(k+1)}}{(k+2)^{n+l-mk+1}}(-1)^m 2^{l} a(p)^{l+m}(t-1)^{l+m}z^{n+lk}\\
&=\sum_{n=0}^{\infty}\sum_{l=0}^{\infty}\sum_{m=0}^{\floor{\frac{n}{k+1}}}{n+l-mk \choose n-m(k+1),l,m}\frac{(k+1)^{n-m(k+1)}}{(k+2)^{n+l-mk+1}}(-1)^m 2^{l} a(p)^{l+m}(t-1)^{l+m}z^{n+lk}\\
&=\sum_{l=0}^{\infty}\sum_{n=lk}^{\infty}\sum_{m=0}^{\floor{\frac{n-lk}{k+1}}}{n-l(k-1)-mk \choose n-lk-m(k+1),l,m}\frac{(k+1)^{n-lk-m(k+1)}}{(k+2)^{n-l(k-1)-mk+1}}(-1)^m 2^{l} a(p)^{l+m}(t-1)^{l+m}z^{n}\\
&=\sum_{n=0}^{\infty}\sum_{l=0}^{\floor{\frac{n}{k}}}\sum_{m=0}^{\floor{\frac{n-lk}{k+1}}}{n-l(k-1)-mk \choose n-lk-m(k+1),l,m}\frac{(k+1)^{n-lk-m(k+1)}}{(k+2)^{n-l(k-1)-mk+1}}(-1)^m 2^{l} a(p)^{l+m}(t-1)^{l+m}z^{n}\\
&=\sum_{n=0}^{\infty}C_n(t)z^n.
\end{align*}
Therefore,
\begin{align*}
&\frac{1}{(k+2)-(k+1)z-a(p)(t-1)z^k(2-z)}=\sum_{n=0}^{\infty}C_n(t)z^n.
\end{align*}
Substituting in \eqref{three:pgf6} and comparing the coefficients of $z^n$, we get
\begin{align*}
\phi_{n}^{\prime}(t)&=a(p)\left[\sum_{s=0}^{n-k-2}\hspace{-.2cm}(n-s+1)C_s(t)\phi_{n-k-s}(t)\hspace{-.1cm}-\hspace{-.3cm}\sum_{s=0}^{n-k-1}\hspace{-.2cm}b_2(n-s)C_s(t)\phi_{n-k-s-1}(t)+qp\hspace{-.3cm}\sum_{s=0}^{n-k-2}\hspace{-.2cm}(n-s)C_s(t)\phi_{n-k-s-2}(t)\right]\nonumber\\
&=a(p)\sum_{i=1}^{3}a_i\sum_{s=0}^{d_i}b_i(n-s)C_s(t)\phi_{n-k-s-i+1}(t).
\end{align*}
This proves the result.\qed
\begin{remark}

Observe that 
\begin{align}
C_s(t)&=\sum_{l=0}^{\floor{s/k}}\sum_{r=0}^{\floor{\frac{s-lk}{k+1}}}{s-l(k-1)-rk \choose s-lk-r(k+1),l,r}\frac{(k+1)^{s-lk-r(k+1)}}{(k+2)^{s-l(k-1)-rk+1}}(-1)^r 2^l [a(p)(t-1)]^{l+r}\nonumber\\
&=\sum_{l=0}^{\floor{s/k}}\sum_{r=0}^{\floor{\frac{s-lk}{k+1}}}\sum_{m=0}^{l+r} {s-l(k-1)-rk \choose s-lk-r(k+1,l,r)}{l+r \choose m}\frac{(k+1)^{s-lk-r(k+1)}}{(k+2)^{s-l(k-1)-rk+1}} (-1)^{l-m} 2^l a(p)^{l+r}t^m\nonumber\allowdisplaybreaks\\
&=\sum_{l=0}^{\floor{s/k}}\sum_{m=0}^{\floor{\frac{s+l}{k+1}}}\sum_{r=m-l}^{\floor{\frac{s-lk}{k+1}}} {s-l(k-1)-rk \choose s-lk-r(k+1,l,r)}{l+r \choose m} \frac{(k+1)^{s-lk-r(k+1)}}{(k+2)^{s-l(k-1)-rk+1}}(-1)^{l-m} 2^l a(p)^{l+r}t^m\nonumber\\
&=\sum_{m=0}^{\floor{s/k}}\left[\sum_{l=m(k+1)-s}^{\floor{s/k}}\sum_{r=m-l}^{\floor{\frac{s-lk}{k+1}}} {s-l(k-1)-rk \choose s-lk-r(k+1,l,r)}{l+r \choose m}\frac{(k+1)^{s-lk-r(k+1)}}{(k+2)^{s-l(k-1)-rk+1}} (-1)^{l-m} 2^l a(p)^{l+r}\right]t^m\nonumber\\
&=\sum_{m=0}^{\floor{s/k}}B_s(m)t^m,\label{three:cst}
\end{align}
where
$$\displaystyle{B_s(m)=\sum_{l=m(k+1)-s}^{\floor{s/k}}\sum_{r=m-l}^{\floor{\frac{s-lk}{k+1}}}{s-l(k-1)-rk \choose s-lk-r(k+1),l,r}{l+r \choose m}\frac{(k+1)^{s-lk-r(k+1)}}{(k+2)^{s-l(k-1)-rk+1}}(-1)^{l-m}2^l a(p)^{l+r}}.$$
Further, the expression can be expressed as
$$C_n(t)=c_0+c_1(t-1)+\dotsb+c_{\floor{n/k}} (t-1)^{\floor{n/k}}=\sum_{m=0}^{\floor{n/k}}B_n(m)t^m,$$
Though the form looks complicated, we only need $C_n(1)=\sum_{m=0}^{\floor{n/k}}B_n(m)=c_0$, $C_n^\prime(1)=\sum_{m=0}^{\floor{n/k}}mB_n(m)=c_1$ and $C_n^{\prime\prime}(1)=\sum_{m=0}^{\floor{n/k}}m(m-1)B_n(m)=2c_2$ to derive the approximation results and they are easy to compute.
\end{remark}

\section{Proofs}\label{three:pr}
\noindent
{\bf Proof of Theorem \ref{three:twth2}}. We know that
\begin{equation}
\phi_{n}(t)= \sum_{m=0}^{\floor{n/k}}p_{m,n}t^m\label{three:pgfp}
\end{equation}
and 
\begin{equation}
\phi_{n}^{\prime}(t)=\sum_{m=0}^{\floor{n/k}}m p_{m,n}t^{m-1}=\sum_{m=0}^{\floor{n/k}-1}(m+1)p_{m+1,n}t^m\quad\text{and}\quad C_s(t)=\sum_{m=0}^{\floor{s/k}}B_s(m)t^m.\label{three:pgfp1}
\end{equation}
Substituting \eqref{three:pgfp} and \eqref{three:pgfp1} in the recursive relation derived in Lemma \ref{three:relpgf}, we have
\begin{align*}
&\sum_{m=0}^{\floor{n/k}-1}(m+1)p_{m+1,n}t^m\\
&=a(p)\hspace{-0.05cm}\sum_{i=1}^{3}\hspace{-0.05cm}a_i\hspace{-0.05cm}\sum_{s=0}^{d_i}b_i(n-s)\hspace{-0.05cm}\left(\sum_{m=0}^{\infty}B_s(m){\bf 1}\left(m \hspace{-0.05cm}\le\hspace{-0.05cm} \floor{\frac{s}{k}}\right)t^m\right)\hspace{-0.05cm}\left(\sum_{m=0}^{\infty} p_{m,n-k-s-i+1}{\bf 1}\left(m \hspace{-0.05cm}\le\hspace{-0.05cm} \floor{\frac{n-k-s-i+1}{k}}\right)t^m\right)\\
&=\sum_{m=0}^{\infty}\left[a(p)\sum_{i=1}^{3}a_i\sum_{s=0}^{d_i}b_i(n-s)\sum_{l=0}^{m}p_{l,n-k-s-i+1}B_s(m-l){\bf 1}\left(l \le \floor{\frac{n-k-s-i+1}{k}}\right){\bf 1}\left(m-l \le \floor{\frac{s}{k}}\right)\right] t^m
\end{align*}
\vskip -0.05cm
\noindent
Multiplying by $(1-wbt)$ and collecting the coefficients of $t^m$, we get
\begin{align*}
&(m+1) p_{m+1,n}{\bf 1}\left(m \le \floor{n/k}-1\right)-wbm p_{m,n}{\bf 1}\left(m \le \floor{n/k}\right)\\
&~~~=a(p)\sum_{i=1}^{3}a_i\sum_{s=0}^{d_i}b_i(n-s)\sum_{l=0}^{m}p_{l,n-k-s-i+1}B_s(m-l){\bf 1}\left(l \le \floor{\frac{n\hspace{-0.07cm}-\hspace{-0.07cm}k\hspace{-0.07cm}-\hspace{-0.07cm}s\hspace{-0.07cm}-\hspace{-0.07cm}i\hspace{-0.07cm}+\hspace{-0.07cm}1}{k}}\right){\bf 1}\left(m-l \le \floor{\frac{s}{k}}\right)\\
&~~~~~-w b a(p)\sum_{i=1}^{3}a_i\sum_{s=0}^{d_i}b_i(n-s)\hspace{-0.1cm}\sum_{l=0}^{m-1}p_{l,n-k-s-i+1}B_s(m-l-1){\bf 1}\left(l \le \floor{\frac{n\hspace{-0.07cm}-\hspace{-0.07cm}k\hspace{-0.07cm}-\hspace{-0.07cm}s\hspace{-0.07cm}-\hspace{-0.07cm}i\hspace{-0.07cm}+\hspace{-0.07cm}1}{k}}\right){\bf 1}\left(m-l-1 \le \floor{\frac{s}{k}}\right)\hspace{-0.07cm}.
\end{align*}
\noindent
Let $g \in {\cal G}_{S_{k_1,k_2}^n}$, then
\vspace{-0.31cm}
\begin{align*}
&\sum_{m=0}^{\infty}g(m+1)\left[(m+1) p_{m+1,n}{\bf 1}\left(m \le \floor{n/k}-1\right)-wbm p_{m,n}{\bf 1}\left(m \le \floor{n/k}\right)\right]\\
&=\sum_{m=0}^{\infty}g(m+1)\left[a(p)\sum_{i=1}^{3}a_i\sum_{s=0}^{d_i}b_i(n-s)\sum_{l=0}^{m}p_{l,n-k-s-i+1}B_s(m-l){\bf 1}\left(l \le \floor{\frac{n\hspace{-0.07cm}-\hspace{-0.07cm}k\hspace{-0.07cm}-\hspace{-0.07cm}s\hspace{-0.07cm}-\hspace{-0.07cm}i\hspace{-0.07cm}+\hspace{-0.07cm}1}{k}}\right){\bf 1}\left(m-l \le \floor{\frac{s}{k}}\right)\right.\\
&~~\left.-w b a(p)\sum_{i=1}^{3}a_i\sum_{s=0}^{d_i}b_i(n-s)\sum_{l=0}^{m-1}p_{l,n-k-s-i+1}B_s(m-l-1){\bf 1}\left(l \le \floor{\frac{n\hspace{-0.07cm}-\hspace{-0.07cm}k\hspace{-0.07cm}-\hspace{-0.07cm}s\hspace{-0.07cm}-\hspace{-0.07cm}i\hspace{-0.07cm}+\hspace{-0.07cm}1}{k}}\right){\bf 1}\left(m-l-1 \le \floor{\frac{s}{k}}\right)\right]\allowdisplaybreaks\\
&=a(p)\sum_{i=1}^{3}a_i\sum_{s=0}^{d_i}b_i(n-s)\sum_{l=0}^{\infty}\sum_{m=l}^{\infty}g(m+1)p_{l,n-k-s-i+1}B_s(m-l){\bf 1}\left(l \le \floor{\frac{n\hspace{-0.07cm}-\hspace{-0.07cm}k\hspace{-0.07cm}-\hspace{-0.07cm}s\hspace{-0.07cm}-\hspace{-0.07cm}i\hspace{-0.07cm}+\hspace{-0.07cm}1}{k}}\right){\bf 1}\left(m-l \le \floor{\frac{s}{k}}\right)\\
&~~-w b a(p)\hspace{-0.07cm}\sum_{i=1}^{3}\hspace{-0.07cm}a_i\hspace{-0.07cm}\sum_{s=0}^{d_i}\hspace{-0.07cm}b_i(n-s)\hspace{-0.07cm}\sum_{l=0}^{\infty}\hspace{-0.07cm}\sum_{m=l+1}^{\infty}\hspace{-0.1cm}g(m\hspace{-0.07cm}+\hspace{-0.07cm}1)p_{l,n-k-s-i+1}B_s(m\hspace{-0.07cm}-\hspace{-0.07cm}l\hspace{-0.07cm}-\hspace{-0.07cm}1){\bf 1}\hspace{-0.07cm}\left(\hspace{-0.07cm}l \hspace{-0.07cm}\le\hspace{-0.07cm} \floor{\frac{n\hspace{-0.07cm}-\hspace{-0.07cm}k\hspace{-0.07cm}-\hspace{-0.07cm}s\hspace{-0.07cm}-\hspace{-0.07cm}i\hspace{-0.07cm}+\hspace{-0.07cm}1}{k}}\hspace{-0.07cm}\right)\hspace{-0.07cm}{\bf 1}\hspace{-0.07cm}\left(\hspace{-0.07cm}m\hspace{-0.07cm}-\hspace{-0.07cm}l\hspace{-0.07cm}-\hspace{-0.07cm}1 \hspace{-0.07cm}\le\hspace{-0.07cm} \floor{\frac{s}{k}}\hspace{-0.07cm}\right).
\end{align*}
This implies
\vspace{-0.31cm}
\begin{align*}
&\sum_{m=0}^{\floor{n/k}}\left[wbm g(m+1)-m g(m)\right]p_{m,n}+a(p)\sum_{i=1}^{3}a_i\sum_{s=0}^{d_i}b_i(n-s)\sum_{l=0}^{\floor{\frac{n-k-s-i+1}{k}}}\sum_{m=0}^{\floor{s/k}}g(m+l+1)p_{l,n-k-s-i+1}B_s(m)\\
&~~~~~~~~~~~~~~~~~~~~~~~~~~~~~~~~~~~~-w b a(p)\sum_{i=1}^{3}a_i\sum_{s=0}^{d_i}b_i(n-s)\sum_{l=0}^{\floor{\frac{n-k-s-i+1}{k}}}\sum_{m=0}^{\floor{s/k}}g(m+l+2)p_{l,n-k-s-i+1}B_s(m)=0.
\end{align*}
Interchanging $m$ and $l$ for second and third terms and using $e^{U(m+1)-U(m)}=a+bm$, we get
\vspace{-0.3cm}
\begin{align}
&\sum_{m=0}^{\floor{n/k}}\left[we^{U(m+1)-U(m)} g(m+1)-m g(m)\right]p_{m,n}-aw \sum_{m=0}^{\floor{n/k}}g(m+1)p_{m,n}\nonumber\\
&~~~~~~~~~~~~~~~~+a(p)\sum_{i=1}^{3}a_i\sum_{s=0}^{d_i}b_i(n-s)\sum_{l=0}^{\floor{s/k}}B_s(l)\sum_{m=0}^{\floor{\frac{n-k-s-i+1}{k}}}g(m+l+1)p_{m,n-k-s-i+1}\nonumber\\
&~~~~~~~~~~~~~~~~-w b a(p)\sum_{i=1}^{3}a_i\sum_{s=0}^{d_i}b_i(n-s)\sum_{l=0}^{\floor{s/k}}B_s(l)\sum_{m=0}^{\floor{\frac{n-k-s-i+1}{k}}}g(m+l+2)p_{m,n-k-s-i+1}=0.\nonumber
\end{align}
This implies
\vspace{-0.3cm}
\begin{align}
&\sum_{m=0}^{\floor{n/k}}\left[we^{U(m+1)-U(m)} g(m+1)-m g(m)\right]p_{m,n}-aw \sum_{m=0}^{\floor{n/k}}g(m+1)p_{m,n}\nonumber\\
&~~~~~~~~~~~~~~~~+a(p)\sum_{i=1}^{3}a_i\sum_{s=0}^{d_i}b_i(n-s)\sum_{l=0}^{\floor{s/k}}B_s(l){\mathbb E}\left[g(S_{k_1,k_2}^{n-k-s-i+1}+l+1)\right]\nonumber\\
&~~~~~~~~~~~~~~~~-w b a(p)\sum_{i=1}^{3}a_i\sum_{s=0}^{d_i}b_i(n-s)\sum_{l=0}^{\floor{s/k}}B_s(l){\mathbb E}\left[g(S_{k_1,k_2}^{n-k-s-i+1}+l+2)\right]=0.\label{b11}
\end{align}
It is known ${\mathbb E}[X]={\mathbb E}[{\mathbb E}(X|Y)]$, therefore
\vspace{-0.3cm}
\begin{align}
&\sum_{m=0}^{\floor{n/k}}\left[we^{U(m+1)-U(m)} g(m+1)-m g(m)\right]p_{m,n}-aw \sum_{m=0}^{\floor{n/k}}g(m+1)p_{m,n}\nonumber\\
&~~~~~~~~~~~~~~~~+a(p)\sum_{i=1}^{3}a_i\sum_{s=0}^{d_i}b_i(n-s)\sum_{l=0}^{\floor{s/k}}B_s(l){\mathbb E}\left[{\mathbb E}\left(g(S_{k_1,k_2}^{n-k-s-i+1}+l+1)|S_{k_1,k_2}^n\right)\right]\nonumber\\
&~~~~~~~~~~~~~~~~-w b a(p)\sum_{i=1}^{3}a_i\sum_{s=0}^{d_i}b_i(n-s)\sum_{l=0}^{\floor{s/k}}B_s(l){\mathbb E}\left[{\mathbb E}\left(g(S_{k_1,k_2}^{n-k-s-i+1}+l+2)|S_{k_1,k_2}^n\right)\right]=0.\nonumber
\end{align}
The expression can be rewritten as
\vspace{-0.24cm}
\begin{align}
&\sum_{m=0}^{\floor{n/k}}\Bigr[we^{U(m+1)-U(m)} g(m+1)-m g(m)-awg(m+1)\nonumber\\
&~~+a(p)\sum_{i=1}^{3}a_i\sum_{s=0}^{d_i}b_i(n-s)\sum_{l=0}^{\floor{s/k}}B_s(l){\mathbb E}\left\{g\left(S_{k_1,k_2}^{n-k-s-i+1}+l+1\right)\Bigr|S_{k_1,k_2}^n=m\right\}\nonumber\\
&~~-w b a(p)\sum_{i=1}^{3}a_i\sum_{s=0}^{d_i}b_i(n-s)\sum_{l=0}^{\floor{s/k}}B_s(l){\mathbb E}\left\{g\left(S_{k_1,k_2}^{n-k-s-i+1}+l+2\right)\Bigr|S_{k_1,k_2}^n=m\right\}\Bigr]p_{m,n}=0.\label{three:00}
\end{align}
Hence, Stein operator of $S_{k_1,k_2}^n$ is given by
\begin{align}
{\cal A}_{S_{k_1,k_2}^n}(g(m))&=we^{U(m+1)-U(m)} g(m+1)-m g(m)-awg(m+1)\nonumber\\
&~~+a(p)\sum_{i=1}^{3}a_i\sum_{s=0}^{d_i}b_i(n-s)\sum_{l=0}^{\floor{s/k}}B_s(l){\mathbb E}\left\{g\left(S_{k_1,k_2}^{n-k-s-i+1}+l+1\right)\Bigr|S_{k_1,k_2}^n=m\right\}\nonumber\\
&~~-w b a(p)\sum_{i=1}^{3}a_i\sum_{s=0}^{d_i}b_i(n-s)\sum_{l=0}^{\floor{s/k}}B_s(l){\mathbb E}\left\{g\left(S_{k_1,k_2}^{n-k-s-i+1}+l+2\right)\Bigr|S_{k_1,k_2}^n=m\right\}\nonumber\\
&={\cal A}_Z g(m)+{\cal U}g(m),\label{three:pp}
\end{align}
where ${\cal A}_Z$ is a Stein operator for DGM and ${\cal U}$ is a perturbed operator. Taking expectation of the perturbed operator ${\cal U}$ w.r.t. $S_{k_1,k_2}^n$ defined in \eqref{three:pp}, we have
\begin{align*}
{\mathbb E}\left[{\cal U} g\left(S_{k_1,k_2}^n\right)\right] \hspace{-0.07cm}=\hspace{-0.07cm} &-wa{\mathbb E}\left[g(S_{k_1,k_2}^n+1)\right]\hspace{-0.07cm}+\hspace{-0.07cm}a(p)\sum_{i=1}^{3}a_i\sum_{s=0}^{d_i}b_i(n\hspace{-0.07cm}-\hspace{-0.07cm}s)\hspace{-0.2cm}\sum_{l=0}^{\floor{s/k}}B_s(l){\mathbb E}\left[{\mathbb E}\left(g(S_{k_1,k_2}^{n-k-s-i+1}+l+1)|S_{k_1,k_2}^n\right)\right]\\
&-wb a(p)\sum_{i=1}^{3}a_i\sum_{s=0}^{d_i}b_i(n\hspace{-0.07cm}-\hspace{-0.07cm}s)\hspace{-0.07cm}\sum_{l=0}^{\floor{s/k}}B_s(l){\mathbb E}\left[{\mathbb E}\left(g(S_{k_1,k_2}^{n-k-s-i+1}+l+1)|S_{k_1,k_2}^n\right)\right]
\end{align*}
Now, using ${\mathbb E}[{\mathbb E}(X|Y)]={\mathbb E}(X)$, we get
\begin{align*}
{\mathbb E}\left[{\cal U} g\left(S_{k_1,k_2}^n\right)\right]=&-wa\sum_{m=0}^{\floor{n/k}}g(m+1)p_{m,n}+a(p)\sum_{i=1}^{3}a_i\sum_{s=0}^{d_i}b_i(n-s)\hspace{-0.2cm}\sum_{l=0}^{\floor{s/k}}B_s(l)\hspace{-0.5cm}\sum_{m=0}^{\floor{\frac{n-k-s-i+1}{k}}}\hspace{-0.6cm}g(m+l+1)p_{m,n-k-s-i+1}\\
&-wb a(p)\sum_{i=1}^{3}a_i\sum_{s=0}^{d_i}b_i(n\hspace{-0.07cm}-\hspace{-0.07cm}s)\hspace{-0.07cm}\sum_{l=0}^{\floor{s/k}}B_s(l)\sum_{m=0}^{\floor{\frac{n-k-s-i+1}{k}}}g(m+l+2)p_{m,n-k-s-i+1}
\end{align*}
Observe that $\floor{\frac{n-k-s-i+1}{k}} \le \floor{n/k}$ for all $s$ and $i$, hence, we replace $\floor{\frac{n-k-s-i+1}{k}}$ by $\floor{n/k}$ as $p_{m,n-k-s-i+1}$ become zero outside of its range. Hence,
\begin{align}
{\mathbb E}\left[{\cal U} g\left(S_{k_1,k_2}^n\right)\right] = &-wa\sum_{m=0}^{\floor{n/k}}g(m+1)p_{m,n}+a(p)\sum_{i=1}^{3}a_i\sum_{s=0}^{d_i}b_i(n-s)\hspace{-0.2cm}\sum_{l=0}^{\floor{s/k}}B_s(l)\hspace{-0.1cm}\sum_{m=0}^{\floor{n/k}}\hspace{-0.1cm}g(m+l+1)p_{m,n-k-s-i+1}\nonumber\\
&-wb a(p)\sum_{i=1}^{3}a_i\sum_{s=0}^{d_i}b_i(n\hspace{-0.07cm}-\hspace{-0.07cm}s)\hspace{-0.07cm}\sum_{l=0}^{\floor{s/k}}B_s(l)\sum_{m=0}^{\floor{n/k}}g(m+l+2)p_{m,n-k-s-i+1}\nonumber\\
= &-wa\sum_{m=0}^{\floor{n/k}}g(m\hspace{-0.05cm}+\hspace{-0.05cm}1)p_{m,n}\hspace{-0.05cm}-\hspace{-0.05cm}wb a(p)\sum_{i=1}^{3}a_i\sum_{s=0}^{d_i}b_i(n\hspace{-0.07cm}-\hspace{-0.07cm}s)\hspace{-0.07cm}\sum_{l=0}^{\floor{s/k}}B_s(l)\sum_{m=0}^{\floor{n/k}}\hspace{-0.1cm}\Delta g(m\hspace{-0.05cm}+\hspace{-0.05cm}l\hspace{-0.05cm}+\hspace{-0.05cm}1)p_{m,n-k-s-i+1}\nonumber\\
&+(1-wb)a(p)\sum_{i=1}^{3}a_i\sum_{s=0}^{d_i}b_i(n-s)\hspace{-0.2cm}\sum_{l=0}^{\floor{s/k}}B_s(l)\hspace{-0.1cm}\sum_{m=0}^{\floor{n/k}}\hspace{-0.1cm}g(m+l+1)p_{m,n-k-s-i+1}.\label{three:exp0}
\end{align}
It is known that
\begin{equation}
g(m+l+1) = \sum_{j=1}^{l}\Delta g(m+j) + g(m+1). \label{three:delta}
\end{equation}
Substituting \eqref{three:delta} in \eqref{three:exp0}, we have
\begin{align}
{\mathbb E}\left[{\cal U} g\left(S_{k_1,k_2}^n\right)\right] &= -wa\sum_{m=0}^{\floor{n/k}}g(m\hspace{-0.05cm}+\hspace{-0.05cm}1)p_{m,n}\hspace{-0.05cm}-\hspace{-0.05cm}wb a(p)\sum_{i=1}^{3}a_i\sum_{s=0}^{d_i}b_i(n\hspace{-0.07cm}-\hspace{-0.07cm}s)\hspace{-0.07cm}\sum_{l=0}^{\floor{s/k}}B_s(l)\sum_{m=0}^{\floor{n/k}}\hspace{-0.1cm}\Delta g(m\hspace{-0.05cm}+\hspace{-0.05cm}l\hspace{-0.05cm}+\hspace{-0.05cm}1)p_{m,n-k-s-i+1}\nonumber\\
&~~~+(1-wb)a(p)\sum_{i=1}^{3}a_i\sum_{s=0}^{d_i}b_i(n-s)\hspace{-0.2cm}\sum_{l=0}^{\floor{s/k}}B_s(l)\hspace{-0.1cm}\sum_{m=0}^{\floor{n/k}}\hspace{-0.1cm}g(m+1)p_{m,n-k-s-i+1}\nonumber\\
&~~~+(1-wb)a(p)\sum_{i=1}^{3}a_i\sum_{s=0}^{d_i}b_i(n-s)\hspace{-0.2cm}\sum_{l=0}^{\floor{s/k}}B_s(l)\hspace{-0.1cm}\sum_{m=0}^{\floor{n/k}}\sum_{j=1}^{l}\Delta g(m+j)p_{m,n-k-s-i+1}.\label{three:zzzz}
\end{align}
Now, from Lemma \ref{three:le1} $(ii)$, replacing $n$ by $n-1,~n-2,\dotsc,n-l+1$, we have
\begin{align*}
p_{m,n}-p_{m,n-1}&=-a(p)[p_{m,n-k}-p_{m-1,n-k}-(p_{m,n-k-1}-p_{m-1,n-k-1})\\
&~~~+qp~(p_{m,n-k-2}-p_{m-1,n-k-2})-(p_{m,n-k}-p_{m-1,n-k}){\bf 1}(n-k=0)\\
&~~~+q(p_{m,n-k-1}-p_{m-1,n-k-1}){\bf 1}(n-k-1=0)]\\
p_{m,n-1}-p_{m,n-2}&=-a(p)[p_{m,n-k-1}-p_{m-1,n-k-1}-(p_{m,n-k-2}-p_{m-1,n-k-2})\\
&~~~+qp~(p_{m,n-k-3}-p_{m-1,n-k-3})-(p_{m,n-k-1}-p_{m-1,n-k-1}){\bf 1}(n-k=1)\\
&~~~+q(p_{m,n-k-2}-p_{m-1,n-k-2}){\bf 1}(n-k-1=1)]\\
&\vdots\\
p_{m,n-\ell+1}-p_{m,n-\ell}&=-a(p)[p_{m,n-k-\ell+1}-p_{m-1,n-k-\ell+1}-(p_{m,n-k-\ell}-p_{m-1,n-k-\ell})\\
&~~~+qp~(p_{m,n-k-\ell-1}-p_{m-1,n-k-\ell-1})\\
&~~~-(p_{m,n-k-\ell+1}-p_{m-1,n-k-\ell+1}){\bf 1}(n-k=\ell-1)\\
&~~~+q(p_{m,n-k-\ell}-p_{m-1,n-k-\ell}){\bf 1}(n-k-1=\ell-1)].
\end{align*}
Adding all the above terms, it can be verified that
\begin{align}
p_{m,n}=p_{m,n-l}-a(p)({p}_{m,n,l}^\star-{p}_{m-1,n,l}^\star),\quad {\rm for}~l \ge 1,\label{three:pmfrela}
\end{align}
where 
\begin{equation}
p_{m,n,l}^\star=p_{m,n-k}-p_{m,n-k-l}+\sum_{u=0}^{l-1}[qp p_{m,n-k-u-2}-p_{m,n-k-u}{\bf 1}(u=n-k)+qp_{m,n-k-u-1}{\bf 1}(u=n-k-1)].\label{three:pmfrela1}
\end{equation}
Using \eqref{three:pmfrela} with $\ell=k+s+i-1$ in the third term of \eqref{three:zzzz}, the expression leads to
\begin{align}
{\mathbb E}\left[{\cal U} g\left(S_{k_1,k_2}^n\right)\right] &= -wa\sum_{m=0}^{\floor{n/k}}g(m\hspace{-0.05cm}+\hspace{-0.05cm}1)p_{m,n}\hspace{-0.05cm}-\hspace{-0.05cm}wb a(p)\sum_{i=1}^{3}a_i\sum_{s=0}^{d_i}b_i(n\hspace{-0.07cm}-\hspace{-0.07cm}s)\hspace{-0.07cm}\sum_{l=0}^{\floor{s/k}}B_s(l)\sum_{m=0}^{\floor{n/k}}\hspace{-0.1cm}\Delta g(m\hspace{-0.05cm}+\hspace{-0.05cm}l\hspace{-0.05cm}+\hspace{-0.05cm}1)p_{m,n-k-s-i+1}\nonumber\\
&~~~+\hspace{-0.07cm}(1\hspace{-0.07cm}-\hspace{-0.07cm}wb)a(p)\sum_{i=1}^{3}\hspace{-0.1cm}a_i\sum_{s=0}^{d_i}\hspace{-0.07cm}b_i(n\hspace{-0.07cm}-\hspace{-0.07cm}s)\hspace{-0.2cm}\sum_{l=0}^{\floor{s/k}}\hspace{-0.15cm}B_s(l)\hspace{-0.2cm}\sum_{m=0}^{\floor{n/k}}\hspace{-0.1cm}g(m\hspace{-0.07cm}+\hspace{-0.07cm}1)[p_{m,n}\hspace{-0.07cm}+\hspace{-0.07cm}a(p)(p_{m,n,k+s+i-1}^\star\hspace{-0.07cm}-\hspace{-0.07cm}p_{m-1,n,k+s+i-1}^\star)]\nonumber\\
&~~~+(1-wb)a(p)\sum_{i=1}^{3}a_i\sum_{s=0}^{d_i}b_i(n-s)\hspace{-0.2cm}\sum_{l=0}^{\floor{s/k}}B_s(l)\hspace{-0.1cm}\sum_{m=0}^{\floor{n/k}}\sum_{j=1}^{l}\Delta g(m+j)p_{m,n-k-s-i+1}\nonumber\allowdisplaybreaks\\
& = (1-wb)\left(-\frac{wa}{1-wb}+a(p)\sum_{i=1}^{3}a_i\sum_{s=0}^{d_i}b_i(n-s)\sum_{l=0}^{\floor{s/k}}B_s(l)\right)\sum_{m=0}^{\floor{n/k}}g(m+1)p_{m,n}\nonumber\\
&~~~+(1\hspace{-0.07cm}-\hspace{-0.07cm}wb)a(p)^2\sum_{i=1}^{3}a_i\sum_{s=0}^{d_i}b_i(n\hspace{-0.07cm}-\hspace{-0.07cm}s)\sum_{l=0}^{\floor{s/k}}B_s(l)\sum_{m=0}^{\floor{n/k}}g(m\hspace{-0.07cm}+\hspace{-0.07cm}1)(p_{m,n,k+s+i-1}^\star\hspace{-0.07cm}-\hspace{-0.07cm}p_{m-1,n,k+s+i-1}^\star)\nonumber\\
&~~~-wb a(p)\sum_{i=1}^{3}a_i\sum_{s=0}^{d_i}b_i(n-s)\sum_{l=0}^{\floor{s/k}}B_s(l)\sum_{m=0}^{\floor{n/k}}\Delta g(m+l+1)p_{m,n-k-s-i+1}\nonumber\allowdisplaybreaks\\
&~~~+(1-wb)a(p)\sum_{i=1}^{3}a_i\sum_{s=0}^{d_i}b_i(n-s)\hspace{-0.2cm}\sum_{l=0}^{\floor{s/k}}B_s(l)\hspace{-0.1cm}\sum_{m=0}^{\floor{n/k}}\sum_{j=1}^{l}\Delta g(m+j)p_{m,n-k-s-i+1}. \label{three:exp2}
\end{align}
Using the fact that $\sum_{l=0}^{\floor{s/k}}B_s(l)=C_s(1)=\frac{(k+1)^s}{(k+2)^{s+1}}$ with ${\mathbb E}(Z)={\mathbb E}\big(S_{k_1,k_2}^n\big)$, it can be verified that
\begin{align*}
-\frac{wa}{1-wb} +a(p)\sum_{i=1}^{3}a_i\sum_{s=0}^{d_i}b_i(n-s)\sum_{l=0}^{\floor{s/k}}B_s(l)&=-\frac{wa}{1-wb}+q[1+(n-k-1)p]a(p)\\
&=-{\mathbb E}(Z)+{\mathbb E}\big(S_{k_1,k_2}^n\big)=0.
\end{align*}
The expression \eqref{three:exp2} now becomes
\begin{align}
{\mathbb E}\left[{\cal U} g\left(S_{k_1,k_2}^n\right)\right] &=(1\hspace{-0.07cm}-\hspace{-0.07cm}wb)a(p)^2\sum_{i=1}^{3}a_i\sum_{s=0}^{d_i}b_i(n\hspace{-0.07cm}-\hspace{-0.07cm}s)\sum_{l=0}^{\floor{s/k}}B_s(l)\sum_{m=0}^{\floor{n/k}}g(m\hspace{-0.07cm}+\hspace{-0.07cm}1)(p_{m,n,k+s+i-1}^\star\hspace{-0.07cm}-\hspace{-0.07cm}p_{m-1,n,k+s+i-1}^\star)\nonumber\\
&~~~-wb a(p)\sum_{i=1}^{3}a_i\sum_{s=0}^{d_i}b_i(n-s)\sum_{l=0}^{\floor{s/k}}B_s(l)\sum_{m=0}^{\floor{n/k}}\Delta g(m+l+1)p_{m,n-k-s-i+1}\nonumber\allowdisplaybreaks\\
&~~~+(1-wb)a(p)\sum_{i=1}^{3}a_i\sum_{s=0}^{d_i}b_i(n-s)\hspace{-0.2cm}\sum_{l=0}^{\floor{s/k}}B_s(l)\hspace{-0.1cm}\sum_{m=0}^{\floor{n/k}}\sum_{j=1}^{l}\Delta g(m+j)p_{m,n-k-s-i+1}\nonumber\\
& =-(1-wb)a(p)^2\sum_{i=1}^{3}a_i\sum_{s=0}^{d_i}b_i(n-s)\sum_{l=0}^{\floor{s/k}}B_s(l)\sum_{m=0}^{\floor{n/k}}\Delta g(m+1)p_{m,n,k+s+i-1}^\star\nonumber\\
&~~~~-wb a(p)\sum_{i=1}^{3}a_i\sum_{s=0}^{d_i}b_i(n-s)\sum_{l=0}^{\floor{s/k}}B_s(l)\sum_{m=0}^{\floor{n/k}}\Delta g(m+l+1)p_{m,n-k-s-i+1}\nonumber\\
&~~~+(1-wb)a(p)\sum_{i=1}^{3}a_i\sum_{s=0}^{d_i}b_i(n-s)\sum_{l=0}^{\floor{s/k}}B_s(l)\hspace{-0.1cm}\sum_{m=0}^{\floor{n/k}}\sum_{j=1}^{l}\Delta g(m+j)p_{m,n-k-s-i+1}.\label{three:oneexp}
\end{align}
We know that
\begin{equation}
\Delta g(m+j) = \sum_{v=1}^{j-1}\Delta^2 g(m+v)+\Delta g(m+1). \label{three:ddelta}
\end{equation}
Substituting \eqref{three:ddelta} in \eqref{three:oneexp}, we get
\begin{align}
{\mathbb E}\left[{\cal U} g\left(S_{k_1,k_2}^n\right)\right] & =-(1-wb)a(p)^2\sum_{i=1}^{3}a_i\sum_{s=0}^{d_i}b_i(n-s)\sum_{l=0}^{\floor{s/k}}B_s(l)\sum_{m=0}^{\floor{n/k}}\Delta g(m+1)p_{m,n,k+s+i-1}^\star\nonumber\\
&~~~-wb a(p)\sum_{i=1}^{3}a_i\sum_{s=0}^{d_i}b_i(n-s)\sum_{l=0}^{\floor{s/k}}B_s(l)\sum_{m=0}^{\floor{n/k}}\Delta g(m+1)p_{m,n-k-s-i+1}\nonumber\\
&~~~-wb a(p)\sum_{i=1}^{3}a_i\sum_{s=0}^{d_i}b_i(n-s)\sum_{l=0}^{\floor{s/k}}B_s(l)\sum_{m=0}^{\floor{n/k}}\sum_{j=1}^{l}\Delta^2 g(m+j)p_{m,n-k-s-i+1}\nonumber\allowdisplaybreaks\\
&~~~+(1-wb)a(p)\sum_{i=1}^{3}a_i\sum_{s=0}^{d_i}b_i(n-s)\sum_{l=0}^{\floor{s/k}}lB_s(l)\sum_{m=0}^{\floor{n/k}}\Delta g(m+1)p_{m,n-k-s-i+1}\nonumber\\
&~~~+(1-wb)a(p)\sum_{i=1}^{3}a_i\sum_{s=0}^{d_i}b_i(n-s)\sum_{l=0}^{\floor{s/k}}B_s(l)\sum_{m=0}^{\floor{n/k}}\sum_{j=1}^{l}\sum_{v=1}^{j-1}\Delta^2 g(m+v)p_{m,n-k-s-i+1}.\label{mmmm}
\end{align}
From \eqref{three:pmfrela1}, for $\ell=k+s+i-1$, we have
\begin{align*}
p_{m,n,k+s+i-1}^\star\hspace{-0.08cm}&=\hspace{-0.08cm}p_{m,n-k}\hspace{-0.08cm}-\hspace{-0.08cm}p_{m,n-2k-s-i+1}\hspace{-0.08cm}+\hspace{-0.3cm}\sum_{u=0}^{k+s+i-2}\hspace{-0.25cm}[qp p_{m,n-k-u-2}\hspace{-0.08cm}-\hspace{-0.08cm}p_{m,n-k-u}{\bf 1}(u\hspace{-0.08cm}=\hspace{-0.08cm}n\hspace{-0.08cm}-\hspace{-0.08cm}k)\hspace{-0.08cm}+\hspace{-0.08cm}qp_{m,n-k-u-1}{\bf 1}(u\hspace{-0.08cm}=\hspace{-0.08cm}n\hspace{-0.08cm}-\hspace{-0.08cm}k\hspace{-0.08cm}-\hspace{-0.08cm}1)]
\end{align*}
Using \eqref{three:pmfrela}, for the right hand side terms, we get
\begin{align}
p_{m,n,k+s+i-1}^\star&=p_{m,n}+a(p)[p_{m,n,k}^\star-p_{m-1,n,k}^\star]\nonumber\\
&~~~-\left[p_{m,n}{\bf 1}(s\le n-2k-i+1)+a(p)(p_{m,n,2k+s+i-1}^\star-p_{m-1,n,2k+s+i-1}^\star)\right]\nonumber\\
&~~~+\sum_{u=0}^{k+s+i-2}\Bigr[qp( p_{m,n}{\bf 1}(u \le n-k-2)+a(p)(p_{m,n,k+u+2}^\star-p_{m-1,n,k+u+2}^\star))\nonumber\\
&~~~~~~~~~~~~~~~~~~~~~~-(p_{m,n}{\bf 1}(u=n-k)+a(p)(p_{m,n,k+u}^\star-p_{m-1,n,k+u}^\star){\bf 1}(u=n-k))\nonumber\\
&~~~~~~~~~~~~~~~~~~~~~~+q(p_{m,n}{\bf 1}(u=n-k-1)+a(p)(p_{m,n,k+u+1}^\star-p_{m-1,n,k+u+1}^\star){\bf 1}(u=n-k-1))\Bigr]\nonumber\allowdisplaybreaks\\
&=\left[1\hspace{-0.08cm}-\hspace{-0.08cm}{\bf 1}(s \hspace{-0.08cm}\le\hspace{-0.08cm} n\hspace{-0.08cm}-\hspace{-0.08cm}2k\hspace{-0.08cm}-\hspace{-0.08cm}i\hspace{-0.08cm}+\hspace{-0.08cm}1)\hspace{-0.08cm}+\hspace{-0.08cm}\sum_{u=0}^{k+s+i-2}(qp {\bf 1}(u \hspace{-0.08cm}\le\hspace{-0.08cm} n\hspace{-0.08cm}-\hspace{-0.08cm}k\hspace{-0.08cm}-\hspace{-0.08cm}2)\hspace{-0.08cm}-\hspace{-0.08cm}{\bf 1}(u\hspace{-0.08cm}=\hspace{-0.08cm}n\hspace{-0.08cm}-\hspace{-0.08cm}k)\hspace{-0.08cm}+\hspace{-0.08cm}q{\bf 1}(u\hspace{-0.08cm}=\hspace{-0.08cm}n\hspace{-0.08cm}-\hspace{-0.08cm}k\hspace{-0.08cm}-\hspace{-0.08cm}1))\right]p_{m,n}\nonumber\\
&~~~+a(p)\Bigr[p_{m,n,k}^\star-p_{m,n,2k+s+i-1}^\star+\sum_{u=0}^{k+s+i-2}\Bigr\{qp p_{m,n,k+u+2}^\star-p_{m,n,k+u}^\star{\bf 1}(u=n-k)\nonumber\\
&~~~+q p_{m,n,k+u+1}^\star{\bf 1}(u=n-k-1)\Bigr\}\Bigr]-a(p)\Bigr[p_{m-1,n,k}^\star-p_{m-1,n,2k+s+i-1}^\star\nonumber\\
&~~~+\sum_{u=0}^{k+s+i-2}\Bigr\{qp p_{m-1,n,k+u+2}^\star-p_{m-1,n,k+u}^\star{\bf 1}(u=n-k)+q p_{m-1,n,k+u+1}^\star{\bf 1}(u=n-k-1)\Bigr\}\Bigr]\nonumber\\
&=\left[1\hspace{-0.08cm}-\hspace{-0.08cm}{\bf 1}(s \hspace{-0.08cm}\le\hspace{-0.08cm} n\hspace{-0.08cm}-\hspace{-0.08cm}2k\hspace{-0.08cm}-\hspace{-0.08cm}i\hspace{-0.08cm}+\hspace{-0.08cm}1)\hspace{-0.08cm}+\hspace{-0.08cm}\sum_{u=0}^{k+s+i-2}(qp {\bf 1}(u \hspace{-0.08cm}\le\hspace{-0.08cm} n\hspace{-0.08cm}-\hspace{-0.08cm}k\hspace{-0.08cm}-\hspace{-0.08cm}2)\hspace{-0.08cm}-\hspace{-0.08cm}{\bf 1}(u\hspace{-0.08cm}=\hspace{-0.08cm}n\hspace{-0.08cm}-\hspace{-0.08cm}k)\hspace{-0.08cm}+\hspace{-0.08cm}q{\bf 1}(u\hspace{-0.08cm}=\hspace{-0.08cm}n\hspace{-0.08cm}-\hspace{-0.08cm}k\hspace{-0.08cm}-\hspace{-0.08cm}1))\right]p_{m,n}\nonumber\\
&~~~+a(p)[p_{m,n,k+s+i-1}^{\star\star}-p_{m-1,n,k+s+i-1}^{\star\star}],\label{three:cccc}
\end{align}
where $p_{m,n,l}^{\star\star}=p_{m,n,k}^\star-p_{m,n,k+l}^\star+\sum_{u=0}^{l-1}[qp p_{m,n,k+u+2}^\star-p_{m,n,k+u}^\star{\bf 1}(u=n-k)+qp_{m,n,k+u+1}^\star{\bf 1}(u=n-k-1)].$ Also, observe that we are not using indicator function other than $p_{m,n}$ term since other terms involving $s$ where the PMF will become zero outside of its range according to the value of $s$. Using \eqref{three:pmfrela} and \eqref{three:cccc} in \eqref{mmmm}, we have
\begin{align}
{\mathbb E}\left[{\cal U} g\left(S_{k_1,k_2}^n\right)\right] &= \left\{-(1-wb)a(p)^2\sum_{i=1}^{3}a_i\sum_{s=0}^{d_i}b_i(n-s)\sum_{l=0}^{\floor{s/k}}B_s(l)\Bigr[1-{\bf 1}(s \le n-2k-i+1)\right.\nonumber\\
&~~~\left.+\sum_{u=0}^{k+s+i-2}(qp {\bf 1}(u \le n-k-2)-{\bf 1}(u=n-k)+q{\bf 1}(u=n-k-1))\Bigr]\right\}\sum_{m=0}^{\floor{n/k}}\Delta g(m+1)p_{m,n}\nonumber\\
&~~~-(1\hspace{-0.07cm}-\hspace{-0.07cm}wb)a(p)^3\sum_{i=1}^{3}a_i\sum_{s=0}^{d_i}b_i(n\hspace{-0.07cm}-\hspace{-0.07cm}s)\sum_{l=0}^{\floor{s/k}}B_s(l)\sum_{m=0}^{\floor{n/k}}\Delta g(m\hspace{-0.07cm}+\hspace{-0.07cm}1)(p_{m,n,k+s+i-1}^{\star\star}\hspace{-0.07cm}-\hspace{-0.07cm}p_{m-1,n,k+s+i-1}^{\star\star})\nonumber\\
&~~~-wba(p)\sum_{i=1}^{3}a_i\sum_{s=0}^{d_i}b_i(n-s)\sum_{l=0}^{\floor{s/k}} B_s(l)\sum_{m=0}^{\floor{n/k}}\Delta g(m+1)p_{m,n}\nonumber\allowdisplaybreaks\\
&~~~-wba(p)^2\sum_{i=1}^{3}a_i\sum_{s=0}^{d_i}b_i(n-s)\sum_{l=0}^{\floor{s/k}} B_s(l)\sum_{m=0}^{\floor{n/k}}\Delta g(m+1)(p_{m,n,k+s+i-1}^\star-p_{m-1,n,k+s+i-1}^\star)\nonumber\\
&~~~+(1-wb) a(p)\sum_{i=1}^{3}a_i\sum_{s=0}^{d_i}b_i(n-s)\sum_{l=0}^{\floor{s/k}}lB_s(l)\sum_{m=0}^{\floor{n/k}}\Delta g(m+1)p_{m,n}\nonumber\\
&~~~+(1\hspace{-0.07cm}-\hspace{-0.07cm}wb)a(p)^2\sum_{i=1}^{3}a_i\sum_{s=0}^{d_i}b_i(n\hspace{-0.07cm}-\hspace{-0.07cm}s)\sum_{l=0}^{\floor{s/k}}lB_s(l)\hspace{-0.1cm}\sum_{m=0}^{\floor{n/k}}\Delta g(m\hspace{-0.07cm}+\hspace{-0.07cm}1)(p_{m,n,k+s+i-1}^\star\hspace{-0.07cm}-\hspace{-0.07cm}p_{m-1,n,k+s+i-1}^\star)\nonumber\allowdisplaybreaks\\
&~~~-wb a(p)\sum_{i=1}^{3}a_i\sum_{s=0}^{d_i}b_i(n-s)\sum_{l=0}^{\floor{s/k}}B_s(l)\sum_{m=0}^{\floor{n/k}}\sum_{j=1}^{l}\Delta^2 g(m+j)p_{m,n-k-s-i+1}\nonumber\\
&~~~+(1-wb)a(p)\sum_{i=1}^{3}a_i\sum_{s=0}^{d_i}b_i(n-s)\sum_{l=0}^{\floor{s/k}}B_s(l)\sum_{m=0}^{\floor{n/k}}\sum_{j=1}^{l}\sum_{v=1}^{j-1}\Delta^2 g(m+v)p_{m,n-k-s-i+1}\nonumber\\
&= \left\{-(1-wb)a(p)^2\sum_{i=1}^{3}a_i\sum_{s=0}^{d_i}b_i(n-s)\sum_{l=0}^{\floor{s/k}}B_s(l)\Bigr[1-{\bf 1}(s \le n-2k-i+1)\right.\nonumber\\
&~~~+\sum_{u=0}^{k+s+i-2}(qp {\bf 1}(u \le n-k-2)-{\bf 1}(u=n-k)+q{\bf 1}(u=n-k-1))\Bigr]\nonumber\\
&~~~-wba(p)\sum_{i=1}^{3}a_i\sum_{s=0}^{d_i}b_i(n-s)\sum_{l=0}^{\floor{s/k}} B_s(l)\nonumber\\
&~~~\left.+(1-wb) a(p)\sum_{i=1}^{3}a_i\sum_{s=0}^{d_i}b_i(n-s)\sum_{l=0}^{\floor{s/k}}lB_s(l)\right\}\sum_{m=0}^{\floor{n/k}}\Delta g(m+1)p_{m,n}\nonumber\allowdisplaybreaks\\
&~~~-(1\hspace{-0.07cm}-\hspace{-0.07cm}wb)a(p)^3\sum_{i=1}^{3}a_i\sum_{s=0}^{d_i}b_i(n\hspace{-0.07cm}-\hspace{-0.07cm}s)\sum_{l=0}^{\floor{s/k}}B_s(l)\sum_{m=0}^{\floor{n/k}}\Delta g(m\hspace{-0.07cm}+\hspace{-0.07cm}1)(p_{m,n,k+s+i-1}^{\star\star}\hspace{-0.07cm}-\hspace{-0.07cm}p_{m-1,n,k+s+i-1}^{\star\star})\nonumber\\
&~~~-wba(p)^2\sum_{i=1}^{3}a_i\sum_{s=0}^{d_i}b_i(n-s)\sum_{l=0}^{\floor{s/k}} B_s(l)\sum_{m=0}^{\floor{n/k}}\Delta g(m+1)(p_{m,n,k+s+i-1}^\star-p_{m-1,n,k+s+i-1}^\star)\nonumber\\
&~~~+(1\hspace{-0.07cm}-\hspace{-0.07cm}wb)a(p)^2\sum_{i=1}^{3}a_i\sum_{s=0}^{d_i}b_i(n\hspace{-0.07cm}-\hspace{-0.07cm}s)\sum_{l=0}^{\floor{s/k}}lB_s(l)\hspace{-0.1cm}\sum_{m=0}^{\floor{n/k}}\Delta g(m\hspace{-0.07cm}+\hspace{-0.07cm}1)(p_{m,n,k+s+i-1}^\star\hspace{-0.07cm}-\hspace{-0.07cm}p_{m-1,n,k+s+i-1}^\star)\nonumber\allowdisplaybreaks\\
&~~~-wb a(p)\sum_{i=1}^{3}a_i\sum_{s=0}^{d_i}b_i(n-s)\sum_{l=0}^{\floor{s/k}}B_s(l)\sum_{m=0}^{\floor{n/k}}\sum_{j=1}^{l}\Delta^2 g(m+j)p_{m,n-k-s-i+1}\nonumber\\
&~~~+(1-wb)a(p)\sum_{i=1}^{3}a_i\sum_{s=0}^{d_i}b_i(n-s)\sum_{l=0}^{\floor{s/k}}B_s(l)\sum_{m=0}^{\floor{n/k}}\sum_{j=1}^{l}\sum_{v=1}^{j-1}\Delta^2 g(m+v)p_{m,n-k-s-i+1}.\label{three:bigexp}
\end{align}
Using $\varphi=\mathrm{Var}(Z)-\mathrm{Var}\big(S_{k_1,k_2}^n\big)$ and ${\mathbb E}(Z)={\mathbb E}\big(S_{k_1,k_2}^n\big)$ with some algebraic calculations, it can be verified that
\begin{align*}
&-(1-wb)a(p)^2\sum_{i=1}^{3}a_i\sum_{s=0}^{d_i}b_i(n-s)\sum_{l=0}^{\floor{s/k}}B_s(l)\Bigr[1-{\bf 1}(s \le n-2k-i+1)\\
&+\sum_{u=0}^{k+s+i-2}(qp {\bf 1}(u \le n-k-2)-{\bf 1}(u=n-k)+q{\bf 1}(u=n-k-1))\Bigr]\\
&-wba(p)\sum_{i=1}^{3}a_i\sum_{s=0}^{d_i}b_i(n-s)\sum_{l=0}^{\floor{s/k}} B_s(l)+(1-wb) a(p)\sum_{i=1}^{3}a_i\sum_{s=0}^{d_i}b_i(n-s)\sum_{l=0}^{\floor{s/k}}lB_s(l)=-\varphi(1-wb).
\end{align*}
In particular, the above expression is zero if $\mathrm{Var}(Z)=\mathrm{Var}\big(S_{k_1,k_2}^n\big)$. Therefore, from \eqref{three:bigexp},
\vspace{-0.001cm}
\begin{align*}
{\mathbb E}\left[{\cal U} g\left(S_{k_1,k_2}^n\right)\right] &=-(1\hspace{-0.07cm}-\hspace{-0.07cm}wb)a(p)^3\sum_{i=1}^{3}a_i\sum_{s=0}^{d_i}b_i(n\hspace{-0.07cm}-\hspace{-0.07cm}s)\sum_{l=0}^{\floor{s/k}}B_s(l)\sum_{m=0}^{\floor{n/k}}\Delta g(m\hspace{-0.07cm}+\hspace{-0.07cm}1)(p_{m,n,k+s+i-1}^{\star\star}\hspace{-0.07cm}-\hspace{-0.07cm}p_{m-1,n,k+s+i-1}^{\star\star})\\
&~~~-wba(p)^2\sum_{i=1}^{3}a_i\sum_{s=0}^{d_i}b_i(n-s)\sum_{l=0}^{\floor{s/k}} B_s(l)\sum_{m=0}^{\floor{n/k}}\Delta g(m+1)(p_{m,n,k+s+i-1}^\star-p_{m-1,n,k+s+i-1}^\star)\\
&~~~+(1\hspace{-0.07cm}-\hspace{-0.07cm}wb)a(p)^2\sum_{i=1}^{3}a_i\sum_{s=0}^{d_i}b_i(n\hspace{-0.07cm}-\hspace{-0.07cm}s)\sum_{l=0}^{\floor{s/k}}lB_s(l)\hspace{-0.1cm}\sum_{m=0}^{\floor{n/k}}\Delta g(m\hspace{-0.07cm}+\hspace{-0.07cm}1)(p_{m,n,k+s+i-1}^\star\hspace{-0.07cm}-\hspace{-0.07cm}p_{m-1,n,k+s+i-1}^\star)\\
&~~~-\hspace{-0.08cm}wb a(p)\hspace{-0.08cm}\sum_{i=1}^{3}a_i\hspace{-0.08cm}\sum_{s=0}^{d_i}b_i(n\hspace{-0.08cm}-\hspace{-0.08cm}s)\hspace{-0.18cm}\sum_{l=0}^{\floor{s/k}}\hspace{-0.08cm}B_s(l)\hspace{-0.18cm}\sum_{m=0}^{\floor{n/k}}\hspace{-0.08cm}\sum_{j=1}^{l}\hspace{-0.08cm}\hspace{-0.08cm}\Delta^2 g(m\hspace{-0.08cm}+\hspace{-0.08cm}j)p_{m,n-k-s-i+1}\hspace{-0.08cm}-\hspace{-0.08cm}\varphi(1\hspace{-0.08cm}-\hspace{-0.08cm}wb)\hspace{-0.18cm}\sum_{m=0}^{\floor{n/k}}\hspace{-0.15cm}\Delta g(m\hspace{-0.08cm}+\hspace{-0.08cm}1)p_{m,n}\\
&~~~+(1-wb)a(p)\sum_{i=1}^{3}a_i\sum_{s=0}^{d_i}b_i(n-s)\sum_{l=0}^{\floor{s/k}}B_s(l)\sum_{m=0}^{\floor{n/k}}\sum_{j=1}^{l}\sum_{v=1}^{j-1}\Delta^2 g(m+v)p_{m,n-k-s-i+1}.
\end{align*}
Observe that
\begin{align}
d_i \le n-k-1\quad {\rm and}\quad |b_i(n-s)| \le n-s+1, \quad {\rm for~all}~s,i. \label{three:inq}
\end{align}
Now, observe that
\begin{align}
\left|{\mathbb E}\left[{\cal U} g\left(S_{k_1,k_2}^n\right)\right]\right| \hspace{-0.07cm}&\le\hspace{-0.07cm} \|\Delta g\|\left\{|1\hspace{-0.07cm}-\hspace{-0.07cm}wb|a(p)^3\sum_{i=1}^{3}|a_i|\sum_{s=0}^{n-k-1}(n\hspace{-0.07cm}-\hspace{-0.07cm}s\hspace{-0.07cm}+\hspace{-0.07cm}1)\sum_{l=0}^{\floor{s/k}}B_s(l)\sum_{m=0}^{\floor{n/k}}|p_{m,n,k+s+i-1}^{\star\star}\hspace{-0.07cm}-\hspace{-0.07cm}p_{m-1,n,k+s+i-1}^{\star\star}|\right.\nonumber\\
&~~~+|wb|a(p)^2\sum_{i=1}^{3}|a_i|\sum_{s=0}^{n-k-1}(n-s+1)\sum_{l=0}^{\floor{s/k}} B_s(l)\sum_{m=0}^{\floor{n/k}}|p_{m,n,k+s+i-1}^\star-p_{m-1,n,k+s+i-1}^\star|\nonumber\allowdisplaybreaks\\
&~~~+|1\hspace{-0.07cm}-\hspace{-0.07cm}wb|a(p)^2\sum_{i=1}^{3}|a_i|\sum_{s=k}^{n-k-1}(n-s+1)\sum_{l=0}^{\floor{s/k}}lB_s(l)\hspace{-0.1cm}\sum_{m=0}^{\floor{n/k}}|p_{m,n,k+s+i-1}^\star\hspace{-0.07cm}-\hspace{-0.07cm}p_{m-1,n,k+s+i-1}^\star|\nonumber\\
&~~~+\hspace{-0.08cm}|wb| a(p)\sum_{i=1}^{3}|a_i|\hspace{-0.08cm}\sum_{s=k}^{n-k-1}(n-s+1)\hspace{-0.08cm}\sum_{l=0}^{\floor{s/k}}lB_s(l)\hspace{-0.08cm}\sum_{m=0}^{\floor{n/k}}\hspace{-0.08cm}|p_{m-1,n-k-s-i+1}-p_{m,n,n-k-s-i+1}|\nonumber\\
&~~~\left.+\hspace{-0.07cm}|1\hspace{-0.08cm}-\hspace{-0.08cm}wb|a(p)\hspace{-0.07cm}\sum_{i=1}^{3}|a_i|\hspace{-0.08cm}\sum_{s=2k}^{n-k-1}\hspace{-0.1cm}(n\hspace{-0.08cm}-\hspace{-0.08cm}s\hspace{-0.08cm}+\hspace{-0.08cm}1)\hspace{-0.08cm}\sum_{l=0}^{\floor{s/k}}\frac{l(l\hspace{-0.07cm}-\hspace{-0.07cm}1)}{2}B_s(l)\sum_{m=0}^{\floor{n/k}}|p_{m-1,n-k-s-i+1}\hspace{-0.08cm}-\hspace{-0.08cm}p_{m,n,n-k-s-i+1}|\hspace{-0.07cm}\right\}\nonumber\\
&~~~+\|\Delta g\||\varphi(1-wb)|\sum_{m=0}^{\floor{n/k}}p_{m,n}.\label{bvbv}
\end{align}
Note that
\begin{align}
\sum_{m=0}^{\floor{n/k}}|p_{m,n-l}-p_{m-1,n-l}|=2 d_{TV}(S_{k_1,k_2}^{n-l},S_{k_1,k_2}^{n-l}+1)\le 2\label{11}
\end{align}
with different values of $l \ge 0$. Therefore,
\begin{align}
\sum_{m=0}^{\floor{n/k}}|p^\star_{m,n,l}-p^\star_{m-1,n,l}|&=\sum_{m=0}^{\floor{n/k}}\Bigr| (p_{m,n-k}-p_{m-1,n-k})-(p_{m,n-k-l}-p_{m-1,n-k-l})\nonumber\\
&~~~+\sum_{u=0}^{l-1}\Bigr[qp(p_{m,n-k-u-2}-p_{m,n-k-u-2})-(p_{m,n-k-u}-p_{m,n-k-u}){\bf 1}(u=n-k)\nonumber\\
&~~~~~~~~~~~~~~+q(p_{m,n-k-u-1}-p_{m,n-k-u-1}){\bf 1}(u=n-k-1)\Bigr]\Bigr|\le 2\Bigr[2+(1+q+qp)l\Bigr].\label{zaq}
\end{align}
Also, note that 
\vspace{-0.175cm}
\begin{align}
C_s(1)&=\sum_{l=0}^{\floor{s/k}} B_s(l)=\frac{(k+1)^s}{(k+2)^{s+1}},~~C_s^\prime(1)=\sum_{l=0}^{\floor{s/k}}lB_s(l)=a(p)((s-k+2)k+2)\frac{(k+1)^{s-k-1}}{(k+2)^{s-k+2}}\label{12}\\
C_{s}^{\prime\prime}(1)&=\sum_{l=0}^{\floor{s/k}}\frac{l(l-1)}{2}B_s(l)=a(p)^2(8+k(k(s-2k-1)(s-2k+8)+8(s+2))) \frac{(k+1)^{s-2k-2}}{(k+2)^{s-2k+3}}\label{13}
\end{align}
Next, using \eqref{zaq} and \eqref{12} in the second term of \eqref{bvbv}
\vspace{-0.175cm}
\begin{align}
&|wb|a(p)^2\sum_{i=1}^{3}|a_i|\sum_{s=0}^{n-k-1}(n-s+1)\sum_{l=0}^{\floor{s/k}} B_s(l)\sum_{m=0}^{\floor{n/k}}|p_{m,n,k+s+i-1}^\star-p_{m-1,n,k+s+i-1}^\star|\nonumber\\
&\le (2+qp)|wb|a(p)^2\sum_{s=0}^{n-k-1}(n-s+1) \frac{(k+1)^s}{(k+2)^{s+1}} (2+(1+q+qp)(s+k+2))\times 2\nonumber\\
&=2(2\hspace{-0.08cm}+\hspace{-0.08cm}qp)|wb|a(p)^2\left[\sum_{s=0}^{n-k-1}(n\hspace{-0.08cm}-\hspace{-0.08cm}s\hspace{-0.08cm}+\hspace{-0.08cm}1) \frac{(k+1)^s}{(k+2)^{s+1}} (2\hspace{-0.08cm}+\hspace{-0.08cm}(1\hspace{-0.08cm}+\hspace{-0.08cm}q\hspace{-0.08cm}+\hspace{-0.08cm}qp)(k\hspace{-0.08cm}+\hspace{-0.08cm}2))+(1\hspace{-0.08cm}+\hspace{-0.08cm}q\hspace{-0.08cm}+\hspace{-0.08cm}qp)\sum_{s=0}^{n-k-1}s(n\hspace{-0.08cm}-\hspace{-0.08cm}s\hspace{-0.08cm}+\hspace{-0.08cm}1) \frac{(k+1)^s}{(k+2)^{s+1}})\right]\nonumber\\
&=2(2+qp)|wb|a(p)^2\left[(n-k)(2+(1+q+qp)(k+2))+(1+q+qp)\left((n-2k-2)(k+1)-\frac{(k+1)^{n-k+1}}{(k+2)^{n-k-1}}\right)\right]\nonumber\\
&=2(2+qp)|wb|a(p)^2\left[(n-k)\delta+\delta^* c_{n,k}^{(1)}\right].\label{25}
\end{align}
Now, using \eqref{11} and \eqref{12} in the forth term of \eqref{bvbv}
\vspace{-0.175cm}
\begin{align}
&|wb| a(p)\sum_{i=1}^{3}|a_i|\sum_{s=k}^{n-k-1}(n-s+1)\sum_{l=0}^{\floor{s/k}}l B_s(l)\sum_{m=0}^{\floor{n/k}}\hspace{-0.08cm}|p_{m-1,n-k-s-i+1}-p_{m,n,n-k-s-i+1}|\nonumber\\
&\le(2+qp)|wb| a(p)^2\sum_{s=k}^{n-k-1}(n-s+1)((s-k+2)k+2)\frac{(k+1)^{s-k-1}}{(k+2)^{s-k+2}}\times 2\nonumber\\
&= 2(2+qp)|wb| a(p)^2\sum_{s=k}^{n-k-1}(n-s+1)((s-k+2)k+2)\frac{(k+1)^{s-k-1}}{(k+2)^{s-k+2}}\nonumber\\
&=2(2+qp)|wb| a(p)^2\left[n-3k +k\left(\frac{k+1}{k+2}\right)^{n-2k}\right]=2(2+qp)|wb| a(p)^2 c_{n,k}^{(2)}.\label{24}
\end{align}
Now, using \eqref{zaq} and \eqref{12} in the third term of \eqref{bvbv}
\vspace{-0.175cm}
\begin{align}
&|1\hspace{-0.07cm}-\hspace{-0.07cm}wb|a(p)^2\sum_{i=1}^{3}|a_i|\sum_{s=k}^{n-k-1}(n-s+1)\sum_{l=0}^{\floor{s/k}}lB_s(l)\hspace{-0.1cm}\sum_{m=0}^{\floor{n/k}}|p_{m,n,k+s+i-1}^\star\hspace{-0.07cm}-\hspace{-0.07cm}p_{m-1,n,k+s+i-1}^\star|\nonumber\\
&\le (2+qp)|1\hspace{-0.07cm}-\hspace{-0.07cm}wb|a(p)^3\sum_{s=k}^{n-k-1}(n-s+1)((s-k+2)k+2)\frac{(k+1)^{s-k-1}}{(k+2)^{s-k+2}}\times (2+(1+q+qp)(s+k+2))\times 2\nonumber\\
&=2(2+qp)|1-wb|a(p)^3\left[(2+(1+q+qp)(k+2))\sum_{s=k}^{n-k-1}(n-s+1)((s-k+2)k+2)\frac{(k+1)^{s-k-1}}{(k+2)^{s-k+2}}\right.\nonumber\\
&~~~ \left.~~~~~~~~~~~~~~~~~~~~~~~~~~~~~~~+(1+q+qp)\sum_{s=k}^{n-k-1}s(n-s+1)((s-k+2)k+2)\frac{(k+1)^{s-k-1}}{(k+2)^{s-k+2}}\right]\nonumber\\
&=2(2+qp)|1-wb|a(p)^3\left[(2+(1+q+qp)(k+2))\left(n-3k +k\left(\frac{k+1}{k+2}\right)^{n-2k}\right)\right.\nonumber\\
&~~~ \left.~~~~~~~~~~~~~~~~~~~~~~~~~~~~~~~+(1+q+qp)\left(n(3k+1)-(11 k^2+9k+2)+(2nk+k^2+7k+2)\left(\frac{k+1}{k+2}\right)^{n-2k}\right)\right]\nonumber\\
&=2(2+qp)|1-wb|a(p)^3\left[\delta \times c_{n,k}^{(2)}+\delta^*\times c_{n,k}^{(4)}\right].\label{23}
\end{align}
Now, using \eqref{11} and \eqref{13} in the fifth term of \eqref{bvbv}
\begin{align}
&|1-wb|a(p)\sum_{i=1}^{3}|a_i|\sum_{s=2k}^{n-k-1}(n-s+1)\sum_{l=0}^{\floor{s/k}}\frac{l(l-1)}{2}B_s(l)\sum_{m=0}^{\floor{n/k}}|p_{m-1,n-k-s-i+1}-p_{m,n,n-k-s-i+1}|\nonumber\\
&\le(2+qp)|1-wb|a(p)^3 \sum_{s=2k}^{n-k-1}(n-s+1)(8+k(k(s-2k-1)(s-2k+8)+8(s+2))) \frac{(k+1)^{s-2k-2}}{(k+2)^{s-2k+3}}\times 2\nonumber\\
&=2(2+qp)|1-wb|a(p)^3\left[n-5k+k(nk+6k+4-k^2)\frac{(k+1)^{n-3k-1}}{(k+2)^{n-3k+1}}\right]=2(2+qp)|1-wb|a(p)^3c_{n,k}^{(3)}.\label{22}
\end{align}
Next, using the definition of $p_{m,n,\ell}^{\star\star}$ with \eqref{zaq}, it can be verified that
\begin{align}
\sum_{m=0}^{\floor{n/k}}|p_{m,n,\ell}^{\star\star}\hspace{-0.07cm}-\hspace{-0.07cm}p_{m-1,n,\ell}^{\star\star}|&=\sum_{m=0}^{\floor{n/k}}\Bigr|(p_{m,n,k}^\star-p_{m-1,n,k}^\star)-(p_{m,n,k+\ell}^\star-p_{m-1,n,k+\ell}^\star)\nonumber\\
&~~~+\sum_{u=0}^{\ell-1}[qp(p_{m,n,k+u+2}^\star-p_{m-1,n,k+u+2}^\star)\nonumber\\
&~~~~~~~~~~~~~~~~~~~~~-(p_{m,n,k+u}^\star-p_{m-1,n,k+u}^\star){\bf 1}(u=n-k)\nonumber\\
&~~~~~~~~~~~~~~~~~~~~~+q(p_{m,n,k+u+1}^\star-p_{m-1,n,k+u+1}^\star){\bf 1}(u=n-k-1)]\Bigr|\nonumber\allowdisplaybreaks\\
&\le 2\Bigr[2+(1+q+qp)k+2+(1+q+qp)(k+\ell)\nonumber\\
&~~~~~~~~~+(1+q+qp)\sum_{u=0}^{\ell-1}(2+(1+q+qp)(k+u+2))\Bigr]\nonumber\\
&=2\Bigr[4+\delta^*(2k+\ell)+\frac{\ell \delta^*}{2}(4+(\ell+2k+3)\delta^*)\Bigr]\nonumber
\end{align}
Therefore, for $i \le 3$, we have
\begin{align}
\sum_{m=0}^{\floor{n/k}}|p_{m,n,k+s+i-1}^{\star\star}\hspace{-0.07cm}-\hspace{-0.07cm}p_{m-1,n,k+s+i-1}^{\star\star}|&\le 2\Bigr[4+\delta^*(3k+s+2)+\frac{(k+s+2) \delta^*}{2}(4+(s+3k+5)\delta^*)\Bigr]\nonumber\\
&=2\Bigr[4\hspace{-0.07cm}+\hspace{-0.07cm}\delta^*\left(5k\hspace{-0.07cm}+\hspace{-0.07cm}6\hspace{-0.07cm}+\hspace{-0.07cm}\frac{3k^2\hspace{-0.07cm}+\hspace{-0.07cm}11k\hspace{-0.07cm}+\hspace{-0.07cm}10}{2}\delta^*\right)\hspace{-0.07cm}+\hspace{-0.07cm}\delta^*\left(3\hspace{-0.07cm}+\hspace{-0.07cm}\frac{4k\hspace{-0.07cm}+\hspace{-0.07cm}7}{2}\delta^*\right)s\hspace{-0.07cm}+\hspace{-0.07cm}\frac{{\delta^*}^2}{2}s^2\Bigr]\nonumber\\
&=2\Bigr[4+\delta_1+\delta_2s+\frac{{\delta^*}^2}{2}s^2\Bigr].\label{xsw}
\end{align}
Now, using \eqref{12} and \eqref{xsw} in the first term of \eqref{bvbv}
\begin{align}
&|1\hspace{-0.07cm}-\hspace{-0.07cm}wb|a(p)^3\sum_{i=1}^{3}|a_i|\sum_{s=0}^{n-k-1}(n-s+1)\sum_{l=0}^{\floor{s/k}}B_s(l)\sum_{m=0}^{\floor{n/k}}|p_{m,n,k+s+i-1}^{\star\star}\hspace{-0.07cm}-\hspace{-0.07cm}p_{m-1,n,k+s+i-1}^{\star\star}|\nonumber\\
&\le 2(2+qp)|1-wb|a(p)^3\sum_{s=0}^{n-k-1}(n-s+1)\frac{(k+1)^s}{(k+2)^{s+1}}\Bigr[4+\delta_1+\delta_2s+\frac{{\delta^*}^2}{2}s^2\Bigr]\nonumber\\
&=2(2+qp)|1-wb|a(p)^3\left[(n-k)(4+\delta_1)+\left((n-2k-2)(k+1)+\frac{(k+1)^{n-k+1}}{(k+2)^{n-k-1}}\right)\delta_2\right.\nonumber\\
&~~~~\left.+\frac{{\delta^*}^2}{2}\left(n(2k+3)-6k^2-16k-10+(2n+2k+5)\frac{(k+1)^{n-k}}{(k+2)^{n-k-1}}\right)\right]\nonumber\\
&=2(2+qp)|1-wb|a(p)^3\left[(n-k)(4+\delta_1)+\delta_2 c_{n,k}^{(1)}+\frac{{\delta^*}^2}{2}c_{n,k}^{(5)}\right].\label{21}
\end{align}
Hence, for $g \in {\cal G}_{Z} \cap {\cal G}_{S_{k_1,k_2}^n}$ and using \eqref{25}, \eqref{24}, \eqref{23}, \eqref{22}, and \eqref{21} in \eqref{bvbv}, we get
\begin{align*}
\left|{\mathbb E}\left[{\cal U} g\left(S_{k_1,k_2}^n\right)\right]\right| & \le \|\Delta g\|\left\{2(2+qp)a(p)^2\left(|1-wb| h_1(n,k,p)a(p)+|wb|h_2(n,k,p)\right)+|\varphi(1-wb)|\right\}.
\end{align*}
This proves result.\qed\\

\noindent
{\bf Proof of Theorem \ref{three:noniidth}}. From \eqref{three:00}, we have
\begin{align}
0={\mathbb E}\left[{\cal A}_{S_{k_1,k_2}^n} g\left(S_{k_1,k_2}^n\right)\right]&=\sum_{m=0}^{\floor{n/k}}\Bigr[we^{U(m+1)-U(m)} g(m+1)-m g(m)-awg(m+1)\nonumber\\
&~~+a(p)\sum_{i=1}^{3}a_i\sum_{s=0}^{d_i}b_i(n-s)\sum_{l=0}^{\floor{s/k}}B_s(l){\mathbb E}\left\{g\left(S_{k_1,k_2}^{n-k-s-i+1}+l+1\right)\Bigr|S_{k_1,k_2}^n=m\right\}\nonumber\\
&~~-w b a(p)\sum_{i=1}^{3}a_i\sum_{s=0}^{d_i}b_i(n\hspace{-0.07cm}-\hspace{-0.07cm}s)\sum_{l=0}^{\floor{s/k}}B_s(l){\mathbb E}\left\{g\left(S_{k_1,k_2}^{n-k-s-i+1}+l+2\right)\Bigr|S_{k_1,k_2}^n\hspace{-0.07cm}=\hspace{-0.07cm}m\right\}\Bigr]p_{m,n}\nonumber\\
&=\sum_{m=0}^{\floor{n/k}}\Bigr[we^{U(m+1)-U(m)} g(m+1)-m g(m)\Bigr]p_{m,n}-aw\sum_{m=0}^{\floor{n/k}}g(m+1)p_{m,n}\nonumber\\
&~~+a(p)\sum_{i=1}^{3}a_i\sum_{s=0}^{d_i}b_i(n-s)\sum_{l=0}^{\floor{s/k}}B_s(l)\sum_{m=0}^{\floor{n/k}}g\left(m+l+1\right)p_{m,n-k-s-i+1}\nonumber\\
&~~-w b a(p)\sum_{i=1}^{3}a_i\sum_{s=0}^{d_i}b_i(n-s)\sum_{l=0}^{\floor{s/k}}B_s(l)\sum_{m=0}^{\floor{n/k}}g\left(m+l+2\right)p_{m,n-k-s-i+1}\nonumber\\
&={\mathbb E}[{\cal A}_Zg(S_{k_1,k_2}^n)]+{\mathbb E}[{\cal U}g(S_{k_1,k_2}^n)].\label{v10}
\end{align}
where
\begin{align*}
{\mathbb E}[{\cal U}g(S_{k_1,k_2}^n)]&=-aw\sum_{m=0}^{\floor{n/k}}g(m+1)p_{m,n}+a(p)\sum_{i=1}^{3}a_i\sum_{s=0}^{d_i}b_i(n-s)\sum_{l=0}^{\floor{s/k}}B_s(l)\sum_{m=0}^{\floor{n/k}}g\left(m+l+1\right)p_{m,n-k-s-i+1}\nonumber\\
&~~~~~-w b a(p)\sum_{i=1}^{3}a_i\sum_{s=0}^{d_i}b_i(n-s)\sum_{l=0}^{\floor{s/k}}B_s(l)\sum_{m=0}^{\floor{n/k}}g\left(m+l+2\right)p_{m,n-k-s-i+1}.\nonumber
\end{align*}
Following the steps from \eqref{three:exp0} to \eqref{three:exp2} and using $\sum_{i=1}^{3}a_i\sum_{s=0}^{d_i}b_i(n-s)\sum_{l=0}^{\floor{s/k}}B_s(l)=q(1+(n-k-1)p)$, we have
\begin{align}
{\mathbb E}\left[{\cal U} g\left(S_{k_1,k_2}^n\right)\right] & = (1-wb)\left(-\frac{wa}{1-wb}+a(p)\sum_{i=1}^{3}a_i\sum_{s=0}^{d_i}b_i(n-s)\sum_{l=0}^{\floor{s/k}}B_s(l)\right)\sum_{m=0}^{\floor{n/k}}g(m+1)p_{m,n}\nonumber\\
&~~~+(1\hspace{-0.07cm}-\hspace{-0.07cm}wb)a(p)^2\sum_{i=1}^{3}a_i\sum_{s=0}^{d_i}b_i(n\hspace{-0.07cm}-\hspace{-0.07cm}s)\sum_{l=0}^{\floor{s/k}}B_s(l)\sum_{m=0}^{\floor{n/k}}g(m\hspace{-0.07cm}+\hspace{-0.07cm}1)(p_{m,n,k+s+i-1}^\star\hspace{-0.07cm}-\hspace{-0.07cm}p_{m-1,n,k+s+i-1}^\star)\nonumber\\
&~~~-wb a(p)\sum_{i=1}^{3}a_i\sum_{s=0}^{d_i}b_i(n-s)\sum_{l=0}^{\floor{s/k}}B_s(l)\sum_{m=0}^{\floor{n/k}}\Delta g(m+l+1)p_{m,n-k-s-i+1}\nonumber\\
&~~~+(1-wb)a(p)\sum_{i=1}^{3}a_i\sum_{s=0}^{d_i}b_i(n-s)\hspace{-0.2cm}\sum_{l=0}^{\floor{s/k}}B_s(l)\hspace{-0.1cm}\sum_{m=0}^{\floor{n/k}}\sum_{j=1}^{l}\Delta g(m+j)p_{m,n-k-s-i+1}\nonumber\allowdisplaybreaks\\
& = -wa\sum_{m=0}^{\floor{n/k}}g(m+1)p_{m,n}+(1-wb)q(1+(n-k-1)p)a(p)\sum_{m=0}^{\floor{n/k}}g(m+1)p_{m,n}\nonumber\\
&~~~-(1\hspace{-0.07cm}-\hspace{-0.07cm}wb)a(p)^2\sum_{i=1}^{3}a_i\sum_{s=0}^{d_i}b_i(n\hspace{-0.07cm}-\hspace{-0.07cm}s)\sum_{l=0}^{\floor{s/k}}B_s(l)\sum_{m=0}^{\floor{n/k}}\Delta g(m\hspace{-0.07cm}+\hspace{-0.07cm}1)p_{m,n,k+s+i-1}^\star\nonumber\\
&~~~-wb a(p)\sum_{i=1}^{3}a_i\sum_{s=0}^{d_i}b_i(n-s)\sum_{l=0}^{\floor{s/k}}B_s(l)\sum_{m=0}^{\floor{n/k}}\Delta g(m+l+1)p_{m,n-k-s-i+1}\nonumber\\
&~~~+(1-wb)a(p)\sum_{i=1}^{3}a_i\sum_{s=0}^{d_i}b_i(n-s)\hspace{-0.2cm}\sum_{l=0}^{\floor{s/k}}B_s(l)\hspace{-0.1cm}\sum_{m=0}^{\floor{n/k}}\sum_{j=1}^{l}\Delta g(m+j)p_{m,n-k-s-i+1},\label{v11}
\end{align}
Substituting \eqref{v11} in \eqref{v10} and using \eqref{panj}, we get
\begin{align}
0={\mathbb E}\left[{\cal A}_{S_{k_1,k_2}^n} g\left(S_{k_1,k_2}^n\right)\right]& = \sum_{m=0}^{\floor{n/k}}\left[wbm g(m+1)-m g(m)+(1-wb)q(1+(n-k-1)p)a(p)g(m+1)\right]p_{m,n}\nonumber\\
&~~~-(1\hspace{-0.07cm}-\hspace{-0.07cm}wb)a(p)^2\sum_{i=1}^{3}a_i\sum_{s=0}^{d_i}b_i(n\hspace{-0.07cm}-\hspace{-0.07cm}s)\sum_{l=0}^{\floor{s/k}}B_s(l)\sum_{m=0}^{\floor{n/k}}\Delta g(m\hspace{-0.07cm}+\hspace{-0.07cm}1)p_{m,n,k+s+i-1}^\star\nonumber\\
&~~~-wb a(p)\sum_{i=1}^{3}a_i\sum_{s=0}^{d_i}b_i(n-s)\sum_{l=0}^{\floor{s/k}}B_s(l)\sum_{m=0}^{\floor{n/k}}\Delta g(m+l+1)p_{m,n-k-s-i+1}\nonumber\\
&~~~+(1-wb)a(p)\sum_{i=1}^{3}a_i\sum_{s=0}^{d_i}b_i(n-s)\hspace{-0.2cm}\sum_{l=0}^{\floor{s/k}}B_s(l)\hspace{-0.1cm}\sum_{m=0}^{\floor{n/k}}\sum_{j=1}^{l}\Delta g(m+j)p_{m,n-k-s-i+1}.\label{three:none}
\end{align}
Consider now
\begin{align}
\sum_{m=0}^{\floor{n/k}}\Delta g(m+1)p_{m,n,k+s+i-1}^\star&=\sum_{m=0}^{\floor{n/k}}\Delta g(m+1)\Bigr[p_{m,n-k}-p_{m,n-2k-s-i+1}{\bf 1}(s\le  n-2k-s-i+1)\nonumber\\
&~~~~+\hspace{-0.07cm}\sum_{u=0}^{k+s+i-2}[qp p_{m,n-k-u-2\hspace{-0.07cm}}\hspace{-0.07cm}-p_{m,n-k-u}{\bf 1}(u\hspace{-0.07cm}=\hspace{-0.07cm}n\hspace{-0.07cm}-\hspace{-0.07cm}k)\hspace{-0.07cm}+\hspace{-0.07cm}qp_{m,n-k-u-1}{\bf 1}(u\hspace{-0.07cm}=\hspace{-0.07cm}n\hspace{-0.07cm}-\hspace{-0.07cm}k-\hspace{-0.07cm}1)]\Bigr]\nonumber\allowdisplaybreaks\\
&={\mathbb E}\Bigr[\Delta g(S_{k_1,k_2}^{n-k}+1)-\Delta g(S_{k_1,k_2}^{n-2k-s-i+1}+1)+\sum_{u=0}^{k+s+i-2}[qp \Delta g(S_{k_1,k_2}^{n-k-u-2}+1)\nonumber\\
&~~~~-\Delta g(S_{k_1,k_2}^{n-k-u}+1){\bf 1}(u=n-k)+q\Delta g(S_{k_1,k_2}^{n-k-u-1}+1){\bf 1}(u=n-k-1)\Bigr]\nonumber\\
&={\mathbb E}\Bigr[{\mathbb E}\Bigr[\Delta g(S_{k_1,k_2}^{n-k}+1)-\Delta g(S_{k_1,k_2}^{n-2k-s-i+1}+1)+\sum_{u=0}^{k+s+i-2}[qp \Delta g(S_{k_1,k_2}^{n-k-u-2}+1)\nonumber\\
&~~~~-\Delta g(S_{k_1,k_2}^{n-k-u}+1){\bf 1}(u\hspace{-0.07cm}=\hspace{-0.07cm}n\hspace{-0.07cm}-\hspace{-0.07cm}k)\hspace{-0.07cm}+\hspace{-0.07cm}q\Delta g(S_{k_1,k_2}^{n-k-u-1}+1){\bf 1}(u\hspace{-0.07cm}=\hspace{-0.07cm}n\hspace{-0.07cm}-\hspace{-0.07cm}k\hspace{-0.07cm}-\hspace{-0.07cm}1)\Bigr|S_{k_1,k_2}^n\Bigr]\Bigr]\nonumber\\
&=\sum_{m=0}^{\floor{n/k}}{\mathbb E}\Bigr[\Delta g(S_{k_1,k_2}^{n-k}+1)\hspace{-0.07cm}-\hspace{-0.07cm}\Delta g(S_{k_1,k_2}^{n-2k-s-i+1}+1)\hspace{-0.07cm}+\hspace{-0.17cm}\sum_{u=0}^{k+s+i-2}[qp \Delta g(S_{k_1,k_2}^{n-k-u-2}+1)\nonumber\\
&~~~~-\hspace{-0.07cm}\Delta g(S_{k_1,k_2}^{n-k-u}\hspace{-0.1cm}+\hspace{-0.07cm}1){\bf 1}(u\hspace{-0.07cm}=\hspace{-0.07cm}n\hspace{-0.07cm}-\hspace{-0.07cm}k)\hspace{-0.07cm}+\hspace{-0.07cm}q\Delta g(S_{k_1,k_2}^{n-k-u-1}\hspace{-0.1cm}+\hspace{-0.07cm}1){\bf 1}(u\hspace{-0.07cm}=\hspace{-0.07cm}n\hspace{-0.07cm}-\hspace{-0.07cm}k\hspace{-0.07cm}-\hspace{-0.07cm}1)\Bigr|S_{k_1,k_2}\hspace{-0.07cm}=\hspace{-0.07cm}m\Bigr]\hspace{-0.05cm}p_{m,n}\nonumber\\
&=\sum_{m=0}^{\floor{n/k}}{\mathbb E}\Bigr[\Delta g(S_{k_1,k_2}^{\star n-k-s-i+1}+1)\Bigr|S_{k_1,k_2}^n=m\Bigr]p_{m,n},\label{b12}
\end{align}
where 
\begin{align}
\left(\Delta g\left(S_{k_1,k_2}^{\star n-l}\hspace{-0.07cm}+\hspace{-0.07cm}1\right)|S_{k_1,k_2}^n\hspace{-0.07cm}=\hspace{-0.07cm}m\right)\hspace{-0.07cm}&=\hspace{-0.07cm}\Bigr[\Bigr(\Delta g(S_{k_1,k_2}^{n-k}\hspace{-0.07cm}+\hspace{-0.07cm}1)\hspace{-0.07cm}-\hspace{-0.07cm}\Delta g(S_{k_1,k_2}^{n-k-l}+1)\hspace{-0.07cm}+\hspace{-0.07cm}\sum_{u=0}^{l-1}[qp \Delta g(S_{k_1,k_2}^{n-k-u-2}\hspace{-0.07cm}+\hspace{-0.07cm}1)\hspace{-0.07cm}\nonumber\\
&~~-\hspace{-0.07cm}\Delta g(S_{k_1,k_2}^{n-k-u}+1){\bf 1}(u\hspace{-0.07cm}=\hspace{-0.07cm}n\hspace{-0.07cm}-\hspace{-0.07cm}k)\hspace{-0.07cm}+\hspace{-0.07cm}q\Delta g(S_{k_1,k_2}^{n-k-u-1}+1){\bf 1}(u\hspace{-0.07cm}=\hspace{-0.07cm}n\hspace{-0.07cm}-\hspace{-0.07cm}k\hspace{-0.07cm}-\hspace{-0.07cm}1)\Bigr)\Bigr|S_{k_1,k_2}=m\Bigr].\label{three:gg}
\end{align}
Following the steps similar from \eqref{b11} to \eqref{three:pp} and using \eqref{b12} in \eqref{three:none}, the Stein operator now becomes
\begin{align}
{\cal A}_{S_{k_1,k_2}^n} g\left(S_{k_1,k_2}^n\right) &= w\left[\frac{(1-wb)}{w}q(1+(n-k-1)p)a(p)+bm\right] g(m+1)-m g(m)\nonumber\\
&~~~-(1\hspace{-0.05cm}-\hspace{-0.05cm}wb)a(p)^2\sum_{i=1}^{3}a_i\sum_{s=0}^{d_i}b_i(n-s)\sum_{l=0}^{\floor{s/k}}B_s(l){\mathbb E}\left\{\Delta g(S_{k_1,k_2}^{\star n-k-s-i+1}\hspace{-0.05cm}+\hspace{-0.05cm}1)|S_{k_1,k_2}^n\hspace{-0.05cm}=\hspace{-0.05cm}m\right\}\nonumber\\
&~~~-wb a(p)\sum_{i=1}^{3}a_i\sum_{s=0}^{d_i}b_i(n-s)\sum_{l=0}^{\floor{s/k}}B_s(l){\mathbb E}\left\{\Delta g(S_{k_1,k_2}^{n-k-s-i+1}+l+1)|S_{k_1,k_2}^n=m\right\}\nonumber\\
&~~~+(1-wb)a(p)\sum_{i=1}^{3}a_i\sum_{s=0}^{d_i}b_i(n-s)\hspace{-0.2cm}\sum_{l=0}^{\floor{s/k}}B_s(l)\hspace{-0.1cm}\sum_{j=1}^{l}{\mathbb E}\left\{\Delta g(S_{k_1,k_2}^{n-k-s-i+1}+j)|S_{k_1,k_2}^n=m\right\}\nonumber\\
&={\cal A}_Zg(m)+\bar{U}_1g(m),\label{three:qq1}
\end{align}
where ${\cal A}_Z$ is a Stein operator for DGM with $a=((1-wb)/w)q(1+(n-k-1)p)a(p)$. Now, from the definition of $\bar{U}_1$ in \eqref{three:qq1} with \eqref{three:inq}, we have
\begin{align}
\|\bar{U}_1g\|&\le|\bar{U}_1 g| \nonumber\\
&\le|1\hspace{-0.05cm}-\hspace{-0.05cm}wb|a(p)^2\sum_{i=1}^{3}|a_i|\sum_{s=0}^{d_i}|b_i(n-s)|\sum_{l=0}^{\floor{s/k}}B_s(l)|{\mathbb E}\left\{\Delta g(S_{k_1,k_2}^{\star n-k-s-i+1}\hspace{-0.05cm}+\hspace{-0.05cm}1)|S_{k_1,k_2}^n\hspace{-0.05cm}=\hspace{-0.05cm}m\right\}|\nonumber\\
&~~~+|wb| a(p)\sum_{i=1}^{3}|a_i|\sum_{s=0}^{d_i}|b_i(n-s)|\sum_{l=0}^{\floor{s/k}}B_s(l)|{\mathbb E}\left\{\Delta g(S_{k_1,k_2}^{n-k-s-i+1}+l+1)|S_{k_1,k_2}^n=m\right\}|\nonumber\\
&~~~+|1-wb|a(p)\sum_{i=1}^{3}|a_i|\sum_{s=0}^{d_i}|b_i(n-s)|\hspace{-0.2cm}\sum_{l=0}^{\floor{s/k}}B_s(l)\hspace{-0.1cm}\sum_{j=1}^{l}|{\mathbb E}\left\{\Delta g(S_{k_1,k_2}^{n-k-s-i+1}+j)|S_{k_1,k_2}^n=m\right\}|\nonumber\allowdisplaybreaks\\
&\le (2+qp)a(p)\|\Delta g\|\left\{|1-wb|a(p)\sum_{s=0}^{n-k-1}(n-s+1)\sum_{l=0}^{\floor{s/k}}B_s(l)\times(2+(1+q+qp)(k+s+2))\right.\nonumber\\
&~~~~~~~~~~~~~~~~~~~~~~~~~~~~~~~~~~~~~~\left.+|wb|\sum_{s=k}^{n-k-1}(n-s+1)\sum_{l=0}^{\floor{s/k}}B_s(l)+|1-wb|\sum_{s=0}^{n-k-1}(n-s+1)\sum_{l=0}^{\floor{s/k}}l B_s(l)\right\}\nonumber\\
&\le (2+qp)a(p)\|\Delta g\|\left\{|1-wb|a(p)\sum_{s=0}^{n-k-1}(n-s+1)\frac{(k+1)^{s}}{(k+2)^{s+1}}\times(2+(1+q+qp)(k+2))\right.\nonumber\\
&~~~~~~~~~~~~~~~~~~~~~~~~~~~~~~~~~~~~~~+|1-wb|a(p)\sum_{s=0}^{n-k-1}(n-s+1)\frac{(k+1)^{s}}{(k+2)^{s+1}}\times(1+q+qp)s\nonumber\\
&~~~~~~~~~~~~~~~~~~~~~~~~~~~~~~~~~~~~~~+|wb|\sum_{s=0}^{n-k-1}(n-s+1)\frac{(k+1)^s}{(k+2)^{s+1}}\nonumber\\
&~~~~~~~~~~~~~~~~~~~~~~~~~~~~~~~~~~~~~~\left.+|1-wb|a(p)\sum_{s=k}^{n-k-1}(n-s+1)((s-k+2)k+2)\frac{(k+1)^{s-k-1}}{(k+2)^{s-k+2}}\right\}\nonumber\allowdisplaybreaks\\
&\le (2+qp)a(p)\|\Delta g\|\left\{|1-wb|a(p)(n-k)\delta+|1-wb|a(p)\delta^* c_{n,k}^{(1)}+|wb|(n-k)+|1-wb|a(p)c_{n,k}^{(2)}\right\}\nonumber\\
&\le \bar{w}_2 \|\Delta g\| \label{three:perbd}
\end{align}
Following the similar steps from \eqref{three:oneexp} to \eqref{three:bigexp} for the last three terms of \eqref{three:none}, we get
\begin{align}
{\mathbb E}\left[{\cal A}_{S_{k_1,k_2}^n} g\left(S_{k_1,k_2}^n\right)\right] &= \sum_{m=0}^{\floor{n/k}}\left[wbm g(m+1)-m g(m)+(1-wb)q(1+(n-k-1)p)a(p)g(m+1)\right]p_{m,n}\nonumber\\
&~~~ \left\{-(1-wb)a(p)^2\sum_{i=1}^{3}a_i\sum_{s=0}^{d_i}b_i(n-s)\sum_{l=0}^{\floor{s/k}}B_s(l)\Bigr[1-{\bf 1}(s \le n-2k-i+1)\right.\nonumber\\
&~~~+\sum_{u=0}^{k+s+i-2}(qp {\bf 1}(u \le n-k-2)-{\bf 1}(u=n-k)+q{\bf 1}(u=n-k-1))\Bigr]\nonumber\\
&~~~-wba(p)\sum_{i=1}^{3}a_i\sum_{s=0}^{d_i}b_i(n-s)\sum_{l=0}^{\floor{s/k}} B_s(l)\nonumber\\
&~~~\left.+(1-wb) a(p)\sum_{i=1}^{3}a_i\sum_{s=0}^{d_i}b_i(n-s)\sum_{l=0}^{\floor{s/k}}lB_s(l)\right\}\sum_{m=0}^{\floor{n/k}}\Delta g(m+1)p_{m,n}\nonumber\allowdisplaybreaks\\
&~~~-\hspace{-0.07cm}(1\hspace{-0.07cm}-\hspace{-0.07cm}wb)a(p)^3\hspace{-0.12cm}\sum_{i=1}^{3}a_i\hspace{-0.1cm}\sum_{s=0}^{d_i}b_i(n\hspace{-0.07cm}-\hspace{-0.07cm}s)\hspace{-0.1cm}\sum_{l=0}^{\floor{s/k}}B_s(l)\hspace{-0.1cm}\sum_{m=0}^{\floor{n/k}}\Delta g(m\hspace{-0.07cm}+\hspace{-0.07cm}1)(p_{m,n,k+s+i-1}^{\star\star}\hspace{-0.07cm}-\hspace{-0.07cm}p_{m-1,n,k+s+i-1}^{\star\star})\nonumber\\
&~~~-\hspace{-0.07cm}wba(p)^2\hspace{-0.1cm}\sum_{i=1}^{3}a_i\hspace{-0.1cm}\sum_{s=0}^{d_i}b_i(n\hspace{-0.07cm}-\hspace{-0.07cm}s)\hspace{-0.1cm}\sum_{l=0}^{\floor{s/k}} B_s(l)\hspace{-0.1cm}\sum_{m=0}^{\floor{n/k}}\Delta g(m\hspace{-0.07cm}+\hspace{-0.07cm}1)(p_{m,n,k+s+i-1}^\star\hspace{-0.07cm}-\hspace{-0.07cm}p_{m-1,n,k+s+i-1}^\star)\nonumber\\
&~~~+\hspace{-0.07cm}(1\hspace{-0.07cm}-\hspace{-0.07cm}wb)a(p)^2\hspace{-0.1cm}\sum_{i=1}^{3}a_i\hspace{-0.1cm}\sum_{s=0}^{d_i}b_i(n\hspace{-0.07cm}-\hspace{-0.07cm}s)\hspace{-0.1cm}\sum_{l=0}^{\floor{s/k}}\hspace{-0.07cm}lB_s(l)\hspace{-0.1cm}\sum_{m=0}^{\floor{n/k}}\hspace{-0.1cm}\Delta g(m\hspace{-0.07cm}+\hspace{-0.07cm}1)(p_{m,n,k+s+i-1}^\star\hspace{-0.07cm}-\hspace{-0.07cm}p_{m-1,n,k+s+i-1}^\star)\nonumber\\
&~~~-wb a(p)\sum_{i=1}^{3}a_i\sum_{s=0}^{d_i}b_i(n-s)\sum_{l=0}^{\floor{s/k}}B_s(l)\sum_{m=0}^{\floor{n/k}}\sum_{j=1}^{l}\Delta^2 g(m+j)p_{m,n-k-s-i+1}\nonumber\\
&~~~+(1-wb)a(p)\sum_{i=1}^{3}a_i\sum_{s=0}^{d_i}b_i(n-s)\sum_{l=0}^{\floor{s/k}}B_s(l)\sum_{m=0}^{\floor{n/k}}\sum_{j=1}^{l}\sum_{v=1}^{j-1}\Delta^2 g(m+v)p_{m,n-k-s-i+1}\nonumber\allowdisplaybreaks\\
&= \sum_{m=0}^{\floor{n/k}}\left[wbm g(m+1)-m g(m)+(1-wb)q(1+(n-k-1)p)a(p)g(m+1)\right]p_{m,n}\nonumber\\
&~~~+\left\{-(1-wb)a(p)^2\sum_{i=1}^{3}a_i\sum_{s=0}^{d_i}b_i(n-s)\sum_{l=0}^{\floor{s/k}}B_s(l)\Bigr[1-{\bf 1}(s \le n-2k-i+1)\right.\nonumber\\
&~~~+\sum_{u=0}^{k+s+i-2}(qp {\bf 1}(u \hspace{-0.07cm}\le \hspace{-0.07cm}n\hspace{-0.07cm}-\hspace{-0.07cm}k\hspace{-0.07cm}-\hspace{-0.07cm}2)\hspace{-0.07cm}-\hspace{-0.07cm}{\bf 1}(u\hspace{-0.07cm}=\hspace{-0.07cm}n\hspace{-0.07cm}-\hspace{-0.07cm}k)+q{\bf 1}(u\hspace{-0.07cm}=\hspace{-0.07cm}n\hspace{-0.07cm}-\hspace{-0.07cm}k\hspace{-0.07cm}-\hspace{-0.07cm}1))\Bigr]-wba(p)\sum_{i=1}^{3}a_i\sum_{s=0}^{d_i}\nonumber\\
&~~~\times b_i(n\hspace{-0.07cm}-\hspace{-0.07cm}s)\sum_{l=0}^{\floor{s/k}} B_s(l)\left.+(1\hspace{-0.07cm}-\hspace{-0.07cm}wb) a(p)\sum_{i=1}^{3}a_i\sum_{s=0}^{d_i}b_i(n\hspace{-0.07cm}-\hspace{-0.07cm}s)\sum_{l=0}^{\floor{s/k}}lB_s(l)\right\}\sum_{m=0}^{\floor{n/k}}\Delta g(m\hspace{-0.07cm}+\hspace{-0.07cm}1)p_{m,n}\nonumber\\
&~~~+(1\hspace{-0.07cm}-\hspace{-0.07cm}wb)a(p)^3\sum_{i=1}^{3}a_i\sum_{s=0}^{d_i}b_i(n\hspace{-0.07cm}-\hspace{-0.07cm}s)\sum_{l=0}^{\floor{s/k}}B_s(l)\sum_{m=0}^{\floor{n/k}}\Delta^2 g(m\hspace{-0.07cm}+\hspace{-0.07cm}1)p_{m,n,k+s+i-1}^{\star\star}\nonumber\\
&~~~+wba(p)^2\sum_{i=1}^{3}a_i\sum_{s=0}^{d_i}b_i(n-s)\sum_{l=0}^{\floor{s/k}} B_s(l)\sum_{m=0}^{\floor{n/k}}\Delta^2 g(m+1)p_{m,n,k+s+i-1}^\star\nonumber\\
&~~~-(1\hspace{-0.07cm}-\hspace{-0.07cm}wb)a(p)^2\sum_{i=1}^{3}a_i\sum_{s=0}^{d_i}b_i(n\hspace{-0.07cm}-\hspace{-0.07cm}s)\sum_{l=0}^{\floor{s/k}}lB_s(l)\hspace{-0.1cm}\sum_{m=0}^{\floor{n/k}}\Delta^2 g(m\hspace{-0.07cm}+\hspace{-0.07cm}1)p_{m,n,k+s+i-1}^\star\nonumber\allowdisplaybreaks\\
&~~~-wb a(p)\sum_{i=1}^{3}a_i\sum_{s=0}^{d_i}b_i(n-s)\sum_{l=0}^{\floor{s/k}}B_s(l)\sum_{m=0}^{\floor{n/k}}\sum_{j=1}^{l}\Delta^2 g(m+j)p_{m,n-k-s-i+1}\nonumber\\
&~~~+(1-wb)a(p)\sum_{i=1}^{3}a_i\sum_{s=0}^{d_i}b_i(n-s)\sum_{l=0}^{\floor{s/k}}B_s(l)\sum_{m=0}^{\floor{n/k}}\sum_{j=1}^{l}\sum_{v=1}^{j-1}\Delta^2 g(m+v)p_{m,n-k-s-i+1}.\label{m13}
\end{align}
Now, consider
\begin{align*}
\sum_{m=0}^{\floor{n/k}}\Delta^2 g(m\hspace{-0.07cm}+\hspace{-0.07cm}1)p_{m,n,k+s+i-1}^{\star\star}&=\sum_{m=0}^{\floor{n/k}}\Delta^2 g(m\hspace{-0.07cm}+\hspace{-0.07cm}1)\Bigr[p_{m,n,k}^\star-p_{m,n,2k+s+i-1}^\star+\sum_{u=0}^{k+s+i-2}[qp p_{m,n,k+u+2}^\star\\
&~~~-p_{m,n,k+u}^\star{\bf 1}(u=n-k)+qp_{m,n,k+u+1}^\star{\bf 1}(u=n-k-1)]\Bigr]
\end{align*}
Using \eqref{b12}, we have
\begin{align}
\sum_{m=0}^{\floor{n/k}}\Delta^2 g(m+1)p_{m,n,k+s+i-1}^{\star\star}&=\sum_{m=0}^{\floor{n/k}}{\mathbb E}\Bigr[\Delta^2 g(S_{k_1,k_2}^{\star n-k}+1)-\Delta^2 g(S_{k_1,k_2}^{\star n-2k-s-i+1}+1)\nonumber\\
&~~~+\sum_{u=0}^{k+s+i-2}[qp \Delta^2 g(S_{k_1,k_2}^{\star n-k-u-2}+1)-\Delta^2 g(S_{k_1,k_2}^{\star n-k-u}+1){\bf 1}(u=n-k)\nonumber\\
&~~~~+q\Delta^2 g(S_{k_1,k_2}^{\star n-k-u-1}+1){\bf 1}(u=n-k-1)]|S_{k_1,k_2}^n=m\Bigr]p_{m,n}\nonumber\\
&=\sum_{m=0}^{\floor{n/k}}{\mathbb E}\Bigr[\Delta^2 g(S_{k_1,k_2}^{**n-k-s-i+1}+1)|S_{k_1,k_2}^n=m\Bigr]p_{m,n},\label{m12}
\end{align}
where
\begin{align*}
\Delta^2 g(S_{k_1,k_2}^{**n-k-s-i+1}+1)&=\Delta^2 g(S_{k_1,k_2}^{\star n-k}+1)-\Delta^2 g(S_{k_1,k_2}^{\star n-2k-s-i+1}+1)+\sum_{u=0}^{k+s+i-2}[qp \Delta^2 g(S_{k_1,k_2}^{\star n-k-u-2}+1)\\
&~~~-\Delta^2 g(S_{k_1,k_2}^{\star n-k-u}+1){\bf 1}(u=n-k)+q\Delta^2 g(S_{k_1,k_2}^{\star n-k-u-1}+1){\bf 1}(u=n-k-1)]
\end{align*}
Next, using \eqref{b12} and \eqref{m12} in \eqref{m13}, we have
\begin{align*}
{\mathbb E}\left[{\cal A}_{S_{k_1,k_2}^n} g\left(S_{k_1,k_2}^n\right)\right]&= \sum_{m=0}^{\floor{n/k}}\left[wbm g(m+1)-m g(m)+(1-wb)q(1+(n-k-1)p)a(p)g(m+1)\right]p_{m,n}\nonumber\\
&~~~+\left\{-(1-wb)a(p)^2\sum_{i=1}^{3}a_i\sum_{s=0}^{d_i}b_i(n-s)\sum_{l=0}^{\floor{s/k}}B_s(l)\Bigr[1-{\bf 1}(s \le n-2k-i+1)\right.\nonumber\\
&~~~+\sum_{u=0}^{k+s+i-2}(qp {\bf 1}(u \hspace{-0.07cm}\le\hspace{-0.07cm} n\hspace{-0.07cm}-\hspace{-0.07cm}k\hspace{-0.07cm}-\hspace{-0.07cm}2)\hspace{-0.07cm}-\hspace{-0.07cm}{\bf 1}(u\hspace{-0.07cm}=\hspace{-0.07cm}n\hspace{-0.07cm}-\hspace{-0.07cm}k)\hspace{-0.07cm}+\hspace{-0.07cm}q{\bf 1}(u\hspace{-0.07cm}=\hspace{-0.07cm}n\hspace{-0.07cm}-\hspace{-0.07cm}k\hspace{-0.07cm}-\hspace{-0.07cm}1))\Bigr]\hspace{-0.07cm}-\hspace{-0.07cm}wba(p)\sum_{i=1}^{3}a_i\sum_{s=0}^{d_i}\nonumber\\
&~~~\times b_i(n\hspace{-0.07cm}-\hspace{-0.07cm}s)\sum_{l=0}^{\floor{s/k}} B_s(l)\left.\hspace{-0.07cm}+\hspace{-0.07cm}(1\hspace{-0.07cm}-\hspace{-0.07cm}wb) a(p)\sum_{i=1}^{3}a_i\sum_{s=0}^{d_i}b_i(n\hspace{-0.07cm}-\hspace{-0.07cm}s)\sum_{l=0}^{\floor{s/k}}lB_s(l)\right\}\sum_{m=0}^{\floor{n/k}}\Delta g(m\hspace{-0.07cm}+\hspace{-0.07cm}1)p_{m,n}\nonumber\\
&~~~+\hspace{-0.07cm}(1\hspace{-0.07cm}-\hspace{-0.07cm}wb)a(p)^3\hspace{-0.1cm}\sum_{i=1}^{3}a_i\hspace{-0.1cm}\sum_{s=0}^{d_i}b_i(n\hspace{-0.07cm}-\hspace{-0.07cm}s)\hspace{-0.1cm}\sum_{l=0}^{\floor{s/k}}\hspace{-0.1cm}B_s(l)\hspace{-0.1cm}\sum_{m=0}^{\floor{n/k}}\hspace{-0.07cm}{\mathbb E}\left\{\Delta^2 g(S_{k_1,k_2}^{**n-k-s-i+1}\hspace{-0.07cm}+\hspace{-0.07cm}1)|S_{k_1,k_2}^n\hspace{-0.07cm}=\hspace{-0.07cm}m\right\}p_{m,n}\nonumber\\
&~~~+\hspace{-0.07cm}wba(p)^2\sum_{i=1}^{3}a_i\sum_{s=0}^{d_i}b_i(n\hspace{-0.07cm}-\hspace{-0.07cm}s)\sum_{l=0}^{\floor{s/k}} B_s(l)\sum_{m=0}^{\floor{n/k}}{\mathbb E}\left\{\Delta^2 g(S_{k_1,k_2}^{*n-k-s-i+1}\hspace{-0.07cm}+\hspace{-0.07cm}1)|S_{k_1,k_2}^n\hspace{-0.07cm}=\hspace{-0.07cm}m\right\}p_{m,n}\nonumber\allowdisplaybreaks\\
&~~~-\hspace{-0.07cm}(1\hspace{-0.07cm}-\hspace{-0.07cm}wb)a(p)^2\hspace{-0.1cm}\sum_{i=1}^{3}a_i\hspace{-0.1cm}\sum_{s=0}^{d_i}b_i(n\hspace{-0.07cm}-\hspace{-0.07cm}s)\hspace{-0.1cm}\sum_{l=0}^{\floor{s/k}}\hspace{-0.1cm}lB_s(l)\hspace{-0.1cm}\sum_{m=0}^{\floor{n/k}}\hspace{-0.07cm}{\mathbb E}\left\{\Delta^2 g(S_{k_1,k_2}^{*n-k-s-i+1}\hspace{-0.07cm}+\hspace{-0.07cm}1)|S_{k_1,k_2}^n\hspace{-0.07cm}=\hspace{-0.07cm}m\right\}p_{m,n}\nonumber\\
&~~~-\hspace{-0.07cm}wb a(p)\sum_{i=1}^{3}a_i\sum_{s=0}^{d_i}b_i(n\hspace{-0.07cm}-\hspace{-0.07cm}s)\hspace{-0.07cm}\sum_{l=0}^{\floor{s/k}}B_s(l)\hspace{-0.07cm}\sum_{m=0}^{\floor{n/k}}\sum_{j=1}^{l}{\mathbb E}\left\{\Delta^2 g(S_{k_1,k_2}^{n-k-s-i+1}\hspace{-0.07cm}+\hspace{-0.07cm}j)|S_{k_1,k_2}^n\hspace{-0.07cm}=\hspace{-0.07cm}m\right\}p_{m,n}\nonumber\\
&~~~+\hspace{-0.07cm}(1\hspace{-0.07cm}-\hspace{-0.07cm}wb)a(p)\hspace{-0.1cm}\sum_{i=1}^{3}\hspace{-0.07cm}a_i\hspace{-0.07cm}\sum_{s=0}^{d_i}b_i(n\hspace{-0.07cm}-\hspace{-0.07cm}s)\hspace{-0.25cm}\sum_{l=0}^{\floor{s/k}}\hspace{-0.13cm}B_s(l)\hspace{-0.25cm}\sum_{m=0}^{\floor{n/k}}\hspace{-0.07cm}\sum_{j=1}^{l}\hspace{-0.07cm}\sum_{v=1}^{j-1}\hspace{-0.07cm}{\mathbb E}\hspace{-0.07cm}\left\{\hspace{-0.07cm}\Delta^2 g(S_{k_1,k_2}^{n-k-s-i+1}\hspace{-0.07cm}+\hspace{-0.07cm}v)|S_{k_1,k_2}^n\hspace{-0.13cm}=\hspace{-0.07cm}m\hspace{-0.07cm}\right\}\hspace{-0.07cm}p_{m,n}\hspace{-0.05cm}.
\end{align*}
Therefore, the Stein operator now becomes
\vspace{-0.31cm}
\begin{align*}
{\cal A}_{S_{k_1,k_2}^n} g\left(m\right) &=wbm g(m+1)-m g(m)+(1-wb)q(1+(n-k-1)p)a(p)g(m+1)\nonumber\\
&~~~+\left\{-(1-wb)a(p)^2\sum_{i=1}^{3}a_i\sum_{s=0}^{d_i}b_i(n-s)\sum_{l=0}^{\floor{s/k}}B_s(l)\Bigr[1-{\bf 1}(s \le n-2k-i+1)\right.\nonumber\\
&~~~+\sum_{u=0}^{k+s+i-2}(qp {\bf 1}(u \le n-k-2)-{\bf 1}(u=n-k)+q{\bf 1}(u=n-k-1))\Bigr]\nonumber\\
&~~~\left.-\hspace{-0.07cm}wba(p)\sum_{i=1}^{3}a_i\sum_{s=0}^{d_i}b_i(n\hspace{-0.07cm}-\hspace{-0.07cm}s)\sum_{l=0}^{\floor{s/k}} B_s(l)\hspace{-0.07cm}+(\hspace{-0.07cm}1\hspace{-0.07cm}-\hspace{-0.07cm}wb) a(p)\sum_{i=1}^{3}a_i\sum_{s=0}^{d_i}b_i(n\hspace{-0.07cm}-\hspace{-0.07cm}s)\sum_{l=0}^{\floor{s/k}}lB_s(l)\right\}\Delta g(m\hspace{-0.07cm}+\hspace{-0.07cm}1)\nonumber\\
&~~~+(1\hspace{-0.07cm}-\hspace{-0.07cm}wb)a(p)^3\sum_{i=1}^{3}a_i\sum_{s=0}^{d_i}b_i(n\hspace{-0.07cm}-\hspace{-0.07cm}s)\sum_{l=0}^{\floor{s/k}}B_s(l){\mathbb E}\left\{\Delta^2 g(S_{k_1,k_2}^{**n-k-s-i+1}\hspace{-0.07cm}+\hspace{-0.07cm}1)|S_{k_1,k_2}^n=m\right\}\nonumber\\
&~~~+wba(p)^2\sum_{i=1}^{3}a_i\sum_{s=0}^{d_i}b_i(n-s)\sum_{l=0}^{\floor{s/k}} B_s(l){\mathbb E}\left\{\Delta^2 g(S_{k_1,k_2}^{*n-k-s-i+1}\hspace{-0.07cm}+\hspace{-0.07cm}1)|S_{k_1,k_2}^n=m\right\}\\
&~~~-(1\hspace{-0.07cm}-\hspace{-0.07cm}wb)a(p)^2\sum_{i=1}^{3}a_i\sum_{s=0}^{d_i}b_i(n\hspace{-0.07cm}-\hspace{-0.07cm}s)\sum_{l=0}^{\floor{s/k}}lB_s(l){\mathbb E}\left\{\Delta^2 g(S_{k_1,k_2}^{*n-k-s-i+1}\hspace{-0.07cm}+\hspace{-0.07cm}1)|S_{k_1,k_2}^n=m\right\}\nonumber\\
&~~~-wb a(p)\sum_{i=1}^{3}a_i\sum_{s=0}^{d_i}b_i(n-s)\sum_{l=0}^{\floor{s/k}}B_s(l)\sum_{j=1}^{l}{\mathbb E}\left\{\Delta^2 g(S_{k_1,k_2}^{n-k-s-i+1}\hspace{-0.07cm}+\hspace{-0.07cm}j)|S_{k_1,k_2}^n=m\right\}\nonumber\\
&~~~+\hspace{-0.07cm}(1\hspace{-0.07cm}-\hspace{-0.07cm}wb)a(p)\sum_{i=1}^{3}a_i\sum_{s=0}^{d_i}b_i(n\hspace{-0.07cm}-\hspace{-0.07cm}s)\sum_{l=0}^{\floor{s/k}}B_s(l)\sum_{j=1}^{l}\sum_{v=1}^{j-1}{\mathbb E}\left\{\Delta^2 g(S_{k_1,k_2}^{n-k-s-i+1}\hspace{-0.07cm}+\hspace{-0.07cm}v)|S_{k_1,k_2}^n=m\right\},
\end{align*}
Now, taking expectation w.r.t. $M_{k_1,k_2}^n$, we have
\begin{align}
{\mathbb E}\hspace{-0.05cm}\left[\hspace{-0.07cm}{\cal A}_{S_{k_1,k_2}^n} g\left({M}_{k_1,k_2}^n\right)\hspace{-0.07cm}\right]\hspace{-0.07cm} &=\hspace{-0.07cm} (1\hspace{-0.07cm}-\hspace{-0.07cm}wb)q(1\hspace{-0.07cm}+\hspace{-0.07cm}(n\hspace{-0.07cm}-\hspace{-0.07cm}k\hspace{-0.07cm}-\hspace{-0.07cm}1)p)a(p){\mathbb E}[g(M_{k_1,k_2}^n\hspace{-0.13cm}+\hspace{-0.09cm}1)]\hspace{-0.07cm}+\hspace{-0.07cm}b{\mathbb E}[M_{k_1,k_2}^n g(M_{k_1,k_2}^n\hspace{-0.13cm}+\hspace{-0.09cm}1)]\hspace{-0.07cm}-\hspace{-0.07cm}{\mathbb E}[M_{k_1,k_2}^n g(M_{k_1,k_2}^n)]\nonumber\\
&~~~+(1-wb)\left\{-a(p)^2\sum_{i=1}^{3}a_i\sum_{s=0}^{d_i}b_i(n-s)\sum_{l=0}^{\floor{s/k}}B_s(l)\Bigr[1-{\bf 1}(s \le n-2k-i+1)\right.\nonumber\\
&~~~+\hspace{-0.07cm}\sum_{u=0}^{k+s+i-2}(qp {\bf 1}(u \hspace{-0.07cm}\le\hspace{-0.07cm} n\hspace{-0.07cm}-\hspace{-0.07cm}k\hspace{-0.07cm}-\hspace{-0.07cm}2)-{\bf 1}(u\hspace{-0.07cm}=\hspace{-0.07cm}n\hspace{-0.07cm}-\hspace{-0.07cm}k)+q{\bf 1}(u\hspace{-0.07cm}=\hspace{-0.07cm}n\hspace{-0.07cm}-\hspace{-0.07cm}k\hspace{-0.07cm}-1\hspace{-0.07cm}))\Bigr]\hspace{-0.07cm}-\hspace{-0.07cm}\frac{wb}{1-wb}a(p)\sum_{i=1}^{3}a_i\nonumber\\
&~~~~~~\times\sum_{s=0}^{d_i}b_i(n-s)\sum_{l=0}^{\floor{s/k}} B_s(l)\left.+a(p)\sum_{i=1}^{3}a_i\sum_{s=k}^{d_i}b_i(n-s)\sum_{l=0}^{\floor{s/k}}lB_s(l)\right\}{\mathbb E}\left\{\Delta g(M_{k_1,k_2}^n+1)\right\}\nonumber\\
&~~~+\hspace{-0.07cm}(1\hspace{-0.07cm}-\hspace{-0.07cm}wb)a(p)^3\hspace{-0.1cm}\sum_{i=1}^{3}a_i\hspace{-0.1cm}\sum_{s=0}^{d_i}b_i(n\hspace{-0.07cm}-\hspace{-0.07cm}s)\hspace{-0.1cm}\sum_{l=0}^{\floor{s/k}}B_s(l){\mathbb E}\left\{{\mathbb E}\left\{\Delta^2 g(S_{k_1,k_2}^{**n-k-s-i+1}\hspace{-0.07cm}+\hspace{-0.07cm}1)|S_{k_1,k_2}^n\hspace{-0.1cm}=\hspace{-0.1cm}M_{k_1,k_2}^n\right\}\hspace{-0.07cm}\right\}\nonumber\\
&~~~+wba(p)^2\sum_{i=1}^{3}a_i\sum_{s=0}^{d_i}b_i(n-s)\sum_{l=0}^{\floor{s/k}} B_s(l) {\mathbb E}\left\{{\mathbb E}\left\{\Delta^2 g(S_{k_1,k_2}^{*n-k-s-i+1}\hspace{-0.07cm}+\hspace{-0.07cm}1)|S_{k_1,k_2}^n=M_{k_1,k_2}^n\right\}\right\}\nonumber\allowdisplaybreaks\\
&~~~-\hspace{-0.07cm}(1\hspace{-0.07cm}-\hspace{-0.07cm}wb)a(p)^2\hspace{-0.1cm}\sum_{i=1}^{3}a_i\hspace{-0.1cm}\sum_{s=0}^{d_i}b_i(n\hspace{-0.07cm}-\hspace{-0.07cm}s)\hspace{-0.1cm}\sum_{l=0}^{\floor{s/k}}lB_s(l){\mathbb E}\left\{{\mathbb E}\left\{\Delta^2 g(S_{k_1,k_2}^{*n-k-s-i+1}\hspace{-0.07cm}+\hspace{-0.07cm}1)|S_{k_1,k_2}^n\hspace{-0.07cm}=\hspace{-0.07cm}M_{k_1,k_2}^n\right\}\hspace{-0.1cm}\right\}\nonumber\\
&~~~-\hspace{-0.07cm}wb a(p)\sum_{i=1}^{3}a_i\hspace{-0.1cm}\sum_{s=0}^{d_i}b_i(n\hspace{-0.07cm}-\hspace{-0.07cm}s)\hspace{-0.1cm}\sum_{l=0}^{\floor{s/k}}B_s(l)\hspace{-0.1cm}\sum_{j=1}^{l}{\mathbb E}\left\{{\mathbb E}\left\{\Delta^2 g(S_{k_1,k_2}^{n-k-s-i+1}\hspace{-0.07cm}+\hspace{-0.07cm}j)|S_{k_1,k_2}^n\hspace{-0.07cm}=\hspace{-0.07cm}M_{k_1,k_2}^n\right\}\right\}\nonumber\\
&~~~+\hspace{-0.08cm}(1\hspace{-0.07cm}-\hspace{-0.07cm}wb)a(p)\hspace{-0.1cm}\sum_{i=1}^{3}a_i\hspace{-0.1cm}\sum_{s=0}^{d_i}b_i(n\hspace{-0.07cm}-\hspace{-0.07cm}s)\hspace{-0.17cm}\sum_{l=0}^{\floor{s/k}}\hspace{-0.17cm}B_s(l)\hspace{-0.07cm}\sum_{j=1}^{l}\hspace{-0.07cm}\sum_{v=1}^{j-1}\hspace{-0.07cm}{\mathbb E}\hspace{-0.07cm}\left\{\hspace{-0.07cm}{\mathbb E}\hspace{-0.07cm}\left\{\Delta^2\hspace{-0.05cm} g(S_{k_1,k_2}^{n\hspace{-0.07cm}-\hspace{-0.01cm}k\hspace{-0.01cm}-\hspace{-0.01cm}s\hspace{-0.01cm}-\hspace{-0.01cm}i\hspace{-0.01cm}+\hspace{-0.01cm}1}\hspace{-0.1cm}+\hspace{-0.07cm}v)|S_{k_1,k_2}^n\hspace{-0.07cm}=\hspace{-0.07cm}M_{k_1,k_2}^n\right\}\hspace{-0.07cm}\right\}\hspace{-0.07cm}.\label{three:bigrexp}
\end{align}
Next, define
\vspace{-0.31cm}
\begin{align*}
M_i&=M_{k_1,k_2}^n-{\mathbb I}_i \quad{\rm and}\quad N_i=M_{k_1,k_2}^n-\sum_{|s-i|\le k+1} {\mathbb I}_s.
\end{align*}
Then $N_i$ and ${\mathbb I}_i$ are independent. It is given that 
\vspace{-0.31cm}
$${\mathbb E}\big(S_{k_1,k_2}^n\big)={\mathbb E}\big(M_{k_1,k_2}^n\big)\implies q(1+(n-k-1)p)a(p)=\sum_{i=1}^{n-k}{\mathbb E}({\mathbb I}_i).$$
Consider the first three terms of \eqref{three:bigrexp}, we have
\vspace{-0.31cm}
\begin{align}
&(1-wb)q(1+(n-k-1)p)a(p){\mathbb E}[g(M_{k_1,k_2}^n+1)]+wb{\mathbb E}[M_{k_1,k_2}^n g(M_{k_1,k_2}^n+1)]-{\mathbb E}[M_{k_1,k_2}^n g(M_{k_1,k_2}^n)]\nonumber\\
&=(1-wb)\sum_{i=1}^{n-k}{\mathbb E}({\mathbb I}_i){\mathbb E}[g(M_{k_1,k_2}^n+1)]-(1-wb){\mathbb E}[M_{k_1,k_2}^n g(M_{k_1,k_2}^n)]+wb{\mathbb E}[M_{k_1,k_2}^n \Delta g(M_{k_1,k_2}^n)]\nonumber\\
&=(1-wb)\sum_{i=1}^{n-k}\{{\mathbb E}({\mathbb I}_i){\mathbb E}[g(M_{k_1,k_2}^n+1)]-{\mathbb E}[{\mathbb I}_i g(M_i+1)]\}+wb{\mathbb E}[M_{k_1,k_2}^n \Delta g(M_{k_1,k_2}^n)]\nonumber\\
&=(1-wb)\sum_{i=1}^{n-k}{\mathbb E}({\mathbb I}_i) {\mathbb E}[g(M_{k_1,k_2}^n+1)-g(N_i+1)]+(1-wb)\sum_{i=1}^{n-k}{\mathbb E}[{\mathbb I}_i(g(N_i+1)-g(M_i+1))]\nonumber\\
&~~~+wb{\mathbb E}[M_{k_1,k_2}^n \Delta g(M_{k_1,k_2}^n)]\nonumber\allowdisplaybreaks\\
&=(1-wb)\sum_{i=1}^{n-k}{\mathbb E}({\mathbb I}_i)\sum_{|j-i|\le k+1} {\mathbb E}\left( {\mathbb I}_j \Delta g\left(N_i+\sum_{s=i-k-1}^{j-1}{\mathbb I}_s+1\right)\right)\nonumber\\
&~~~-(1-wb)\sum_{i=1}^{n-k}\sum_{\substack{|j-i|\le k+1\\j\neq i}}{\mathbb E}\left({\mathbb I}_i{\mathbb I}_j\Delta g\left(N_i+\sum_{\substack{s=i-k-1\\s\neq i}}^{j-1}{\mathbb I}_s+1\right)\right)+wb{\mathbb E}[M_{k_1,k_2}^n \Delta g(M_{k_1,k_2}^n)].\label{three:sim}
\end{align}
Observe the terms involving ${\mathbb E}(\Delta g(M_{k_1,k_2}^n+1))$ in \eqref{three:bigrexp}, we have
\vspace{-0.31cm}
\begin{align}
a(p)\sum_{i=1}^{3}a_i\sum_{s=0}^{d_i}b_i(n-s)\sum_{l=0}^{\floor{s/k}} B_s(l)=q(1+(n-k-1)p)a(p)=\sum_{i=1}^{n-k}{\mathbb E}({\mathbb I}_i) \label{three:mmm}
\end{align}
and, using \eqref{three:mmm} with $\tau=\mathrm{Var}\big(M_{k_1,k_2}^n\big)-\mathrm{Var}\big(S_{k_1,k_2}^n\big)$, it can be verified that
\vspace{-0.31cm}
\begin{align}
&-a(p)^2\sum_{i=1}^{3}a_i\sum_{s=0}^{d_i}b_i(n-s)\sum_{l=0}^{\floor{s/k}}B_s(l)\Bigr[1-{\bf 1}(s \le n-2k-i+1)+\sum_{u=0}^{k+s+i-2}(qp {\bf 1}(u \le n-k-2)\nonumber\\
&~~~-{\bf 1}(u=n-k)+q{\bf 1}(u=n-k-1))\Bigr]+a(p)\sum_{i=1}^{3}a_i\sum_{s=0}^{d_i}b_i(n-s)\sum_{l=0}^{\floor{s/k}}lB_s(l)\nonumber\\
&=-s_{n,k}=-\sum_{i=1}^{n-k}{\mathbb E}({\mathbb I}_i)\sum_{|j-i|\le k+1}{\mathbb E}({\mathbb I}_j)+\sum_{i=1}^{n-k}\sum_{\substack{|j-i|\le k+1\\j\neq i}}{\mathbb E}({\mathbb I}_i {\mathbb I}_{j})-\tau\label{nbnb}
\end{align}
Now, using \eqref{three:sim}, \eqref{three:mmm} and \eqref{nbnb} in \eqref{three:bigrexp}, we get
\vspace{-0.12cm}
\begin{align}
{\mathbb E}\hspace{-0.05cm}\left[\hspace{-0.07cm}{\cal A}_{S_{k_1,k_2}^n} g\left({M}_{k_1,k_2}^n\right)\hspace{-0.07cm}\right]\hspace{-0.07cm} &=-(1\hspace{-0.05cm}-\hspace{-0.05cm}wb)\hspace{-0.07cm}\sum_{i=1}^{n-k}{\mathbb E}({\mathbb I}_i)\hspace{-0.17cm}\sum_{|j-i|\le k+1}\left[{\mathbb E}({\mathbb I}_j){\mathbb E}(\Delta g(M_{k_1,k_2}^n\hspace{-0.05cm}+\hspace{-0.05cm}1))\hspace{-0.05cm}-\hspace{-0.05cm}{\mathbb E}\left({\mathbb I}_j\Delta g \left(N_i\hspace{-0.05cm}+\hspace{-0.05cm}\sum_{s=i-k-1}^{j-1}{\mathbb I}_s\hspace{-0.05cm}+\hspace{-0.05cm}1\right)\right)\right]\nonumber\\
&~~~+(1\hspace{-0.05cm}-\hspace{-0.05cm}wb)\sum_{i=1}^{n-k}\sum_{\substack{|j-i|\le k+1\\j\neq i}}\left[{\mathbb E}({\mathbb I}_i{\mathbb I}_j){\mathbb E}(\Delta g(M_{k_1,k_2}^n\hspace{-0.05cm}+\hspace{-0.05cm}1))\hspace{-0.05cm}-\hspace{-0.05cm}{\mathbb E}\left({\mathbb I}_i{\mathbb I}_j\Delta g \left(N_i\hspace{-0.05cm}+\hspace{-0.05cm}\sum_{\substack{s=i-k-1\\s\neq i}}^{j-1}{\mathbb I}_s\hspace{-0.05cm}+\hspace{-0.05cm}1\right)\right)\right]\nonumber\allowdisplaybreaks\\
&~~~-wb \sum_{i=1}^{n-k}\left[{\mathbb E}({\mathbb I}_i){\mathbb E}(\Delta g(M_{k_1,k_2}^n+1))-{\mathbb E}\left({\mathbb I}_i\Delta g \left(M_i+1\right)\right)\right]-\tau (1-wb){\mathbb E}(\Delta g(M_{k_1,k_2}^n+1))\nonumber\\
&~~~+\hspace{-0.07cm}(1\hspace{-0.07cm}-\hspace{-0.07cm}wb)a(p)^3\hspace{-0.1cm}\sum_{i=1}^{3}a_i\hspace{-0.1cm}\sum_{s=0}^{d_i}b_i(n\hspace{-0.07cm}-\hspace{-0.07cm}s)\hspace{-0.1cm}\sum_{l=0}^{\floor{s/k}}B_s(l){\mathbb E}\left\{{\mathbb E}\left\{\Delta^2 g(S_{k_1,k_2}^{**n-k-s-i+1}\hspace{-0.07cm}+\hspace{-0.07cm}1)|S_{k_1,k_2}^n\hspace{-0.1cm}=\hspace{-0.1cm}M_{k_1,k_2}^n\right\}\hspace{-0.07cm}\right\}\nonumber\\
&~~~+wba(p)^2\sum_{i=1}^{3}a_i\sum_{s=0}^{d_i}b_i(n-s)\sum_{l=0}^{\floor{s/k}} B_s(l) {\mathbb E}\left\{{\mathbb E}\left\{\Delta^2 g(S_{k_1,k_2}^{*n-k-s-i+1}\hspace{-0.07cm}+\hspace{-0.07cm}1)|S_{k_1,k_2}^n=M_{k_1,k_2}^n\right\}\right\}\nonumber\\
&~~~-\hspace{-0.07cm}(1\hspace{-0.07cm}-\hspace{-0.07cm}wb)a(p)^2\hspace{-0.1cm}\sum_{i=1}^{3}a_i\hspace{-0.1cm}\sum_{s=0}^{d_i}b_i(n\hspace{-0.07cm}-\hspace{-0.07cm}s)\hspace{-0.1cm}\sum_{l=0}^{\floor{s/k}}lB_s(l){\mathbb E}\left\{{\mathbb E}\left\{\Delta^2 g(S_{k_1,k_2}^{*n-k-s-i+1}\hspace{-0.07cm}+\hspace{-0.07cm}1)|S_{k_1,k_2}^n\hspace{-0.07cm}=\hspace{-0.07cm}M_{k_1,k_2}^n\right\}\hspace{-0.1cm}\right\}\nonumber\\
&~~~-\hspace{-0.07cm}wb a(p)\sum_{i=1}^{3}a_i\hspace{-0.1cm}\sum_{s=0}^{d_i}b_i(n\hspace{-0.07cm}-\hspace{-0.07cm}s)\hspace{-0.1cm}\sum_{l=0}^{\floor{s/k}}B_s(l)\hspace{-0.1cm}\sum_{j=1}^{l}{\mathbb E}\left\{{\mathbb E}\left\{\Delta^2 g(S_{k_1,k_2}^{n-k-s-i+1}\hspace{-0.07cm}+\hspace{-0.07cm}j)|S_{k_1,k_2}^n\hspace{-0.07cm}=\hspace{-0.07cm}M_{k_1,k_2}^n\right\}\right\}\nonumber\\
&~~~+\hspace{-0.08cm}(1\hspace{-0.07cm}-\hspace{-0.07cm}wb)a(p)\hspace{-0.1cm}\sum_{i=1}^{3}a_i\hspace{-0.1cm}\sum_{s=0}^{d_i}b_i(n\hspace{-0.07cm}-\hspace{-0.07cm}s)\hspace{-0.17cm}\sum_{l=0}^{\floor{s/k}}\hspace{-0.17cm}B_s(l)\hspace{-0.07cm}\sum_{j=1}^{l}\hspace{-0.07cm}\sum_{v=1}^{j-1}\hspace{-0.07cm}{\mathbb E}\hspace{-0.07cm}\left\{\hspace{-0.07cm}{\mathbb E}\hspace{-0.07cm}\left\{\Delta^2\hspace{-0.05cm} g(S_{k_1,k_2}^{n\hspace{-0.07cm}-\hspace{-0.01cm}k\hspace{-0.01cm}-\hspace{-0.01cm}s\hspace{-0.01cm}-\hspace{-0.01cm}i\hspace{-0.01cm}+\hspace{-0.01cm}1}\hspace{-0.1cm}+\hspace{-0.07cm}v)|S_{k_1,k_2}^n\hspace{-0.07cm}=\hspace{-0.07cm}M_{k_1,k_2}^n\right\}\hspace{-0.07cm}\right\}\hspace{-0.07cm}.\nonumber
\end{align}
Therefore,
\begin{align}
\left|{\mathbb E}\hspace{-0.05cm}\left[{\cal A}_{S_{k_1,k_2}^n} g\left({M}_{k_1,k_2}^n\right)\right]\right|\hspace{-0.07cm} &=|1\hspace{-0.05cm}-\hspace{-0.05cm}wb|\sum_{i=1}^{n-k}{\mathbb E}({\mathbb I}_i)\sum_{|j-i|\le k+1}\left|{\mathbb E}({\mathbb I}_j){\mathbb E}(\Delta g(M_{k_1,k_2}^n\hspace{-0.05cm}+\hspace{-0.05cm}1))\hspace{-0.05cm}-\hspace{-0.05cm}{\mathbb E}\left({\mathbb I}_j\Delta g \left(N_i\hspace{-0.05cm}+\hspace{-0.2cm}\sum_{s=i-k-1}^{j-1}{\mathbb I}_s\hspace{-0.05cm}+\hspace{-0.05cm}1\right)\right)\right|\nonumber\\
&~~~+|1\hspace{-0.05cm}-\hspace{-0.05cm}wb|\sum_{i=1}^{n-k}\sum_{\substack{|j-i|\le k+1\\j\neq i}}\left|{\mathbb E}({\mathbb I}_i{\mathbb I}_j){\mathbb E}(\Delta g(M_{k_1,k_2}^n\hspace{-0.05cm}+\hspace{-0.05cm}1))\hspace{-0.05cm}-\hspace{-0.05cm}{\mathbb E}\left({\mathbb I}_i{\mathbb I}_j\Delta g \left(N_i\hspace{-0.05cm}+\hspace{-0.2cm}\sum_{\substack{s=i-k-1\\s\neq i}}^{j-1}{\mathbb I}_s\hspace{-0.05cm}+\hspace{-0.05cm}1\right)\right)\right|\nonumber\allowdisplaybreaks\\
&~~~+|wb| \sum_{i=1}^{n-k}\left|{\mathbb E}({\mathbb I}_i){\mathbb E}(\Delta g(M_{k_1,k_2}^n\hspace{-0.07cm}+\hspace{-0.07cm}1))\hspace{-0.07cm}-\hspace{-0.07cm}{\mathbb E}\left({\mathbb I}_i\Delta g \left(M_i\hspace{-0.07cm}+\hspace{-0.07cm}1\right)\right)\right|\hspace{-0.07cm}+\hspace{-0.07cm}|\tau (1\hspace{-0.07cm}-\hspace{-0.07cm}wb)||{\mathbb E}(\Delta g(M_{k_1,k_2}^n\hspace{-0.07cm}+\hspace{-0.07cm}1))|\nonumber\\
&+\hspace{-0.07cm}\|\Delta g\|\hspace{-0.12cm}\left\{|1\hspace{-0.07cm}-\hspace{-0.07cm}wb|a(p)^3\sum_{i=1}^{3}|a_i|\hspace{-0.1cm}\sum_{s=0}^{n-k-1}\hspace{-0.17cm}(n\hspace{-0.07cm}-\hspace{-0.07cm}s\hspace{-0.07cm}+\hspace{-0.07cm}1)\hspace{-0.1cm}\sum_{l=0}^{\floor{s/k}}\hspace{-0.12cm}B_s(l)\hspace{-0.1cm}\sum_{m=0}^{\floor{n/k}}|p_{m,n,k+s+i-1}^{\star\star}\hspace{-0.07cm}-\hspace{-0.07cm}p_{m-1,n,k+s+i-1}^{\star\star}|\right.\nonumber\\
&~~~+|wb|a(p)^2\sum_{i=1}^{3}|a_i|\sum_{s=0}^{n-k-1}(n-s+1)\sum_{l=0}^{\floor{s/k}} B_s(l)\sum_{m=0}^{\floor{n/k}}|p_{m,n,k+s+i-1}^\star-p_{m-1,n,k+s+i-1}^\star|\nonumber\\
&~~~+|1\hspace{-0.07cm}-\hspace{-0.07cm}wb|a(p)^2\sum_{i=1}^{3}|a_i|\sum_{s=k}^{n-k-1}(n-s+1)\sum_{l=0}^{\floor{s/k}}lB_s(l)\hspace{-0.1cm}\sum_{m=0}^{\floor{n/k}}|p_{m,n,k+s+i-1}^\star\hspace{-0.07cm}-\hspace{-0.07cm}p_{m-1,n,k+s+i-1}^\star|\nonumber\\
&~~~+\hspace{-0.08cm}|wb| a(p)\sum_{i=1}^{3}|a_i|\hspace{-0.08cm}\sum_{s=k}^{n-k-1}(n-s+1)\hspace{-0.08cm}\sum_{l=0}^{\floor{s/k}}lB_s(l)\hspace{-0.08cm}\sum_{m=0}^{\floor{n/k}}\hspace{-0.08cm}|p_{m-1,n-k-s-i+1}-p_{m,n,n-k-s-i+1}|\nonumber\\
&~~~\left.+\hspace{-0.07cm}|1\hspace{-0.07cm}-\hspace{-0.07cm}wb|a(p)\hspace{-0.1cm}\sum_{i=1}^{3}\hspace{-0.07cm}|a_i|\hspace{-0.25cm}\sum_{s=2k}^{n-k-1}\hspace{-0.23cm}(n\hspace{-0.07cm}-\hspace{-0.07cm}s\hspace{-0.07cm}+\hspace{-0.07cm}1)\hspace{-0.2cm}\sum_{l=0}^{\floor{s/k}}\hspace{-0.1cm}\frac{l(l\hspace{-0.07cm}-\hspace{-0.07cm}1)}{2}B_s(l)\hspace{-0.14cm}\sum_{m=0}^{\floor{n/k}}\hspace{-0.17cm}|p_{m-1,n-k-s-i+1}\hspace{-0.07cm}-\hspace{-0.07cm}p_{m,n,n-k-s-i+1}|\hspace{-0.1cm}\right\}\hspace{-0.12cm}.\label{three:bigrexp2}
\end{align}
Define
$$S_i=M_{k_1,k_2}^n-\sum_{|j-i|\le 2k+2}{\mathbb I}_j.$$
Then, $S_i$ and $\{{\mathbb I}_j:|j-i|\le k+1\}$ are independent. Now, consider
\begin{align}
&\left|{\mathbb E}({\mathbb I}_j){\mathbb E}(\Delta g(M_{k_1,k_2}^n+1))-{\mathbb E}\left({\mathbb I}_j\Delta g \left(N_i+\sum_{s=i-k-1}^{j-1}{\mathbb I}_s+1\right)\right)\right|\nonumber\\
&= \left|{\mathbb E}({\mathbb I}_j){\mathbb E}\{\Delta g(M_{k_1,k_2}^n+1)-\Delta g(S_i+1)\}+{\mathbb E}\left\{{\mathbb I}_j\left(\Delta g(S_i+1)-\Delta g \left(N_i+\sum_{s=i-k-1}^{j-1}{\mathbb I}_s+1\right)\right)\right\}\right|\nonumber\\
&\le {\mathbb E}({\mathbb I}_j)\left|{\mathbb E}\{\Delta g(M_{k_1,k_2}^n+1)-\Delta g(S_i+1)\}\right|+\left|{\mathbb E}\left\{{\mathbb I}_j\left(\Delta g(S_i+1)-\Delta g \left(N_i+\sum_{s=i-k-1}^{j-1}{\mathbb I}_s+1\right)\right)\right\}\right|\allowdisplaybreaks\nonumber\\
&= \left|\sum_{u=i-2k-2}^{j-1}{\mathbb E}\left({\mathbb I}_j {\mathbb I}_u \Delta^2 g\left(S_i+\sum_{s=i-2k-2}^{u-1}{\mathbb I}_s+\sum_{s=i+k+2}^{i+2k+2}{\mathbb I}_s+1\right)\right)+\sum_{v=i+k+2}^{i+2k+2}{\mathbb E}\left({\mathbb I}_j {\mathbb I}_v \Delta^2 g\left(S_i+\hspace{-0.2cm}\sum_{s=i+k+2}^{v-1}{\mathbb I}_s+1\right)\right)\right|\nonumber\\
&~~~+{\mathbb E}({\mathbb I}_j)\left|\sum_{|u-i|\le 2k+2}{\mathbb E}\left({\mathbb I}_u\Delta^2 g\left(S_i+\sum_{s=i-2k-2}^{u-1}{\mathbb I}_s+1\right)\right)\right|\nonumber\\
&\le 2 \|\Delta g\|\left\{\left(\sum_{u=i-2k-2}^{j-1}+\sum_{u=i+k+2}^{i+2k+2}\right){\mathbb E}({\mathbb I}_j {\mathbb I}_u)+{\mathbb E}({\mathbb I}_j)\sum_{|u-i|\le 2k+2}{\mathbb E}({\mathbb I}_u)\right\}.\label{1}
\end{align}
Similarly,
\begin{align}
&\left|{\mathbb E}({\mathbb I}_i{\mathbb I}_j){\mathbb E}(\Delta g(M_{k_1,k_2}^n\hspace{-0.15cm}+\hspace{-0.1cm}1))\hspace{-0.1cm}-\hspace{-0.1cm}{\mathbb E}\hspace{-0.1cm}\left(\hspace{-0.1cm}{\mathbb I}_i{\mathbb I}_j\Delta g \hspace{-0.1cm}\left(\hspace{-0.1cm}N_i\hspace{-0.1cm}+\hspace{-0.3cm}\sum_{\substack{s=i-k-1\\s\neq i}}^{j-1}\hspace{-0.3cm}{\mathbb I}_s\hspace{-0.1cm}+\hspace{-0.1cm}1\hspace{-0.1cm}\right)\hspace{-0.1cm}\right)\hspace{-0.1cm}\right|\hspace{-0.1cm}\le\hspace{-0.1cm} 2\|\Delta g\|\hspace{-0.1cm}\left\{\hspace{-0.1cm}\left(\hspace{-0.1cm}\sum_{\substack{u=i-2k-2\\u\neq i}}^{j-1}\hspace{-0.15cm}+\hspace{-0.15cm}\sum_{u=i+k+2}^{i+2k+2}\hspace{-0.1cm}\right)\hspace{-0.1cm}{\mathbb E}({\mathbb I}_i{\mathbb I}_j{\mathbb I}_u)\hspace{-0.1cm}+\hspace{-0.1cm}{\mathbb E}({\mathbb I}_i{\mathbb I}_j)\hspace{-0.6cm}\sum_{|u-i|\le 2k+2}\hspace{-0.34cm}{\mathbb E}({\mathbb I}_u)\hspace{-0.1cm}\right\}\label{2}
\end{align}
and
\begin{align}
\left|{\mathbb E}({\mathbb I}_i){\mathbb E}(\Delta g(M_{k_1,k_2}^n+1))-{\mathbb E}\left({\mathbb I}_i\Delta g \left(M_i+1\right)\right)\right|&\le 2 \|\Delta g\|\left\{\sum_{\substack{|j-i|\le k+1\\j\neq i}}{\mathbb E}({\mathbb I}_i{\mathbb I}_j)+{\mathbb E}({\mathbb I}_i)\sum_{|j-i|\le k+1}{\mathbb E}({\mathbb I}_j)\right\}.\label{3}
\end{align}
Now, using \eqref{25}, \eqref{24}, \eqref{23}, \eqref{22}, \eqref{21}, \eqref{1}, \eqref{2} and \eqref{3} in \eqref{three:bigrexp2}, we get 
\begin{align*}
|{\mathbb E}\hspace{-0.05cm}\Bigr[{\cal A}_{S_{k_1,k_2}^n} g\Bigr({M}_{k_1,k_2}^n\Bigr)\Bigr]|&\le 2\|\Delta g\|\Bigg\{|1-wb|\Bigg[\sum_{i=1}^{n-k}\sum_{|j-i|\le k+1}\Bigg(\Bigr(\sum_{u=i-2k-2}^{j-1}+\sum_{u=i+k+2}^{i+2k+2}\Bigr)[{\mathbb E}({\mathbb I}_i){\mathbb E}({\mathbb I}_j {\mathbb I}_u)+{\mathbb E}({\mathbb I}_i{\mathbb I}_j{\mathbb I}_u)]\\
&~~~+[{\mathbb E}({\mathbb I}_i){\mathbb E}({\mathbb I}_j)+{\mathbb E}({\mathbb I}_i{\mathbb I}_j)]\sum_{|u-i|\le 2k+2}{\mathbb E}({\mathbb I}_u)\Bigg)+\frac{|\tau|}{2}+(2+qp)h_1(n,k,p)a(p)^3\Bigg]\\
&~~~+|wb|\Bigg[\sum_{i=1}^{n-k}\Bigg(\sum_{|j-i|\le k+1}{\mathbb E}({\mathbb I}_i{\mathbb I}_j)+{\mathbb E}({\mathbb I}_i)\sum_{|j-i|\le k+1}{\mathbb E}({\mathbb I}_j)\Bigg)+(2+qp)h_2(n,k,p)a(p)^2\Bigg] \Bigg\}.
\end{align*}
Hence,
\begin{equation}
\left|{\mathbb E}\left[{\cal A}_{S_{k_1,k_2}^n} g\left({M}_{k_1,k_2}^n\right)\right]\right|\le \bar{\varepsilon}^*\|\Delta g\|.\label{three:45}
\end{equation}
Using \eqref{three:perbd} and \eqref{three:45} with Lemma \ref{three:perturbation}, we get the required result.\qed

\small
\singlespacing


\begin{thebibliography}{99}
\bibitem{AKI} {Aki, S.} (1997). On sooner and later problems between success and failure runs. {\em Advances in Combinatorial Methods and Applications to Probability and Statistics} (Ed., N. Balakrishnan), 385-400, Borkh\"{a}user, Boston.
\bibitem{AKH} {Aki, S., Kuboki, H. and Hirano, K.} (1984). On discrete distributions of order $k$. {\em Ann. Inst. Statist. Math.} {\bf 36}, 431–440.
\bibitem{ABK} {Antzoulakos, D. L., Bersimis, S., Koutras, M. V.} (2003). Waiting times associated with the sum of success run lengths. In: Lindqvist, B., Doksum, K. (Eds.), {\em Mathematical and Statistical Methods in Reliability}. World Scientific, Singapore, 141–157.
\bibitem{AC} {Antzoulakos, D. L. and Chadjiconstantinidis, S.} (2001). Distributions of numbers of success runs of fixed length in Markov dependent trials. {\em Ann. Inst. Statist. Math.} {\bf 53}, 599–619.
\bibitem{BC} {Balakrishnan, N. and Chan, P. S.} (1999). Two-stage start-up demonstration testing, In {\em Statistical and Probabilistic Models in Reliability} (Eds., D. C. Ionescu and N. Limnios), Bikh\"{a}user, Bosten, 251-263.
\bibitem{BK} {Balakrishnan, N. and Koutras, M. V.} (2002). Runs and Scans with Applications, {\em John Wiley $\&$ Sons}, New York.
\bibitem{Ba} {Barbour, A. D.} (1990). Stein’s method for diffusion approximations. {\em Probab. Theory and Related Fields}, {\bf 84}, 297-322.
\bibitem{BCX} {Barbour, A. D., \v{C}ekanavi\v{c}ius, V. and Xia, A.} (2007). On Stein’s method and perturbations. {\em ALEA}, {\bf 3}, 31-53.
\bibitem{BAC} {Barbour, A. D. and Chen, L. H. Y.} (2014). Stein's (magic) method. {Preprint}:arXiv:1411.1179.
\bibitem{BCL} {Barbour, A. D., Chen, L. H. Y. and Loh, W. L.} (1992). Compound Poisson approximation for nonnegative random variables via Stein’s method. {\em Ann. Prob.}, {\bf 20}, 1843-1866.
\bibitem{BHJ} {Barbour, A. D., Holst, L. and Janson, S.} (1992). Poisson Approximation. Oxford University Press, Oxford.
\bibitem{CV} {\v{C}ekanavi\v{c}ius, V.} (2016). Approximation Methods in Probability Theory. {\em Springer, Cham,} Universitext.
\bibitem{CR} {\v{C}ekanavi\v{c}ius, V. and Roos, B.} (2004). Two-parametric compound binomial approximations. {\em Lith. Math. J.}, {\bf 44}, 354-373.
\bibitem{CGS} {Chen, L. H. Y., Goldstein, L. and Shao, Q.-M.} (2011). Normal Approximation by Stein’s Method. {\em Springer}, Heidelberg.
\bibitem{DAP} {Dafnis, S. D., Antzoulakos, D. L. and Philippou, A. N.} (2010). Distribution related to $(k_1,k_2)$ events. {\em J. Stat. Plan. Inference}, {\bf 140}, 1691-1700.
\bibitem{DZ} {Diaconis, P. and Zabell, S.} (1991). Closed form summation for classical distributions. variations on a theme of de Moivre. {\em Statist. Sci.}, {\bf 6}, 284–302.
\bibitem{ER} {Eichelsbacher, P. and Reinert, G.} (2008). Stein’s method for discrete Gibbs measures. {\em Ann. Appl. Probab.}, {\bf 18}, 1588-1618.
\bibitem{FK} {Fu, J. C. and Koutras, M. V.} (1994). Distribution theory of runs: a Markov chain approach, {\em J. Amer. Statist. Assoc.}, {\bf 89}, 1050-1058.
\bibitem{G} {G\"{o}tze, F.} (1991). On the rate of convergence in the multivariate CLT. {\em Ann. Prob.}. {\bf 19}, 724–739.
\bibitem{HT} {Huang, W. T. and Tsai, C. S.} (1991). On a modified binomial distribution of order $k$. {\em Statist. Prob. Lett.}. {\bf 11}, 125-131.
\bibitem{koutras} {Koutras, M. V.} (1997). Waiting time distributions associated with runs of fixed length in two-state Markov chains, {\em Ann. Inst. Statist. Math.}, {\bf 49}, 123-139.
\bibitem{KU1} {Kumar, A. N. and Upadhye, N. S.} (2016). On perturbations of Stein operator. Comm. Statist. Theory Methods., {\bf 46}, 9284-9302. 
\bibitem{KU2} {Kumar, A. N. and Upadhye, N. S.} (2016). Pseudo-binomial Approximation to $(k_1, k_2) $-runs. {\em arXiv preprint  	arXiv:1609.07847}.
\bibitem{LRS} {Ley, C., Reinert G. and Swan, Y.} (2014). Stein’s method for comparison of univariate distributions. {\em Probab. Surv.}, {\bf 14}, 1-52.
\bibitem{MPP} {Makri, F. S., Philippou, A. N. and Psillakis, Z. M.} (2007). Shortest and longest length of success runs in binary sequences. {\em J. Statist. Plann. Inference} 137.7, 2226–2239.
\bibitem{NP} {Nourdin, I. and Peccati, G.} (2012). Normal Approximations with Malliavin Calculus. From Stein’s Method to Universality. {\em Cambridge Tracts in Mathematics No. 192.}
\bibitem{Panj95} Panjer, H. H. and Wang, Sh. (1995). Computational aspects of Sundt's generalized class. \emph{ASTIN Bull.}, {\bf 25}, 5-17.
\bibitem{PGP} {Philippou, A. N., Georghiou, C. and Philippou, G. N.} (1983). A generalized distribution and some of its properties. {\em Statist. Prob. Lett.} {\bf 1}, 171–175.
\bibitem{PM} {Philippou, A. N. and Makri, A.} (1986). Success, runs and longest runs. {\em Statist. Prob. Lett.} {\bf 4}, 211–215.
\bibitem{R1} {Reinert, G.} (2005). Three general approaches to Stein's method. {An introduction to Stein's method}, {\em Lect. Notes Ser. Inst. Math. Sci. Natl. Univ. Singap.}, {\bf 4}, Singapore University Press, Singapore, 183-221.
\bibitem{stein} {Stein, C.} (1972). A bound for the error in the normal approximation to the distribution of a sum of dependent random variables, {\em Proc. Sixth Berkeley Symp. Math. Statist. Prob.}, {\bf 2}, 583-602. University of California Press, Berkeley.
\bibitem{UCV} {Upadhye, N. S., \v{C}ekanavi\v{c}ius V. and Vellaisamy, P.} (2014). On Stein operators for discrete approximations. {\em Bernoulli}, {\bf 23}, 2828-2859.

\bibitem{V} {Vellaisamy, P.} (2004). Poisson approximation for $(k_1, k_2)$-events via the Stein–Chen method. {\em J. Appl. Probab.} {\bf 41}, 1081–1092.
\bibitem{VUC} {Vellaisamy, P., Upadhye, N. S., and \v{C}ekanavi\v{c}ius, V.} (2013). On Negative Binomial Approximation. {\em Theory Probab. Appl.}, {\bf 57}(1), 97-109.
\end{thebibliography}
\end{document}